\documentclass[amssymb,amscd,amsfonts,refcheck,11pt,graphicx,verbatim,righttag]{amsart}

\usepackage{amsmath,amscd, amsthm, amsfonts, amssymb,mathrsfs,multirow,graphicx,color}

\usepackage{amsfonts}
\usepackage{amssymb}
\usepackage[OT2,OT1]{fontenc}
\usepackage{caption}
\usepackage{longtable}

 \allowdisplaybreaks[4]
\theoremstyle{plain}

\newcommand{\N}{{\mathbb N}}
\newcommand{\G}{{\mathcal G}}
\newcommand{\LL}{{\mathcal L}}
\newcommand{\GG}{{\mathscr G}}
\newcommand{\Z}{{\mathbb Z}}
\newcommand{\R}{{\mathbb R}}
\newcommand{\RR}{{\mathcal R}}

\newcommand{\ZZ}{{\mathcal Z}}
\newcommand{\E}{{\mathcal E}}
\newcommand{\B}{{\mathcal B}}
\newcommand{\BB}{{\mathscr B}}
\newcommand{\I}{{\mathcal I}}

\newcommand{\A}{{\mathcal A}}

\allowdisplaybreaks

\newtheorem{thm}{Theorem}[section]
\newtheorem{cor}[thm]{Corollary}
\newtheorem{prop}[thm]{Proposition}
\newtheorem{rem}[thm]{Remark}
\newtheorem{lem}[thm]{Lemma}
\newtheorem{defi}[thm]{Definition}
\newtheorem{exam}[thm]{Example}

\def \proof{\bigbreak\noindent{\it Proof.~~}}

\begin{document}
\title{The spectrum of period-doubling Hamiltonian}

\author{Qinghui LIU}
\address[Q.H. LIU]{
Department of  Computer Science,
Beijing Institute of Technology,
Beijing 100081, P. R. China.}
\email{qhliu@bit.edu.cn}
\author{Yanhui QU}
\address[Y.H. QU]{Department of  Mathematical Science, Tsinghua University, Beijing 100084, P. R. China.}
\email{yhqu@tsinghua.edu.cn}
\author{Xiao Yao}
\address[X. Yao]{School of Mathematical Sciences and LPMC, Nankai University, Tianjin, 300071,  P. R. China.}
\email{yaoxiao@nankai.edu.cn}

\begin{abstract}
In this paper, we show the following: the Hausdorff dimension of the spectrum of period-doubling Hamiltonian is bigger than $\log \alpha/\log 4$, where $\alpha$ is the Golden number; there exists a dense uncountable subset
of the spectrum such that for each energy in this set, the related trace orbit is unbounded,
which is in contrast with  a recent  result of Carvalho (Nonlinearity 33, 2020);  we  give a complete characterization for the structure of gaps and the gap labelling of the spectrum.
All of these results are consequences of an  intrinsic coding of the spectrum we construct in this paper.

\end{abstract}

\maketitle


\section{Introduction}


\subsection{Background and motivation}\

Since the pioneer works \cite{KKT,OPRSS}, the discrete Schr\"odinger operators
with substitutional  potentials have been extensively studied. In these two works, a crucial
concept--trace map, were introduced, which became a standard and powerful tool for studying the spectral
properties of this class of operators. Using trace map technique, Bovier and Ghez \cite{BG} showed that if
 the potential is generated by primitive substitution, then, under certain technical condition, the spectrum of the related operator
 is of zero Lebesgue measure. Liu et. al. \cite{LTWW} removed this technical condition. This is
  the most general result in this direction. Lenz \cite{L02} also proved the same result by using cocycle dynamical
   technique, which is quite different from the previous two papers.

   Three operators are heavily studied: the Fibonacci, Thue-Morse and period-doubling Hamiltonians.
   The potentials of these  Hamiltonians are generated by Fibonacci, Thue-Morse and period-doubling  substitutions, respectively. We note that Fibonacci substitution is invertible, while the Thue-Morse and period-doubling  substitutions are not. The differences in substitution rules for these  models have great impact on the corresponding dynamical systems induced by the trace polynomials. For the Fibonacci case, the induced dynamical system is a diffeomorphism on a cubic surface, which has a  close relation with
    H\'enon map on $\mathbb{C}^2$. There has been much progress for understanding the dynamics of H\'enon map in the past three decades. As a  resonance, many powerful tools have been developed in the study of the spectrum of the Fibonacci Hamiltonian, see \cite{Can,DG,DG2,DG3,DGY}. While,  for the Thue-Morse and period-doubling models, the corresponding two dimensional dyanimical systems are completely different: they are not diffeomorphisms at all, see for example \cite{BBG,LQY} and the present work.  This makes them  extremely hard to study, both from  real and complex analytic point of views. To the best of our knowledge,   there is no even a basic picture for the dynamics of this kind of maps  in the community of complex dynamical system of higher dimension.

   Fibonacci Hamiltonian
   was introduced in \cite{KKT,OPRSS}, which soon became the most popular model for quasicrystal.
   Let $\{\hat h_n: n\ge1\}$ be the trace polynomials related to Fibonacci Hamiltonian, Casdagli \cite{Cas} defined
    the pesudospectrum of the operator:
    \begin{equation}\label{pseudo-spe}
    \hat B_\infty:=\left\{E\in \R: \{\hat h_n(E): n\ge1\} \text{ is bounded} \right\}
    \end{equation}
    and showed that $\hat B_\infty$ is of zero Lebesgue measure. S\"ut\H{o} \cite{Suto} showed that $\hat B_\infty$
     coincides with the spectrum $\hat \sigma_\lambda$ of Fibonacci Hamiltonian. Dynamically, $\hat B_\infty$ can be redefined as follows(see for example \cite{Cas}). Let
     $$\hat f(x,y,z):=(y,z,yz-x)$$
      be the {\it trace map} of Fibonacci Hamiltonian, then $\hat f:S_\lambda\to S_\lambda$ is a diffeomorphism, where $\lambda>0$ is the coupling constant and
      $$
      S_\lambda:=\{(x,y,z)\in \R^3: x^2+y^2+z^2-xyz=4+\lambda^2\}.
      $$
      Moreover, for any $n\ge0$ and $E\in \R$, one has
     $$\hat f(\hat h_n(E),\hat h_{n+1}(E),\hat h_{n+2}(E))=(\hat h_{n+1}(E),\hat h_{n+2}(E),\hat h_{n+3}(E)).$$
     For any $E$, write the initial condition as $\hat \ell(E):=(2,E,E-\lambda)$. Then
     \begin{equation}\label{B-infty}
    \hat B_\infty=\left\{E\in \R: \{\hat f^n(\hat \ell(E)):n\ge1\} \text{ is bounded} \right\}.
    \end{equation}

     These two papers suggest a way of
     characterizing  the spectrum via boundedness of trace orbits, compare \eqref{pseudo-spe}.
     Along this direction, Damanik \cite{Da} showed
     that if the potential is generated by an invertible primitive
     substitution over two letters, then the spectrum of the related operator still coincides with
     $\hat B_\infty$. In \cite{Cas}, Casdagli  essentially  coded the spectrum by a kind of
      subshift of finite type when the coupling  $\lambda>16$. Raymond \cite{Ray} developed this approach to give coding for
      Sturmian Hamiltonian, which  has Sturmian sequence as its potential. We note that Sturmian
      sequence is not necessarily generated by single substitution. Through this coding, he could compute
      the integrated density of states(IDS) of all the energies in an explicit way when $\lambda>4$. His approach also gave
       a proof for gap opening and gap labelling of Sturmian Hamiltonian.
      The  Hausdorff dimension  of $\hat \sigma_\lambda$ has been well understood,  see \cite{Ray, JL,Can,DEGT,DG,DG2,DGY}. In particular  the following property is shown in \cite{DEGT}:
$$
\lim_{\lambda\to\infty} (\dim_H \hat \sigma_\lambda)\log \lambda=\log (1+\sqrt{2}).
$$
This implies that $\dim_H \hat \sigma_\lambda\to 0$ with the speed $1/\log \lambda$ when $\lambda\to\infty.$
In the remarkable paper \cite{DGY}, almost all interesting  spectral properties of Fibonacci Hamiltonian are established  for all the coupling $\lambda>0$.

Next  we discuss Thue-Morse Hamiltonian(TMH).
This operator has potential generated by Thue-Morse substitution, which is primitive but non invertible, as we have mentioned above.
TMH was studied in many works at 1980's via trace maps, see for example \cite{AAKMP,AP1,AP,B,Lu}.
In particular,  the gap labelling and gap opening properties of the operator were studied in \cite{AP1,B,Lu}.
However,  more detailed  structure of the spectrum (such as a coding) is not known compare with the Fibonacci case.
This is reflected by the fact that the computation of the dimension of the spectrum is very difficult.
Recently it was shown in \cite{LQ} that for any $\lambda>0$, the spectrum $\tilde \sigma_\lambda$ of TMH satisfies
 $$
 \dim_H \tilde \sigma_\lambda\ge \frac{\log 2}{140\log 2.1}=0.00667\cdots
 $$
  Also, a numerical estimation for the box dimension of the spectrum was provided in \cite{PEF}.
Indicated by the results of \cite{Cas,Suto},
one may  expect that the spectrum also coincides with $\tilde B_\infty$, where
\begin{equation}\label{pseudo-spe-tm}
    \tilde B_\infty:=\left\{E\in \R: \{\tilde h_n(E): n\ge1\} \text{ is bounded} \right\},
    \end{equation}
    where $\{\tilde h_n:n\ge0\}$ are the trace polynomials related to TMH.
 This problem was studied in \cite{LQY}.
There,  a Thue-Morse trace map
$\tilde f:\mathbb{R}^2\to \mathbb{R}^2$ is defined as
\begin{equation*}\label{def-2-dim-map}
\tilde f(x,y):=(x^2(y-2)+2,x^2y^2(y-2)+2).
\end{equation*}
Note that $\tilde f$ is not a diffeomorphism. Indeed it is  not even surjective and is  generically 4-to-1.
For any $E\in\mathbb{R}$ and any $n\ge1$,
\begin{equation*}
\tilde f(\tilde h_n(E),\tilde h_{n+1}(E))=(\tilde h_{n+2}(E),\tilde h_{n+3}(E)).
\end{equation*}
Write the initial condition as $\tilde\ell(E):=(E^2-\lambda^2-2,(E^2-\lambda^2)^2-4E^2+2)$, then
$$
\tilde B_\infty=\{E\in \R: \{\tilde f^n(\tilde \ell(E)):n\ge1\} \text{ is bounded}\}.
$$
In \cite{LQY}, it was shown that  $\tilde{\mathscr E}_\infty:=\tilde \sigma_\lambda\setminus  \tilde B_\infty$ is dense in $\tilde \sigma_\lambda$ and uncountable. That means there are many energies in the spectrum with unbounded trace orbit. We remark that, by a bit extra work, one indeed  can show that $\tilde B_\infty\subset \tilde \sigma_\lambda$, see for example Theorem \ref{main-infty-energy}.

Now we come to period-doubling Hamiltonian(PDH), which is the object of this paper.
The potential of PDH is generated by periodic-doubling substitution.
Luck \cite{Lu} first studied the trace maps of PDH and presented several conjectures about the labelling and opening of the gaps.
Bellissard, Bovier and Ghez \cite{BBG} confirmed these conjectures.
Moreover, they showed that the spectrum is of zero Lebesgue measure. Like the Thue-Morse case, further information  on the structure of the spectrum is missing.
The first motivation of this paper is the following:
try to find out some coding map similar to that of Fibonacci Hamiltonian,
then use this coding to study the dimension, the IDS and  gap labeling of the spectrum in a concrete manner,
as Raymond did for Sturmian Hamiltonian(\cite{Ray}).

Most recently, Carvalho \cite{Car} showed that for PDH,
the spectrum also coincides with the set of energies related to bounded trace orbit. However, our experience in \cite{LQY}
seems suggest that the opposite result should hold, just like the Thue-Morse case. Our second motivation is to understand the problem of existence of unbounded trace orbit completely.

\subsection{Main results}\

Now we set up the setting  and state the main results of the paper.

\subsubsection{Basic definitions}\label{sec-def}\

Let $\eta$ be the {\it period-doubling substitution}: $\eta(a)=ab$, $\eta(b)=aa$. It is seen that
$\eta^{2n}(a)$ is both a prefix  and suffix of $\eta^{2(n+1)}(a)$.
Define  a two-sided sequence $\xi$ as
$$\xi:=\lim_{n\rightarrow\infty}\eta^{2n}(a)|\eta^{2n}(a)
=\cdots \xi(-2)\xi(-1)|\xi(0)\xi(1)\cdots.$$
Define the {\it period-doubling potential} $V_\xi=(V_\xi(n))_{n\in\Z}$ by,
$V_\xi(n)=1$ if $\xi(n)=a$ and $V_\xi(n)=-1$ if $\xi(n)=b$ for $n\in\mathbb{Z}$.
Take $\lambda\in \R$ and $\lambda\ne0$.
Let $H_{\lambda V_\xi}$ be the  discrete Schr\"odinger operator acting
 on $\ell^2(\Z)$ with potential $\lambda V_\xi$, i.e.,
for any $n\in\mathbb{Z}$,
$$(H_{\lambda V_\xi}\psi)_n=\psi_{n+1}+\psi_{n-1}+\lambda V_\xi(n)\psi_n, \ \ \ \forall \psi\in\ell^2(\Z).$$
We call $H_{\lambda V_\xi}$ the {\it period-doubling Hamiltonian}(PDH).
Denote by  $\sigma(H_{\lambda V_\xi})$ the  spectrum of $H_{\lambda V_\xi}$.
From now on, we only consider PDH, so we will write
\begin{equation}\label{H-lambda}
H_{\lambda}:=H_{\lambda V_\xi} \ \ \text{ and } \ \ \sigma_\lambda:=\sigma(H_{\lambda V_\xi}).
\end{equation}
 It is a general fact that
$\sigma_{-\lambda}=-\sigma_{\lambda}$. So in the following, without loss of generality,
we always assume $\lambda>0$.

Given $E\in \R,$ define an anti-homomorphism $\tau:{\rm FG}(a,b)\to {\rm SL}(2,\R)$ as
$$
\tau(a):= \left[
\begin{array}{cc}
E-\lambda&-1\\
1&0
\end{array}
\right]\ \ \ \text{ and }\ \ \ \tau(b):= \left[
\begin{array}{cc}
E+\lambda&-1\\
1&0
\end{array}
\right]
$$
and $\tau(a_1\cdots a_n):=\tau(a_n)\cdots\tau(a_1).$ Define
\begin{equation}\label{trace-poly-n}
h_n(E):={\rm tr }(\tau(\eta^n(a))).
\end{equation}
$h_n$ is called  the $n$-th {\it  trace polynomial} related to PDH.
It is seen that ${\rm deg}(h_n)=2^n.$ For each $E\in \R$, we define
 \begin{equation}\label{O-E}
 \mathscr O(E):=\{h_n(E):n\ge0\}
 \end{equation}
and call $\mathscr O(E)$ the {\it trace orbit} of $E$. Define
\begin{equation}\label{B-infty}
B_\infty:=\{E\in \R: \mathscr O(E) \text{ is bounded}\}.
\end{equation}
 We say $E\in \sigma_\lambda$ is of {\it $\infty$-type} if $\mathscr O(E)$ is unbounded. Define
\begin{equation}\label{E-lambda}
\mathscr E_\infty:=\{ E\in \sigma_\lambda: E \text{ is of } \infty\text{-type}\}.
\end{equation}

\subsubsection{Lower bound for the Hausdorff dimension of the spectrum}\

At first we have the following estimate for the dimension of the spectrum:

\begin{thm}\label{main-low-dim}
Let $\alpha={(1+\sqrt{5})}/{2}$.  Then  for any $\lambda>0$,
$$
\dim_H\sigma_{\lambda}\ge\frac{\log\alpha}{\log 4}=0.34712\cdots
$$
\end{thm}

\begin{rem}
{\rm
1) In \cite{PEF}, a numerical estimation for $\dim_H\sigma_\lambda$ is given. The numerical  lower bound of $\dim_H\sigma_\lambda$ is around $0.75$.

 2) Indeed, $\log\alpha$ is the entropy of certain subshift of finite type  and $\log 4$ can be viewed as an upper bound of the Lyapunov exponents, see Remark \ref{rem-dim} for explanation. So the lower bound is kind of Young's dimension formula.

3) As we have mentioned, in \cite{LQ} we obtained a lower bound $0.00667\cdots$ for the Hausdorff dimension of the spectrum of TMH. Let us make a comparison on the related  methods.
 In \cite{LQ}, we showed  that, after a suitable renormalization, the trace polynomials $\{\tilde h_n:n\ge1\}$ are exponentially close to the model family  $\{2\cos 2^nE: n\ge1\}$. Using this closeness, we can construct a nested covering structure such that the ratios of lengths  between the son intervals and the father intervals are bounded from below. By this, we can show that the dimension of the limit set has a universal  lower bound (which is  $0.00667\cdots$, very small and far from optimal). Also by the construction, the limit set is a subset of the spectrum, so we obtain the dimension lower bound for the spectrum. It is desirable to adapt the method to the PDH case. It is true that one can still find a model family $\{2\cos 2^n\sqrt{E}: n\ge1\}$ such that  after a suitable renormalization, the trace polynomials $\{h_n:n\ge1\}$ of PDH are  close to them. However it is much harder to show the exponential convergence. It is even harder to construct a good Cantor subset of the spectrum if it was possible. Here we take a totally different approach. Due to the explicit coding of the spectrum (see Theorem \ref{main-coding-spectrum}), we can construct a nested covering structure with nice separation property. This structure determines a limit  set, which is a subset of the spectrum.
  Every energy in this subset has a bounded  trace orbit and the bound is 2.  This bound together with the recurrence relation of the trace polynomials  imply that the lengths of the bands in level $n$ are bigger than $C4^{-n}$. Now it is not hard to show that the dimension of the  limit set has a good lower bound.

  4) We also note that it is doable to apply the method in the present paper to the TMH to give a much better lower bound for the Hausdorff dimension of the spectrum.
}
\end{rem}

\subsubsection{Existence of unbounded trace orbits}\

\begin{thm}\label{main-infty-energy}
$\ZZ\subset B_\infty\subset \sigma_\lambda$ and  
$\mathscr E_\infty=\sigma_\lambda\setminus B_\infty$. Both $B_\infty$ and  $\mathscr E_\infty$ are dense in $\sigma_\lambda$ and uncountable.
\end{thm}

\begin{rem}
{\rm
 1) As we have mentioned, Carvalho \cite{Car} claimed that $\mathscr E_\infty$ is empty. Here we show the contrary.  We will see soon that  $\infty$-type energies are closely related to the gap edges of the spectrum.
 
 2) In \cite{LQY}, for TMH, we showed that $\tilde{\mathscr E}_\infty$ is dense in $\tilde \sigma_\lambda$ and uncountable. But we did not show that $\tilde B_\infty$ is a subset of $\tilde \sigma_\lambda$. Here, we show that $\{B_\infty, \mathscr E_\infty\}$ form a partition of $\sigma_\lambda.$

 3) Similar to the TMH case, we will construct a period-doubling trace map
$ f:\mathbb{R}^2\to \mathbb{R}^2$ (see \eqref{f}) such that
for any $E\in\mathbb{R}$ and any $n\ge0$,
\begin{equation*}
 f(h_n(E),h_{n+1}(E))=( h_{n+2}(E), h_{n+3}(E)).
\end{equation*}
If we write the initial condition as $\ell(E):=(E-\lambda,E^2-\lambda^2-2)$, then
$$
B_\infty=\{E\in \R: \{f^n( \ell(E)):n\ge1\} \text{ is bounded}\}.
$$
We will see that the dynamical properties of $f$ play an essential role for showing the injectivity of the coding map. }
\end{rem}

Both results are  consequences of an intrinsic coding of the spectrum. So we turn to this coding map.

\subsubsection{Coding of the spectrum}\label{sec-main-coding}\

We will code the spectrum $\sigma_\lambda$ via a subshift of finite type. At first we define such a subshift. We will  explain the idea behind these definitions in Sec. \ref{idea-graph}.

\begin{figure}[htbp]
\begin{center}
\includegraphics[width=0.60\textwidth]{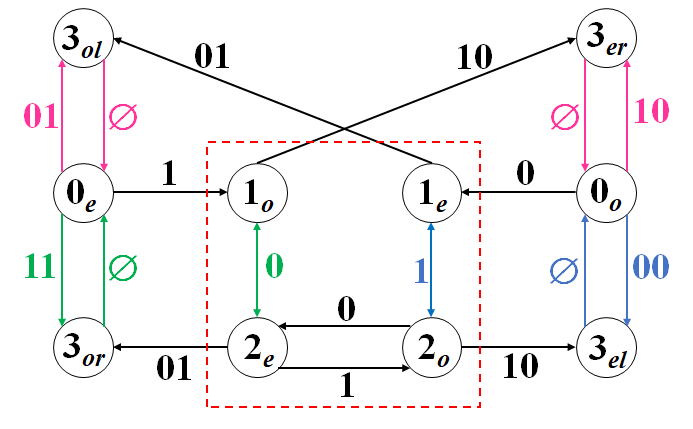}
\caption{The directed graph $\mathbb G$ related to the subshift}\label{fig-graph}
\end{center}
\end{figure}

Define an {\it alphabet} $\A$ as
\begin{equation}\label{alphabet}
\A:=\{0_e,0_o,1_{e},1_{o}, 2_e,2_o,3_{el},3_{er},3_{ol},3_{or}\}.
\end{equation}
Define  the {\it admissible rules} as
\begin{equation}\label{adm}
\begin{split}
 0_e\to \{3_{ol}, 1_o,3_{or}\};\ \  0_o\to\{3_{el},1_e,3_{er}\};\ \
 1_e\to\{3_{ol},2_o\}; \ \ 1_o\to\{2_e,3_{er}\};\\
2_e\to\{1_o,3_{or},2_o\};\ \   \ 2_o\to\{2_e,3_{el},1_e\};
\ \ \ 3_{el}, 3_{er}\to 0_o;\ \ \ \  3_{ol}, 3_{or}\to 0_e.
\end{split}
\end{equation}
See Figure \ref{fig-graph} for the related directed graph $\mathbb G$ (we will explain the extra information such as the colors and the labels of the edges, the dashed red square etc.  later).
Define the {\it adjacency matrix} $A=[a_{\alpha \beta}]$ of $\mathbb G$ as
$$
\begin{cases}
a_{\alpha \beta}=1, \  \text{ if }\ \alpha\to\beta; \\
a_{\alpha \beta}=0, \  \text{ otherwise. }
\end{cases}
$$
It is seen that $A$ is irreducible but not aperiodic, indeed it has period 2. Consequently, $\mathbb G$ is connected. Let $\Omega_A$ be the related subshift of finite type:
\begin{equation*}\label{subshift}
\Omega_A:=\{\omega\in \A^\N: a_{\omega_{j}\omega_{j+1}}=1, j\ge0\}.
\end{equation*}

Define a metric $d$ on $\Omega_A$ as
$
d(\omega,\hat\omega):=2^{-|\omega\wedge \hat\omega|}.
$
It is standard to check that $d$ is indeed a metric and $(\Omega_A,d)$ is compact.
Define the symbolic space
\begin{equation}\label{def-omega-infty}
\Omega_\infty:=\{\omega\in \Omega_A: \omega_0\in \{3_{el},0_e\}\}.
\end{equation}
 $\Omega_\infty$ is a union of two cylinders:  $\Omega_\infty=[3_{el}]\cup [0_e]$, so it is compact.

 Later we will define a total order $\preceq$ on $\Omega_\infty$, which is essentially induced by the ordinary order $\le $ on $\R$. ( We will explain this in Sec. \ref{idea-order-gap} and Sec. \ref{sec-symbolic}).

\begin{thm}\label{main-coding-spectrum}

There exists an order-preserving
homeomorphism $\pi: (\Omega_\infty,\preceq)\to (\sigma_\lambda,\le )$,
where $\le $ is the standard order on $\R$.
\end{thm}

We call $\pi$ the {\it coding map} of $\sigma_\lambda$.

\begin{rem}
{\rm
In some sense,  we construct a ``Markov partition" of the spectrum $\sigma_\lambda$ through the coding given in Theorem \ref{main-coding-spectrum}. The idea of constructing combinatoric models (Markov partition, Young tower, Yoccoz puzzle, inducing scheme, etc.) play a central role for understanding the hyperbolicity of the dynamical system, which can date back to Bowen.
For the spectrum of discrete Schr\"odinger operator, the idea of constructing  certain Markov partition can date back to Casdagli for Fibonacci Hamiltonian(\cite{Cas}) and  Raymond  for Sturmian Hamiltonian (\cite{Ray}). Some recent progress  have been made  by Cantat(\cite{Can}) and Damanik,  Gorodetski,  Yessen (\cite{DG,DG2,DG3,DGY}), which heavily rely on the dynamics of $\hat f$ on the surface $S_\lambda$.

  To some extent, all the constructions of the combinatoric models relies on the bounded distortion (complex bound) arguments.
For the PDH model, the existence of the energies with unbounded trace orbits is the main obstruction for the bounded distortion phenomenon, which does not happen in the Fibonacci model. This makes the construction of combinatoric model extremely hard, we need be very careful in discussion of various combinatorics of the zeros, bands, ends of the bands and related issues.

}
\end{rem}

\begin{rem}
{\rm
Let us make a comparison with \cite{Ray} on the methods of constructing codings.
 For both classes of Hamiltonians, let $\sigma_n$ be the $n$-th periodic approximation of the spectrum $\sigma(H_\lambda)$, one can show that $\{\sigma_n\cup\sigma_{n+1}:n\ge1\}$ is a nested sequence of coverings of the spectrum. By choosing  an optimal covering in each level, it is not hard to construct a symbolic space $\Omega$ and a surjective coding $\pi: \Omega\to \sigma(H_\lambda)$. Usually it is much harder to get a bijective coding. To achieve this, one need some disjoint condition for bands in $\sigma_n$ and $\sigma_{n+1}$. That is why in \cite{Ray}, the technical assumption $\lambda>4$ is needed. We remark that when $\lambda\downarrow 0$, the overlapping  of $\sigma_n$ and $\sigma_{n+1}$ is unavoidable. This  presents an obstruction for a bijective coding. In our case, many  bands in the $n$-th optimal covering do overlap with each other(at least for small $\lambda$).  One can see Figure \ref{2-type-gap} (b) for some intuition of this situation. A big challenge is to show that as $n$ grows, this kind of overlapping  will disappear. By a detailed analysis on the band configurations of $\sigma_n\cup\sigma_{n+1}\cup\sigma_{n+2}$, finally we succeed to show that all the overlapping disappear and   gaps of $\sigma(H_\lambda)$ show up.  Surprisingly, we find that the  edges of this kind of gaps are related to unbounded trace orbits. Indeed, the bijectivity of $\pi$ follows from the study of the unbounded trace orbits.

}
\end{rem}

\subsubsection{Gaps, zeros of trace polynomials and $ \infty$-type energies }\label{subsec-gap-zero-infty}\

With this coding map in hand, we can give a complete description for the gaps of the spectrum. We will see that the gap edges, zeros of trace polynomials and $\infty$-type energies are related in a delicate but beautiful way.

Define the set of zeros of trace polynomials as
\begin{equation}\label{Zero}
  \ZZ^o:=\bigcup_{n \text{ odd}}\{h_n=0\};\ \ \ \ \ZZ^e:=\bigcup_{n \text{ even}}\{h_n=0\};\ \ \ \  \ZZ:=\ZZ^o\cup \ZZ^e.
\end{equation}

We will see that,  if an energy is a gap edge, then its coding is eventually 2-periodic. Figure \ref{fig-graph} suggests that there are 7 types of eventually 2-periodic codings:
\begin{equation}\label{gap-edge}
\begin{cases}
  \mathcal E_{l}^o &:=  \{\omega\in \Omega_\infty:
  \omega=w(2_o1_e)^\infty \text{ for some } w\in \Omega_\ast\}\\
  \mathcal E_{r}^o &:= \{\omega\in \Omega_\infty:
  \omega=w(0_o3_{el})^\infty \text{ for some } w\in \Omega_\ast\}\\
   \mathcal E_{l}^e &:=  \{\omega\in \Omega_\infty:
  \omega=w(0_e3_{or})^\infty \text{ for some } w\in \Omega_\ast\}\\
  \mathcal E_{r}^e&:=  \{\omega\in \Omega_\infty:
  \omega=w(2_e1_o)^\infty \text{ for some } w\in \Omega_\ast\}\\
  \widetilde {\mathcal E}_{l} &:=  \{\omega\in \Omega_\infty:
  \omega=w(0_o3_{er})^\infty \text{ for some } w\in \Omega_\ast\}\\
   \widetilde{\mathcal E}_{r} &:=  \{\omega\in \Omega_\infty:
  \omega=w(0_e3_{ol})^\infty \text{ for some } w\in \Omega_\ast\}\\
   \mathcal F &:=  \{\omega\in \Omega_\infty:
  \omega=w(2_e2_{o})^\infty \text{ for some } w\in \Omega_\ast\}
  \end{cases}
\end{equation}
where $\Omega_\ast$ is the set of finite admissible words.
 Define two special codings as follows:
\begin{equation}\label{two-boundary}
  \omega_\ast:=3_{el}(0_o3_{el})^\infty;\ \ \ \omega^\ast:=(0_e3_{or})^\infty.
\end{equation}

Denote the set of   gaps of $\sigma_\lambda$ by $\GG.$ We have the following complete description on the structure of $\GG.$

\begin{thm}\label{main-gap-zero-inftyenergy}
With the notations above, we have:

i) $\min \sigma_\lambda=\pi(\omega_\ast)$ and $\max \sigma_\lambda=\pi(\omega^\ast)$.

ii) (Type-I gaps): There exists a bijection $\iota^o: \pi(\mathcal E_{l}^o)\to \pi (\mathcal E_{r}^o\setminus \{\omega_\ast\})$   such that $E<\iota^o(E)$ for any $E\in \pi(\mathcal E_{l}^o)$ and
\begin{equation*}\label{}
  \GG_I^o:=\{(E,\iota^o(E)): E\in \pi(\mathcal E_{l}^o)\}\subset \GG.
\end{equation*}
There exists a bijection $\iota^e: \pi(\mathcal E_{r}^e)\to \pi (\mathcal E_{l}^e\setminus \{\omega^\ast\})$   such that $\iota^e(E)<E$ for any $E\in \pi(\mathcal E_{r}^e)$ and
\begin{equation*}\label{}
  \GG_I^e:=\{(\iota^e(E),E): E\in \pi(\mathcal E_{r}^e)\}\subset \GG.
  \end{equation*}

  iii) (Type-II gaps): There exists a bijection $\iota: \pi(\widetilde {\mathcal E}_l)\to \pi (\widetilde {\mathcal E}_r)$   such that $E<\iota(E)$ for any $E\in \pi(\widetilde {\mathcal E}_l)$ and
\begin{equation*}\label{}
  \GG_{II}:=\{(E,\iota(E)): E\in \pi(\widetilde {\mathcal E}_l)\}\subset \GG.
\end{equation*}

 iv) $\GG=\GG_I\cup\GG_{II}$, where $\GG_I:=\GG_I^o\cup\GG_{I}^e.$

 v) (Coding of zeros): $\pi(\mathcal E_{l}^o)=\ZZ^o$ and $\pi(\mathcal E_{r}^e)=\ZZ^e$.

 vi) ($\infty$-type energies):  For any $\mathcal E\in \{\mathcal E_{l}^e,\mathcal E_{r}^o, \widetilde {\mathcal E}_l,\widetilde {\mathcal E}_r\}$, we have
 \begin{equation*}\label{infty-I-II}
   \pi(\mathcal E)\subset \mathscr{E}_\infty \text{ and }\ \  \pi(\mathcal E) \text{ is dense in } \sigma_\lambda.
 \end{equation*}

\end{thm}

We say that the  energy $E\in \sigma_\lambda$ is of {\it type $\infty_I$ (type $\infty_{II}$)} if $E\in \pi(\mathcal E_{l}^e\cup \mathcal E_{r}^o)$ ($ \pi(\widetilde {\mathcal E}_l\cup\widetilde {\mathcal E}_r))$.

\begin{rem}
{\rm
1) For a gap  in $\GG_I$, one of its end point is a zero of trace polynomials, another endpoint if of type $\infty_I$. For a gap in $\GG_{II}$, both endpoints of it are of type $\infty_{II}$. So the codings of all the gap edges are  eventually 2-periodic. However, by ii)-iv), $E\in \pi(\mathcal F)$ is not an edge of any gap. See Theorem \ref{main-gap-label} and Remark \ref{rem-label} for the special role played by $\mathcal F.$

2) $\GG_I^e, \GG_I^o, \GG_{II}$ are corresponding to the green, blue and pink edges of $\mathbb G$, respectively, see Figure \ref{fig-graph}. We also remark that $\GG_I$ are  exact the gaps in \cite{BBG} Theorem 4 which open linearly, while $\GG_{II}$ are  exact the gaps  which open exponentially.

}
\end{rem}

\subsubsection{Integrated density of states and Gap labelling of the spectrum}\label{statement-IDS}\

Now we compute the IDS for every energy in the spectrum. The labels on the edges of the graph $\mathbb G$ play an essential role in this computation.

Define the infinite  binary tree $\Sigma_\infty$
 and the evaluation map
  $\varepsilon:\Sigma_\infty\to [0,1]$ as
 \begin{equation}\label{Sigma-infty}
 \Sigma_\infty:=\{0,1\}^\infty; \ \ \ \varepsilon(\sigma)=\sum_{n\ge1}\frac{\sigma_n}{2^n}.
 \end{equation}

 Define a map $\Pi:\Omega_\infty\to \Sigma_\infty$ as follows: given $\omega\in \Omega_\infty$, follow the infinite path $\omega=\omega_0\omega_1\omega_2\cdots$ in $\mathbb G$,  record the label $w_i$ on the edge $\omega_{i-1}\omega_{i}$ and concatenate them to get
 \begin{equation}\label{def-Pi}
 \Pi(\omega):=i(\omega)w_1w_2w_3\cdots\in \Sigma_\infty,
 \end{equation}
 where $i(\omega)=\emptyset$ if $\omega_0=0_e$ and $i(\omega)=0$ if $\omega_0=3_{el}$. Here we use the convention that $\emptyset w=w\emptyset=w$ for a finite word $w$.

 Denote by $N(E)$ the IDS of $E\in\sigma_\lambda$.
 We have the following:

 \begin{thm}\label{main-comp-IDS}
  $\Pi:\Omega_\infty\to \Sigma_\infty$ is surjective and  for any
  $\omega\in \Omega_\infty$,
\begin{equation*}\label{IDS}
N(\pi(\omega))=\varepsilon(\Pi(\omega)).
\end{equation*}
In other words, we have  the following commutative diagram:
\begin{equation*}\label{commu-diag}
\begin{CD}
\Omega_\infty @>\pi>>\sigma_\lambda\\
@VV\Pi V @VV N V\\
\Sigma_\infty @>\varepsilon>>[0,1]
\end{CD}
\end{equation*}

\end{thm}

 Later, we will show that $\Pi$ is not injective. Indeed, we have
$
\Pi(\widetilde {\mathcal E}_l)=\Pi(\widetilde {\mathcal E}_r),
$
and $\Pi:\Omega_\infty\setminus\widetilde {\mathcal E}_r\to\Sigma_\infty$ is bijective, see Proposition \ref{Pi-basic} for detail. We also point out that $\Pi$ play an essential role in the study of the symbolic space $\Omega_\infty.$

It was conjectured in
\cite{Lu} and proved in \cite{BBG} that all the possible gaps of $\sigma_\lambda$ are open and
 labeled by dyadic numbers or dyadic numbers divided by $3$ in $(0,1)$. However, to the best of our knowledge, the complete characterization of the
 labels of these gaps seems unknown in the literature (compare \cite{BBG} Theorem 4 (iv) and \cite{B2} Sec. 5.4).
 With the coding map in hand, we can determine the exact set of numbers which label the gaps.

  It is known  that $N$ is constant on a gap $G$ of the spectrum and  we denote this number by $N(G)$.

   \begin{thm}\label{main-gap-label}
We have the gap labelling:
\begin{equation*}\label{gap-lab}
  \{N(G): G\in \mathscr G_I\}=\mathscr D\cap (0,1)\ \ \text{ and }\ \  \{N(G): G\in{\mathscr G}_{II}\}=\mathscr E,
\end{equation*}
where $\mathscr D$ is the set of positive dyadic numbers and
\begin{equation}\label{big-E}
\mathscr E=\varepsilon(\Pi(\widetilde {\mathcal E}_l))=[(\mathscr D/{3}\setminus \mathscr D)\cap(0,1)]\setminus \varepsilon(\Pi(\mathcal F)).
\end{equation}
\end{thm}

\begin{rem}\label{rem-label}
Assume $s\in (0,1)$ is a dyadic number divided by 3.
\eqref{big-E} tells us that  $s$ is  a label of some gap of $\sigma_\lambda$ if and only if $s\notin \varepsilon(\Pi(\mathcal F))$.
\end{rem}

\subsection{Sketch of the main ideas}\

Very roughly, our idea is the following: using periodic approximations to construct a family of nested optimal coverings; find out all the possible band configurations between consecutive levels; define the types of bands and  figure out the evolution laws of types;  finally  construct  the symbolic space. We explain more details in the following.

\subsubsection{Nested structure}\label{sec-NS}\

The family of coverings is a nested structure which we define now.

For any $n\ge0$, let $\I_n=\{I_{n,1},\cdots, I_{n,k_n}\}$ be a family of compact, non degenerate intervals. $\I:=\{\I_n:n\ge0\}$ is called a {\it nested structure(NS)}, if

i) $\I$ is {\it optimal}: for any $n\ge0$, $i\ne j$, $I_{n,i}\not\subset I_{n,j}$;

ii) $\I$ is {\it nested}: for any $I\in \I_{n+1}$, there is a unique $J\in \I_n$ such that $I\subset J;$

iii) $\I$ is {\it minimal}: for any $J\in\I_n$, there exists $I\in \I_{n+1}$ such that $I\subset J.$

Assume $\I:=\{\I_n: n\ge0\}$ is a NS. By ii), $\{\bigcup_{i=1}^{k_n}I_{n,i}:n\ge0\}$ is a decreasing sequence of nonempty compact sets. Define
$$
A(\I):=\bigcap_{n\ge0}\bigcup_{i=1}^{k_n}I_{n,i}.
$$
We call $A(\I)$ the {\it limit set} of $\I$.

Assume $\I=\{\I_n:n\ge0\}$ is a NS. $\I$ is called a {\it separating nested structure(SNS)}, if moreover,

iv) $\I$ is {\it separating}: for each $n$, $\I_n$ is a disjoint family.

We will construct a NS such that its limit set is $\sigma_\lambda$. Then we  choose a sub-NS from it such that this sub-NS is a SNS and for which we can estimate the dimension.

\subsubsection{NS via Periodic approximations }\

For all $n\ge 0$, define the $n$-th {\it periodic approximation} of $\sigma_\lambda$ as
\begin{equation}\label{sigman}
\sigma_{\lambda,n}:=\{E\in\R: |h_n(E)|\le 2\}.
\end{equation}
 By Floquet theory, $\sigma_{\lambda,n}$ is made of $2^n$ non-overlapping intervals. Following the convention, we call these intervals bands.    We denote this family of bands by $\B_n$ and  code it by $\Sigma_n=\{0,1\}^n$ as
\begin{equation}\label{B-n}
  \B_n:=\{B_\sigma: \sigma\in \Sigma_n\}.
\end{equation}
See Sec. \ref{sec-code-B-n} for detail.
 It is well-known that  $\sigma_{\lambda,n}\to \sigma_\lambda$ in Hausdorff distance. In general, $\B_n$ is not a covering of $\sigma_\lambda$. However we will show that $\{\B_n\cup\B_{n+1}:n\ge 0\}$ is a family of nested coverings. We will see that $\B_n\cup\B_{n+1}$ is always not optimal. To get an optimal covering, we simply delete those $B\in \B_{n+1}$ which is contained in some $\hat B\in \B_{n}$.
 Denote by $\BB_n$ the resulting optimal covering. Formally
 \begin{equation}\label{def-BB-n}
   \BB_n:=\B_n\cup\{B\in \B_{n+1}: B\not\subset \hat B \text{ for any } \hat B\in \B_n\}.
\end{equation}
Now $\BB:=\{\BB_n: n\ge0\}$ become  a NS and its limit set coincides with $\sigma_\lambda$, see Corollary \ref{NS-BB}. In principle, one can construct a symbolic space and the related coding map of the spectrum  by following the sequence $\{\BB_n:n\ge0\}.$
 However, we warn that when $\lambda>0$ is small, $\BB$ is never  a SNS. As a consequence, it is very difficult to show the injectivity of the coding map.

 To decide which bands in $\B_{n+1}$ should be deleted, we need to study the configurations of bands in $\B_n\cup \B_{n+1}$. This is the content of  Lemma \ref{tech-1}--the first key lemma of the paper.

\subsubsection{Types, graph-directed structure and coding}\label{idea-graph}\

The next step is to study how  the bands in $\BB_{n+1}$ are situated in their father bands in $\BB_n.$ The crucial concept here is the {\it  type} of the band.

By  performing a renormalization process, we observe that every band in $\BB_n$ has certain type.  The type of $B\in \BB_n$ determines its sons  $\hat B\in\BB_{n+1}$ and  their types, see Figure \ref{fig-type} for some intuition. By detailed analysis, we obtain an alphabet $\A$ for the types and a directed graph $\mathbb G$ for the evolution of types.

  Recall that $\BB_n$ is the union of $\B_n$ and a subset  of $\B_{n+1}$. Roughly speaking, there are four types for $B\in \BB_n$: $0,1,2,3$.  For $k=0,1,2$, $B$ has type $k$ if $B\in \B_n$ and $\#\partial B\cap \RR_{n-1}=k$ where $\RR_{n-1}$ is the set of zeros of trace polynomials $h_0,\cdots, h_{n-1}$ (see \eqref{Z-n-R-n}). $B$ has type $3$ if $B\notin \B_{n}$. Since we aim  to give an order-preserving coding for the spectrum, the order of the son bands become important. With this in mind, we observe that the bands with same types but with different parities of levels have different evolution laws: they are reflectional symmetric. This observation suggests that we need to distinguish $0_o, 0_e$ etc.  This explains the definition \eqref{alphabet} except for type 3, for which we need to distinguish left and right further. The evolution laws \eqref{adm} of types are hidden in Lemma \ref{tech-2}--the second key lemma of this paper. Now we can build a related directed graph $\mathbb G$.

  Once $\mathbb G$ is constructed, we can first  code the covering $\BB_n$ by $\Omega_n$ ( the set of admissible words of level $n$, see \eqref{Omega-n-ast})  as
  $\BB_n:=\{I_w: w\in \Omega_n\}$ and then code the spectrum $\sigma_\lambda=A(\BB)$ by the symbolic space $\Omega_\infty$ as
  $\pi(\omega):=\bigcap_{n}I_{\omega|_n}$. As we have mentioned, it is easy to show that $\pi$ is  surjective and continuous. However to show that $\pi$ is injective, we need further information on $\infty$-type energies, see Sec. \ref{idea-order-gap}.

\subsubsection{Sub-SNS and the  lower bound of the dimension}\label{idea-low-dim}\

Now consider the sub-alphabet $\widetilde A=\{1_e,1_o,2_e,2_o\}$ and the sub directed graph $\widetilde {\mathbb G}$ (enclosed by the red square in Figure \ref{fig-graph}) with  restricted admissible rules:
$$
1_e\to 2_o;\ \ 1_o\to 2_e;\ \ 2_e\to 1_o, 2_o;\ \  2_o\to 1_e, 2_e.
$$
Consider the sub-NS $\widetilde \BB:=\{\widetilde \BB_n: n\ge1\}$ defined by
$$
\widetilde \BB_n:=\{I_{0_e1_ow}: w=w_1\cdots w_n, w_i\in \widetilde A\}.
$$
Obviously its limit set $A(\widetilde \BB)$ is a subset of the spectrum. Indeed, each energy in $A(\widetilde \BB)$ has a bounded trace orbit of bound 2. In particular,  $A(\widetilde \BB)\subset B_\infty$. Dynamically, if $E\in A(\widetilde \BB)$, then for any $n\ge0$,
$$
f^n(\ell(E))\in [-2,2]\times [-2,2].
$$
In some sense, $A(\widetilde \BB)$ is the portion of $\sigma_\lambda$ on which the trace polynomials behave in a uniformly hyperbolic way.
We can show that $\widetilde \BB$ is a SNS and there is a uniform estimation for the length of the bands in $\widetilde \BB_n$, see Proposition \ref{SNS}. Now by a general result for SNS(see Proposition \ref{SNS-lower-dim}), we obtain a lower bound for the Hausdorff dimension of the spectrum.

\subsubsection{ Order, Gaps and $\infty$-type energies}\label{idea-order-gap}\

Now we explain the orders on $\Omega_\infty.$ There are two orders on $\Omega_\infty$, both  are induced from orders on $\A$. So we start from the orders on $\A$.

We have two goals for designing an order on $\A$: 1) to remember the order of the son-bands inside a father-band; 2) to remember whether two neighbour son-bands overlap. As we have mentioned several times, if $\lambda$ is small, the overlapping of some son-bands are unavoidable. To obtain a nice coding, it is essential  to remember whether two bands overlap.
 For this purpose, we define the following two ``orders" to compare two bands.

Assume   $I=[a,b], J=[c,d]$. We define
 \begin{equation}\label{order-bands}
 \begin{cases}
   I\prec J,   & \text{ if }\  a<c; b<d,\\
   I<J,& \text{ if }\  b<c.
   \end{cases}
 \end{equation}
 $I<J$ means $I$ and $J$ are disjoint and $J$ is on the right of $I$. While $I\prec J$ means that
  $I\not\subset J, J\not\subset I$ and $J$ is on the right of $I$, however $I$ and $J$ may overlap.
  It is obvious that $I<J$ implies $I\prec J$. Lemmas \ref{type-0}-\ref{type-2} summarize the order configurations of son bands for all types of bands.  These motivate the definitions of two orders $\preceq$ and $\le$ on $\A$ (see \eqref{order} and \eqref{order-strong}).

    Now we can  extend these orders to $\Omega_\infty$ in an obvious way to get the two orders $\preceq$ and $\le$ on $\Omega_\infty$. We will  show that $\preceq$ is a total order on $\Omega_\infty$, which is the order appearing in Theorem \ref{main-coding-spectrum}. $\le$ is stronger than $\preceq$ in the following sense: if $\omega<\omega'$, then  $\pi(\omega)<\pi(\omega')$; however, if $\omega\prec\omega'$ but $\omega\not<\omega'$,  we only have $\pi(\omega)\le\pi(\omega')$. It is easy to show the former statement, which follows from the disjointness  of the bands. It is quite subtle to show the latter statement because the condition means that the bands containing $\pi(\omega)$ and $\pi(\omega')$ overlap with each other. We obtain the inequality by a detailed analysis on the combinatorics of $\ZZ.$ Note  that our ultimate goal is to show the strict inequality, which is even harder. We need to study the dynamics of the trace map of $f$.

  Next we study the gap of $\sigma_\lambda$. Since $(\Omega_\infty, \preceq)$ is totally ordered, the notion of gap for $\Omega_\infty$ is available. We will at first study the gap structure on $\Omega_\infty$ (see Theorem \ref{symbolic-gap}), and then descend them to $\sigma_\lambda$ by the coding map $\pi$.

  Now we discuss $\infty$-type energy. Recall  that $\B_n\cup\B_{n+1}$ is a covering of the spectrum, this implies that for any $E\in \sigma_\lambda$ and any $n\ge0$, $|h_n(E)|\le2$ or $|h_{n+1}(E)|\le2$. Inspired by the Thue-Morse case \cite{LQY}, to produce a unbounded trace orbit, we can ask for  the following mechanism: there exists $N$ such that
  \begin{equation}\label{ge2-le2}
  |h_{N+2k}(E)|\le 2;\ \ \  \  |h_{N+2k+1}(E)|>2, \ \ \text{ for any } k.
  \end{equation}
  Dynamically, this means that for $n$ large enough (if $N$ is even),
  $$
  f^{n}(\ell(E))\in [-2,2]\times [-2,2]^c.
  $$
  If $E$ has coding $\omega$, \eqref{ge2-le2} implies  that necessarily $\omega$ is eventually $0_\ast3_\ast0_\ast3_\ast\cdots$. With this in mind, we finally succeed to show that all the energies with coding eventually $(0_a3_b)^\infty$ are of $\infty$-type.

  As one will see in Sec. \ref{sec-infty-energy}, it is relatively easy to show that $E$ is of $\infty$-type if the coding of $E$ is eventually $(0_o3_{el})^\infty$ or $(0_e3_{or})^\infty$; however it is much harder  to show that $E$ is of $\infty$-type if the coding of $E$ is eventually $(0_o3_{er})^\infty$ or $(0_e3_{ol})^\infty$. For the latter case, we need to study the local dynamics of $f$ around the fixed point $(-1,-1)$. The injectivity of $\pi$ is also a consequence of this study.

\subsubsection{Labelled directed graph and  integrated density of states}\

 At last, we explain the labels on the edges of $\mathbb G$. Assume $B_\sigma\in \BB_n$ has type $\alpha$, $B_\tau\in \BB_{n+1}$ has type $\beta$ and $B_\tau\subset B_\sigma$. Then we will show that $\alpha\to \beta$ and $\sigma$ is a prefix of $\tau$. Write $\tau=\sigma\sigma'$, then the label of the edge $\alpha\beta$ is $\alpha'$.
Now by Lemmas \ref{type-3}-\ref{type-2}. we can label all the edges in $\mathbb G.$

Next we compute the IDS for $E\in \ZZ$. To do this we need to study the structure of $\ZZ$, see Proposition \ref{zero-order}. Once this is computed, since $\ZZ$ is dense in $\sigma_\lambda,$ by the continuity of $N$, we can compute the IDS for all $E\in \sigma_\lambda.$

Now we come to the last point--gap labelling. By Theorem \ref{main-gap-zero-inftyenergy}, we know the coding of the gap edges. To compute the IDS of gaps, we just follow the path
$$
E\to \pi^{-1}(E)\to \Pi\circ \pi^{-1}(E)\to\varepsilon\circ\Pi\circ \pi^{-1}(E).
$$

  {\it Final remark about the proof:} As we will see in all the properties which need  to distinguish  the parities of the level $n$ that the statements for the odd case are completely symmetric with that for the even case, so are the proofs. So in the rest of the paper, usually we only give the proofs for the odd case and leave the routine checking for the even case to the reader.

\subsection{Structure of the paper}\

The rest of this paper is organized as follows. In Sec. \ref{sec-optimal-covering}, we study the  band configurations  of  $\B_n\cup\B_{n+1}$ and give a  characterization for $\BB_n$, we also study the properties of $\ZZ$.  In Sec. \ref{sec-type-evo-order},  we introduce the concept of type and derive the evolution law and the order of types and the labels of the directed edges.
In Sec. \ref{sec-symbolic}, we study  the symbolic space $\Omega_\infty$. In Sec. \ref{sec-coding-Hausdorff}, we prove a weak version of Theorem \ref{main-coding-spectrum}, and  apply it to  prove Theorem \ref{main-low-dim}.
 In Sec. \ref{sec-infty-energy}, by a detailed  study of  $\infty$-type energies, we
 prove Theorems \ref{main-infty-energy}, \ref{main-coding-spectrum} and  \ref{main-gap-zero-inftyenergy}.
In Sec. \ref{sec-ids-gap}, we prove Theorems \ref{main-comp-IDS} and \ref{main-gap-label}. In Sec. \ref{sec-pf-tech-lemma}, we prove the technical lemmas and the related results.

\section{Optimal covering $\BB_n$ and the property  of $\ZZ$}\label{sec-optimal-covering}

At first, we study the  band configurations  of  $\B_n\cup\B_{n+1}$ and give a  characterization for $\BB_n$. Then we study the property of $\ZZ.$

\subsection{The optimal covering $\BB_n$}

\subsubsection{Coverings of $\sigma_\lambda$ via periodic approximations}\

It is known that (\cite{BBG,Lu}) the sequence of trace polynomials of PDH $\{h_n:n\ge 0\}$ defined by \eqref{trace-poly-n} satisfy the
following initial conditions and recurrence relation:
\begin{equation}\label{recurrence}
\begin{cases}
h_0(E)=E-\lambda,\ \  h_1(E)=E^2-\lambda^2-2;\\
h_{n+1}(E)=h_n(E)\left(h_{n-1}^2(E)-2\right)-2,\ \  (\forall n\ge1).
\end{cases}
\end{equation}

We have defined $\sigma_{\lambda,n}$ in \eqref{sigman} and $\B_n$ in \eqref{B-n}.
By Floquet theory (see for example \cite{T} chaper 4), $\sigma_{\lambda,n}$ is the spectrum of $H_{V^{(n)}}$,
where $V^{(n)}$ is $2^n$-periodic and
$$V^{(n)}(j)=\lambda V_\xi(j),\ \ \ j=0,\cdots, 2^{n}-1.$$
 $\sigma_{\lambda,n}$ is called the $n$-th
 {\it periodic approximation} of  $\sigma_\lambda$.
 We know that $\B_n$ is a non-overlapping family. Indeed, we can say more:

 \begin{prop}\label{gap-open-sigma-n}
 For any $n\ge1$,
 all the possible gaps of $\sigma_{\lambda,n}$ are open, i.e. $\B_n$ is a disjoint family.
 \end{prop}

 We will prove this proposition in Sec. \ref{sec-Z-band-edge}. In general,  $\B_n$ is not a covering of $\sigma_\lambda.$ However, we have

\begin{lem}\label{covering-basic}
i) For any $n\ge0$,
\begin{equation*}\label{sigma-n-inclusion}
  \sigma_{\lambda,n+2}\subset\sigma_{\lambda,n}\cup \sigma_{\lambda,n+1}.
\end{equation*}
Thus $\{\sigma_{\lambda,n}\cup \sigma_{\lambda,n+1}:n\ge0\}$ is a decreasing sequence.
Moreover,
\begin{equation}\label{spec-covering}
\sigma_\lambda=\bigcap_{n\ge0}(\sigma_{\lambda,n}\cup\sigma_{\lambda,n+1}).
\end{equation}

ii) We have
\begin{equation}\label{band-length}
  \max_{B\in \B_n}|B| \to 0, \ \ n\to\infty.
\end{equation}
\end{lem}

See \cite{LQ} Lemma A.2. for a proof (there all the results are proven for TMH, but the proof for PDH is the same).

As a consequence, $\B_n\cup\B_{n+1}$ is a covering of $\sigma_\lambda.$

\subsubsection{The coding of   $\B_n$}\label{sec-code-B-n} \

We  mentioned that $\B_n\cup\B_{n+1}$ is not an optimal covering of $\sigma_\lambda$. To construct an optimal covering, we need to throw some bands in $\B_{n+1}$. To proceed,  we  give a coding for $\B_n$.

Assume  $n\ge 1$.
Let $\Sigma_n:=\{0,1\}^n$ and $\mathbb{M}_n:=\{0,1,\cdots,2^n-1\}$.
By convention, $\Sigma_0:=\{\emptyset\}$ and $\mathbb{M}_0=\{0\}$.
 We endow $\Sigma_n$ with the lexicographical order and denote by $\le$.
 If $\sigma\ne 1^n (\sigma\ne 0^n)$, denote  the {\it successor(predecessor)} of $\sigma$
 by $\sigma^+(\sigma^-)$.
For $\sigma=\sigma_1\cdots\sigma_n\in\Sigma_n$ and $k\le n$, we define the $k$-th prefix of $\sigma$ as
$$
\sigma|_k:=\sigma_1\cdots\sigma_{k}.
$$
Assume $n,k\ge0$ and $\sigma\in\Sigma_n, \tau\in \Sigma_k$. Then $\sigma\wedge\tau$ is the maximal common
prefix of $\sigma$ and $\tau$;
 $\sigma\tau\in \sigma_{n+k}$
is the concatenation
of $\sigma$ and $\tau$. If $\sigma$ is a prefix of $\tau$, we write $\sigma\lhd \tau.$
We denote the length of $\sigma$ by $|\sigma|$. By convention, $|\emptyset|=0.$

When $n=0$, define $\varpi_0:\Sigma_0\to \mathbb{M}_0$ as $\varpi_0(\emptyset)=0.$
When $n\ge 1$, define
$\varpi_n:\Sigma_n\to \mathbb{M}_n$ as
\begin{equation}\label{varpi-n}
\varpi_n(\sigma):=\sum_{j=1}^{n}\sigma_j2^{n-j}=2^n\sum_{j=1}^n\frac{\sigma_j}{2^j}, \ \ \ \text{where}\ \ \sigma=\sigma_1\cdots\sigma_n.
\end{equation}
It is direct to check that $\varpi_n:\Sigma_n\to \mathbb{M}_n$ is  bijective and order-preserving.

Now we list the $2^n$ bands in $\B_n$ from left to right with numbering $0,\cdots, 2^n-1.$
Then we denote the $j$-th bands by $B_{\varpi_n^{-1}(j)}$. In this way, we give a coding for
$\B_n.$

For any $n\ge0$, combining with Proposition \ref{gap-open-sigma-n}, we have
\begin{equation*}\label{}
\B_n=\{B_{\sigma}: \sigma\in\Sigma_n\}\ \ \  \text{and}\ \ B_\sigma< B_{\sigma'} \text{ if }  \sigma<\sigma'.
\end{equation*}
Now we have
\begin{equation}\label{sigma-lamb-n}
\sigma_{\lambda,n}=\bigcup_{B\in \B_n} B=\bigcup_{\sigma\in \Sigma_n}B_\sigma.
\end{equation}
We say that all the bands in $\B_n$ are of {\it level}-$n$.

For any $\sigma\in \Sigma_n$, write
\begin{equation*}\label{B-sigma}
  [a_\sigma, b_\sigma]:=B_\sigma.
\end{equation*}

For any $n\ge0$, define
\begin{equation}\label{Z-n-R-n}
  \ZZ_n:=\{E\in \R: h_n(E)=0\}\ \ \text{ and }\ \ \RR_n:=\bigcup_{j=0}^{n}\ZZ_j.
\end{equation}

By Floquet theory, in each band $B_\sigma$, there is exactly one zero of $h_n$, which we denote by $z_\sigma.$ So we have
\begin{equation}\label{coding-Z-n}
  a_\sigma<z_\sigma<b_\sigma;\ \   \ZZ_n=\{z_\sigma: \sigma\in \Sigma_n\}; \ \  z_\sigma< z_\tau \text{ if } \sigma< \tau.
\end{equation}

As an example, we compute $\B_0$ and $\B_1$. By \eqref{recurrence}, we have
\begin{equation}\label{B-0-1}
\begin{cases}
B_\emptyset=[\lambda-2,\lambda+2],\ \ \
B_0=[-\sqrt{\lambda^2+4},-\lambda],\ \ \ B_1=[\lambda,\sqrt{\lambda^2+4}]\\
\B_0=\{B_\emptyset\},\ \ \ \B_1=\{B_0,B_1\}.
\end{cases}
\end{equation}

\subsubsection{The band configurations  of  $\B_n\cup\B_{n+1}$ and characterization of $\BB_n$} \

Recall that we have defined two ``orders" for bands in \eqref{order-bands}. The band configurations of $\B_n\cup\B_{n+1}$ is summarized  in the following technical lemma:

 \begin{lem}\label{tech-1}
 Fix $n\ge0$ and $\sigma\in \Sigma_n$.

 (i) Assume $n$ is odd. Then $b_{\sigma0}=z_\sigma, B_{\sigma0}\subset B_\sigma$ and $a_{\sigma1}\notin \RR_n$.
 Moreover,
 $$
 \begin{cases}
 B_{\sigma0}=[a_{\sigma},z_\sigma],& \text{ if }\  a_\sigma\in \RR_{n-1}\\
 B_{\sigma0}\subset (a_{\sigma},z_\sigma],&\text{ if }\  a_\sigma\notin\RR_{n-1}
 \end{cases};\ \ \ \ \
 \begin{cases}
 B_{\sigma1}\subset(z_\sigma,b_\sigma],&\text{ if }\  b_\sigma\in \RR_{n-1}\\
 [z_\sigma,b_\sigma]\prec B_{\sigma1},&\text{ if }\  b_\sigma\notin\RR_{n-1}
 \end{cases}
 $$

 (ii) Assume $n$ is even. Then $a_{\sigma1}=z_\sigma, B_{\sigma1}\subset B_\sigma$ and $b_{\sigma0}\notin \RR_n$.
 Moreover,
 $$
 \begin{cases}
 B_{\sigma1}=[z_{\sigma},b_\sigma],& \text{ if }\ b_\sigma\in \RR_{n-1}\\
 B_{\sigma1}\subset [z_{\sigma},b_\sigma),&\text{ if }\  b_\sigma\notin\RR_{n-1}
 \end{cases};\ \ \ \ \
 \begin{cases}
 B_{\sigma0}\subset [a_\sigma,z_\sigma),&\text{ if }\  a_\sigma\in \RR_{n-1}\\
  B_{\sigma0}\prec [a_\sigma,z_\sigma],& \text{ if }\ a_\sigma\notin\RR_{n-1}
 \end{cases}
 $$

 iii) $B_{\sigma0}<B_{\sigma1}$. If $\sigma> 0^n$, then $B_{\sigma^-}<B_{\sigma0}$; if $\sigma< 1^n$, then $B_{\sigma1}<B_{\sigma^+}$.
 \end{lem}

Here we use the convention that $\RR_{-1}=\emptyset.$

\begin{figure}[htbp]
\begin{minipage}[t]{0.3\linewidth}
\centering
\includegraphics[width=1.5in]{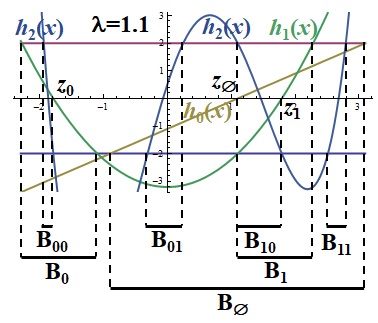}
\caption*{(a)}
\end{minipage}%
\begin{minipage}[t]{0.3\linewidth}
\centering
\includegraphics[width=1.5in]{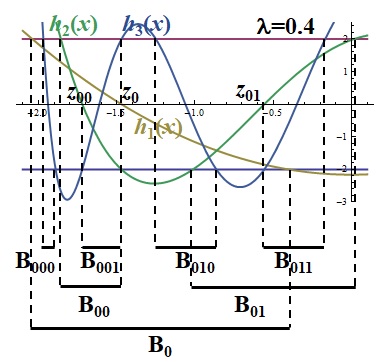}
\caption*{(b)}
\end{minipage}%
\begin{minipage}[t]{0.3\linewidth}
\centering
\includegraphics[width=1.5in]{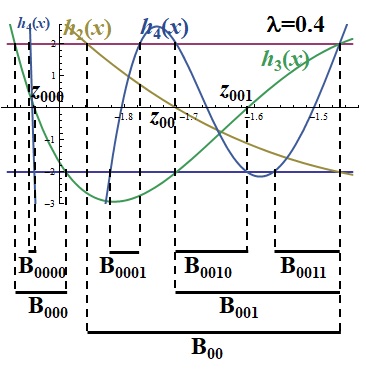}
\caption*{(c)}
\end{minipage}%

\begin{minipage}[t]{0.3\linewidth}
\centering
\includegraphics[width=1.5in]{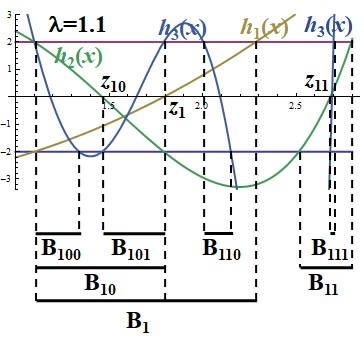}
\caption*{(d)}
\end{minipage}%
\begin{minipage}[t]{0.3\linewidth}
\centering
\includegraphics[width=1.5in]{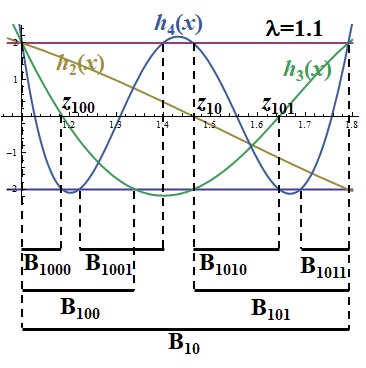}
\caption*{(e)}
\end{minipage}%
\begin{minipage}[t]{0.3\linewidth}
\centering
\includegraphics[width=1.5in]{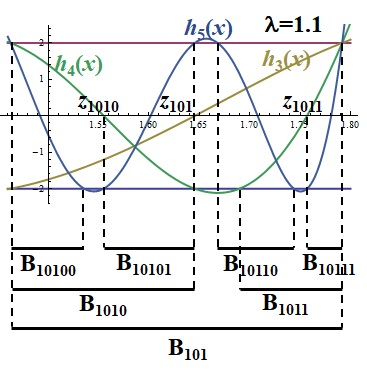}
\caption*{(f)}
\end{minipage}%
\caption{Band configurations of $\B_n\cup \B_{n+1}\cup \B_{n+2}$}\label{fig-type}
\end{figure}

One can see Figure \ref{fig-type} for some concrete examples.

We will prove this lemma in Sec. \ref{sec-pf-Lem2.3}.  For the moment, we will use it to  characterize the bands in  $\BB_{n}$.

Recall that $\BB_n$ is defined by \eqref{def-BB-n}.  By Proposition \ref{gap-open-sigma-n} and Lemma \ref{covering-basic} i), $\BB_n$ is an optimal covering of $\sigma_\lambda.$  By the definition, we always have $\B_n\subset \BB_n$.
Now we give a characterization  for a band $B\in \B_{n+1}$ in $\BB_n.$

\begin{prop}\label{BB-n}
 Assume $\sigma\in \Sigma_{n+1}$. Then $B_\sigma\in \BB_n$ if and only if $B_\sigma\not\subset B_{\sigma|_n}$.
\end{prop}

\begin{proof}
Write $\tau:=\sigma|_n$. If $B_\sigma\subset B_\tau\in \B_n$, by \eqref{def-BB-n}, $B_\sigma\notin\BB_n.$

Now assume $B_\sigma\not\subset B_\tau$. By Lemma \ref{tech-1} iii),
$$
B_{\tau^-}<B_\sigma, \text{ if } \tau> 0^{n};\ \ \ B_\sigma<B_{\tau^+}, \text{ if } \tau<1^{n}.
$$
By the coding and the order of the bands in  $\B_n$, we conclude that $B_\sigma\not\subset B$ for any $B\in \B_n$. So by \eqref{def-BB-n}, $B_\sigma\in \BB_n.$
\hfill $\Box$
\end{proof}

 As an example, by \eqref{B-0-1}, we have
 \begin{equation}\label{initial}
B_1\subset B_\emptyset\ \ \ \text{ and }\ \ \ B_0\prec B_\emptyset.
\end{equation}
Consequently we have (see Figure \ref{fig-type} (a))
\begin{equation*}\label{BB-0}
\BB_0=\{B_0,B_\emptyset\}.
\end{equation*}

By direct computation, we have (see also Figure \ref{fig-type} (a))
\begin{equation*}\label{BB-1}
\BB_1=\{B_0,B_1,B_{01}, B_{11}\}.
\end{equation*}

\subsection{The property of $\ZZ$}\

At first, we have the following simple but useful observation:

\begin{lem}\label{basic-zeros}
i) For  $n\ge0$, if $E\in \ZZ_n$, then
\begin{equation}\label{0-m2-p2}
h_{n+1}(E)=-2 \ \ \text{ and }\ \
h_{k}(E)=2 \text{ for any }  k\ge n+2.
\end{equation}
Consequently, $\ZZ_n\cap \ZZ_m=\emptyset$ if $n\ne m.$

ii) For $n\ge1$ and $\sigma\in \Sigma_n$,
$${\rm int}(B_\sigma)\cap \RR_{n-1}=\emptyset\ \ \text{ and }\ \ {\rm int}(B_\sigma)\cap \RR_{n}=\{z_\sigma\}.$$

iii) $\ZZ\subset B_\infty\cap \sigma_\lambda$ and $\ZZ$ is dense in $\sigma_\lambda.$
\end{lem}
\begin{proof}
i) By \eqref{recurrence}, one check directly that the \eqref{0-m2-p2} holds. WLOG assume $n<m$. By \eqref{0-m2-p2}, $h_{n+1}(\ZZ_n)=\{-2\}$ and $h_{k}(\ZZ_n)=\{2\}$ for $k\ge n+2$. Hence $\ZZ_n\cap \ZZ_m=\emptyset$.

ii) If $E\in \RR_{n-1}$, then by i), $h_n(E)=\pm2$. By Floquet theory, $h_n({\rm int}(B_\sigma))=(-2,2)$. So ${\rm int}(B_\sigma)\cap \RR_{n-1}=\emptyset$. Still by Floquet theory, ${\rm int}(B_\sigma)\cap \ZZ_{n}=\{z_\sigma\}$. Since $\RR_n=\RR_{n-1}\cup\ZZ_n,$ we have  ${\rm int}(B_\sigma)\cap \RR_{n}=\{z_\sigma\}.$

iii) By i), $\ZZ\subset B_\infty.$ By i) and \eqref{spec-covering}, $\ZZ\subset \sigma_\lambda.$ Fix any $E\in \sigma_\lambda$. For any $n\ge0$, by \eqref{spec-covering}, there exists  $\sigma^{(n)}\in \Sigma_n\cup\Sigma_{n+1}$ such that $E\in B_{\sigma^{(n)}}$. Since $z_{\sigma^{(n)}}\in B_{\sigma^{(n)}},$ by \eqref{band-length}, we have
$
|z_{\sigma^{(n)}}-E|\to 0.
$
So $\ZZ$ is dense in $\sigma_\lambda.$
\hfill $\Box$
\end{proof}

\subsubsection{The order-preserving coding of $\ZZ$}\

Recall that $\ZZ$ is defined by \eqref{Zero}.
Define
\begin{equation}\label{Sigma-ast}
\Sigma_\ast:=\bigcup_{n=0}^\infty \Sigma_n.
\end{equation}
By \eqref{coding-Z-n} and Lemma \ref{basic-zeros} i), the map $\sigma\to z_\sigma$ is a bijection between $\Sigma_\ast$ and $\ZZ$. Thus we obtain a coding of $\ZZ$ as
\begin{equation*}\label{coding-ZZ}
\ZZ=\{z_\sigma: \sigma\in \Sigma_\ast\}.
\end{equation*}

In the following, we define an order on $\Sigma_\ast$ such that this coding preserves the orders.

At first we define an order $\preceq$ on $\{0,\emptyset,1\}$ as
$$
0\prec \emptyset \prec 1.
$$
For each  $n\ge1$, let $\preceq$ be the lexicographical order on $ \{0,\emptyset,1\}^n$.

Now we define an order $\preceq$ on $\Sigma_\ast$ as follows.
Given $\sigma\in \Sigma_n$ and $\tau\in \Sigma_m$.
Assume $n\le m$. Define $\tilde \sigma:=\sigma\emptyset^{m-n}$. Then define
\begin{equation}\label{order-Sigma-ast}
\begin{cases}
\sigma\prec \tau\ \ \  &\text{ if }\ \ \ \tilde \sigma \prec \tau,\\
\sigma=\tau\ \ \ &\text{ if } \ \ \  \tilde \sigma=\tau,\\
\tau\prec \sigma\ \ \  &\text{ if }\ \ \ \tau \prec \tilde \sigma.
\end{cases}
\end{equation}
One check directly that $\preceq$ is a total order on $\Sigma_\ast.$

\begin{prop}\label{zero-order}
Assume $\sigma,\tau\in \Sigma_\ast$. Then
$$
z_\sigma< z_\tau \ \ \ \text{ if and only if }\ \ \ \sigma\prec \tau.
$$
\end{prop}

It has the following consequence:

\begin{cor}\label{sgn-prod}
If $n\ge1$ and $\sigma\in \Sigma_n$, then for any $E\in {\rm int}(B_\sigma)$,
$$
\prod_{j=0}^{n-1}h_j(E)
\begin{cases}
<0,& \text{ if }\ \sigma_n=0\\
>0,& \text{ if }\  \sigma_n=1.
\end{cases}
$$
\end{cor}

We will prove Proposition \ref{zero-order} and Corollary \ref{sgn-prod} in Sec. \ref{sec-order-Z}.

\subsubsection{$\ZZ$ and band edges}\

Recall that  $\RR_{-1}=\emptyset.$

\begin{lem}\label{endsymbol}

 Assume $n,t\ge0$ and $\sigma\in \Sigma_n.$
 Then

 i) $b_{\sigma01^t}<a_{\sigma10^t}$ and
$
z_\sigma=
\begin{cases}
b_{\sigma01^t},&  \text{ if }\ n \text{ is odd};\\
a_{\sigma10^t},& \text{ if }\  n \text{ is even}.
\end{cases}
$

ii)
If $a_\sigma\in\RR_{n-1}$, then $a_{\sigma0^t}=a_\sigma$;
if $b_\sigma\in\RR_{n-1}$, then $b_{\sigma1^t}=b_\sigma$.

iii)  $a_\sigma\in \RR_{n-1}$ if and only if $a_{\sigma0}\in \RR_n$.  $b_\sigma\in \RR_{n-1}$ if and only if $b_{\sigma1}\in \RR_n$.

iv) If  $\sigma< 1^n$, then $\#(\{b_\sigma, a_{\sigma^+}\}\cap \mathcal R_{n-1})=1.$
\end{lem}

%
%

 We will prove this lemma in Sec. \ref{sec-Z-band-edge}.

\subsection{Two useful basic facts}\

We record two basic facts, which are frequently used in the whole paper.

The following lemma is fundamental for us.

\begin{lem}\label{fundamental}
For any $n\ge0$ and $E\in \R$,
\begin{equation}\label{fdit}
h_{n+1}(E)-(h_n^2(E)-2)=(-1)^n2\lambda\prod_{c\in\RR_n}(E-c)=(-1)^n2\lambda\prod_{j=0}^n h_j(E).
\end{equation}

\end{lem}

\proof
We prove it by induction.
For $n=0$, we have
$$
h_1(E)-(h_0^2(E)-2)=2\lambda(E-\lambda)=2\lambda h_0(E).
$$
Then \eqref{fdit} holds for $n=0$.

Assume \eqref{fdit} holds for $n=k\ge0$. Now consider $n=k+1$.
By \eqref{recurrence} and induction hypothesis,   we have
\begin{eqnarray*}
&&h_{k+1}(E)-(h_{k}^2(E)-2)\\
&=&h_{k}(E)(h_{k-1}^2(E)-2)-2-(h_{k}^2(E)-2)\\
&=&-h_{k}(E)\left(h_{k}(E)-(h_{k-1}^2(E)-2)\right)\\
&=&(-1)^k2\lambda\prod_{j=0}^k h_j(E)=(-1)^k2\lambda\prod_{c\in\RR_k}(E-c).
\end{eqnarray*}
Hence, the result holds for $n=k+1$.
By induction, the result follows.
\hfill $\Box$

The following is  standard Floquet theory:

\begin{prop}\label{basic-facts}
i) For each $\sigma\in \Sigma_n$, $h_n: B_\sigma\to [-2,2]$ is a homeomorphism. If $\sigma_n=0$ ($\sigma_n=1$), then
 $h_n'(E)<0$ ($h_n'(E)>0$) for $E\in {\rm int}(B_\sigma)$. Consequently
 $$
 h_n(a_\sigma)=2,\ \ h_n(b_\sigma)=-2 \ \ (h_n(a_\sigma)=-2,\ \ h_n(b_\sigma)=2).
 $$

 ii) If $\sigma\ne 1^n$, then the gap $(b_\sigma, a_{\sigma^+})$ is open if and only if
  $h_n'(b_\sigma)\ne0$ or $h_n'(a_{\sigma^+})\ne0.$

  iii) $E$ is an endpoint of some band of level $n$ iff $h_n(E)=\pm2.$ If $h_n(E)=-2$
  and $h_n'(E)<0$($h_n'(E)>0$), then
  $E=b_\sigma$ ($a_\sigma$) for some $\sigma\in \Sigma_n.$ If $h_n(E)=2$
  and $h_n'(E)<0$($h_n'(E)>0$), then
  $E=a_\sigma$ ($b_\sigma$) for some $\sigma\in \Sigma_n.$

\end{prop}


\section{Types, evolution law, order and labels}\label{sec-type-evo-order}

In this section, at first we study the band configurations of $\B_n\cup\B_{n+1}\cup\B_{n+2}$. Then we introduce the concept of type, which is fundamental for this paper. Next we derive the evolution law and the order of types and the labels of the directed edges. At last, we show that $\{\BB_n:n\ge0\}$ is a NS with limit set $\sigma_\lambda.$

\subsection{band configurations of $\B_n\cup\B_{n+1}\cup\B_{n+2}$}\

Now we present the second technical lemma of our paper, which summarizes all the possible band configurations of $\B_n\cup\B_{n+1}\cup\B_{n+2}$:

\begin{lem}\label{tech-2}
Fix $n\ge0$ and $\sigma\in \Sigma_n$.

i) Assume $n$ is odd. Then
 $ a_{\sigma10},b_{\sigma10}\notin \RR_{n+1}$
  and
  $$
  B_{\sigma00}\subset B_\sigma,\ B_{\sigma10}\subset {\rm int}(B_\sigma);\ \
  B_{\sigma01}\subset B_{\sigma0};\ \
  B_{\sigma0}<B_{\sigma10}\prec B_{\sigma1}.
  $$
  Moreover we have
  $$
 \begin{cases}
 B_{\sigma00}\subset B_{\sigma0},&\text{ if }\  a_\sigma\in \RR_{n-1}\\
  B_{\sigma00}\prec B_{\sigma0},&\text{ if }\  a_\sigma\notin\RR_{n-1}
 \end{cases}
 \ \ \ \ \
 \begin{cases}
 B_{\sigma11}\subset B_{\sigma1}\subset B_\sigma,&\text{ if }\  b_\sigma\in \RR_{n-1}\\
B_\sigma\prec B_{\sigma11},&\text{ if }\  b_\sigma\notin\RR_{n-1}
 \end{cases}
 $$
    If $a_{\sigma}\notin \RR_{n-1}$, then  $a_{\sigma00}, b_{\sigma00}\notin \RR_{n+1}$ and  $B_{\sigma00}\subset {\rm int}(B_\sigma)$.

   ii)  Assume $n$ is even. Then
 $ a_{\sigma01},b_{\sigma01}\notin \RR_{n+1}$
  and
  $$ B_{\sigma11}\subset B_\sigma,\
  B_{\sigma01}\subset {\rm int}(B_\sigma);\ \  B_{\sigma10}\subset B_{\sigma1};\ \
  B_{\sigma0}\prec B_{\sigma01}< B_{\sigma1}.
  $$
  Moreover we have
  $$
 \begin{cases}
 B_{\sigma11}\subset B_{\sigma1},&\text{ if }\  b_\sigma\in \RR_{n-1}\\
  B_{\sigma1}\prec B_{\sigma11},& \text{ if }\ b_\sigma\notin\RR_{n-1}
 \end{cases}
 \ \ \ \ \
 \begin{cases}
 B_{\sigma00}\subset B_{\sigma0}\subset B_\sigma,& \text{ if }\ a_\sigma\in \RR_{n-1}\\
B_{\sigma00}\prec B_{\sigma},& \text{ if }\ a_\sigma\notin\RR_{n-1}
 \end{cases}
 $$
    If $b_{\sigma}\notin \RR_{n-1}$, then  $a_{\sigma11}, b_{\sigma11}\notin \RR_{n+1}$ and  $B_{\sigma11}\subset {\rm int}(B_\sigma)$.
 \end{lem}

 Combine Lemmas \ref{tech-1} and \ref{tech-2}, we have

 \begin{cor}\label{tech-cor}

 i) Assume $n\ge1$ and $\sigma\in \Sigma_n.$ Then
 $$
 B_\sigma\subset B_{\sigma|_{n-1}}\Leftrightarrow B_\sigma\cap \RR_{n-1}\ne \emptyset\Leftrightarrow \partial B_\sigma\cap \RR_{n-1}\ne \emptyset.
 $$

 ii) Assume $n\ge1$ and $\sigma\in \Sigma_n$.  Then
\begin{equation*}\label{}
  B_\sigma\subset B_{\sigma|_{n-1}} \text{ or } \  B_\sigma\subset B_{\sigma|_{n-2}}.
\end{equation*}
(For $n=1$, $B_{\sigma|_0}:=B_\emptyset$ and $B_{\sigma|_{-1}}:=\R$).

iii) Assume $n\ge0$ and $\sigma\in \Sigma_n,\tau\in \Sigma_n\cup \Sigma_{n+1}\cup\Sigma_{n+2}$.  If $B_\tau\subset B_\sigma$, then $\sigma\lhd \tau.$

 \end{cor}

We will prove  Lemma \ref{tech-2} and Corollary \ref{tech-cor} in  Sec. \ref{sec-pf-lem3.1-cor3.2}. Now we apply them to derive the concept of type and figure out the evolution law.

\subsection{The types of the bands}\

Lemmas \ref{tech-1} and \ref{tech-2} suggest that for a band $B_\sigma$ of level $n$, once we know that whether $a_\sigma, b_\sigma\in \RR_{n-1}$, the local configurations of $B_\sigma$ and $B_{\sigma\tau}$ with $|\tau|\le2$ are completely determined. This observation motivates the following definition for types of bands.

\begin{defi}  \label{def-types}
Assume $n\ge0$,  $0\le \kappa\le 2$ and  $B\in\BB_n$.

At first we assume $B\in \mathcal{B}_n$.

We say $B$ has type $\kappa_o(\kappa_e)$, if $n$ is odd (even),  and
$\#\partial B\cap \RR_{n-1}=\kappa$.

(Recall that $\RR_{-1}=\emptyset.$)

Next we  assume $B\notin \B_n.$ Then $B\in \B_{n+1}$. There exists a unique
$\sigma\in \Sigma_{n+1}$ such that
 $B=B_\sigma$. Since $B\in \BB_n$, by Proposition \ref{BB-n}, $B\not\subset B_{\sigma|_n}$. By Corollary
 \ref{tech-cor} ii), $ B\subset B_{\sigma|_{n-1}}$.

 If $B_{\sigma|_n}\subset B_{\sigma|_{n-1}}$.

 We say $B$ has type $3_{ol} (3_{el})$ if $n$ is odd (even) and $B\prec B_{\sigma|_n}.$

  We say $B$ has type $3_{or} (3_{er})$ if $n$ is odd (even) and $ B_{\sigma|_n}\prec B.$

If $B_{\sigma|_n}\not\subset B_{\sigma|_{n-1}}$.

 We say $B$ has type $3_{ol} (3_{el})$ if $n$ is odd (even) and $ B_{\sigma|_n}\prec B_{\sigma|_{n-1}}.$

  We say $B$ has type $3_{or} (3_{er})$ if $n$ is odd (even) and $ B_{\sigma|_{n-1}}\prec B_{\sigma|_{n}}.$
\end{defi}

By the definition above, each band in $\BB_n$ has one and only one type.

\begin{exam}\label{exam-type}
{\rm
 Recall that  $\BB_0=\{B_0,B_\emptyset\}$. By the definition, $B_\emptyset $ has type $0_e$, $B_0$ has type $3_{el}$, see Figure \ref{fig-type} (a).
  $\BB_1=\{B_0,B_1,B_{01}, B_{11}\}$.  $B_0$ has type $0_o$,
 $B_1$ has  type $1_o$,  $B_{01}$ have type $3_{ol}$ and  $B_{11}$ have type $3_{or}$, see Figure \ref{fig-type} (a).
}
\end{exam}

\subsection{The evolution law  and the order of the types, the label of the directed edges }\

At first, we study type $3_\ast$ bands.

\begin{lem}\label{type-3}
 Assume $B_\sigma\in \BB_n$ has type $3_{ol}$ or $3_{or}$ ($3_{el}$ or $3_{er}$).
  Then $B_{\sigma}$ contains exact one band in $\BB_{n+1}$:
$B_\sigma$ itself with type $0_e$ ($0_o$).
\end{lem}

\proof\ We only show the odd case. Assume $B_\sigma\in \BB_n$ has type $3_{ol}$ or $3_{or}$,
 then $n$ is odd and $\sigma\in \Sigma_{n+1}$.  By Proposition \ref{BB-n},
 $B_\sigma\not\subset B_{\sigma|_n}$.
  By Corollary \ref{tech-cor} i),
  $ B_{\sigma}\cap \RR_n=\emptyset.$ Thus $B_\sigma\in \B_{n+1}\subset \BB_{n+1}$ has type $0_e$.
  By the definition of $\BB_{n+1},$ $B_\sigma$ is the only band in $\BB_{n+1}$ and
  contained in $B_\sigma.$
\hfill $\Box$

Before continuing, we do the following observation.
Take  $B_\sigma\in\B_n$ and $B_\tau\in \BB_{n+1}$. Assume
 $B_\tau\subset B_\sigma$. Then
  $\sigma\in \Sigma_n$,  $\tau\in \Sigma_{n+1}$ or $\Sigma_{n+2}$. By Corollary \ref{tech-cor} iii), $ \sigma\lhd \tau$. So
 \begin{equation}\label{tau}
 \tau\in \{\sigma0,\sigma1,\sigma00,\sigma01,\sigma10,\sigma11\}.
 \end{equation}

Next, we study type $0_\ast$  bands.

\begin{lem}\label{type-0}
Assume $B_\sigma\in \BB_n$.

i) If $B_\sigma$ has type $0_o$, then $B_{\sigma}$
contains exact three bands in $\BB_{n+1}$ with the order:
$$
B_{\sigma00}\prec B_{\sigma0}<B_{\sigma10},
$$
and with the types $3_{el}, 1_e, 3_{er}$, respectively.  Moreover,
$$b_{\sigma0}=z_\sigma;\ \ B_{\sigma00},B_{\sigma10}\subset{\rm int}(B_\sigma).$$

ii)  If $B_\sigma$ has type $0_e$, then $B_{\sigma}$ contains exact three bands in $\BB_{n+1}$ with the order:
$$
B_{\sigma01}< B_{\sigma1}\prec B_{\sigma11},
$$
and with the type $3_{ol}, 1_o, 3_{or}$, respectively.  Moreover,
$$a_{\sigma1}=z_\sigma;\ \ B_{\sigma01},B_{\sigma11}\subset{\rm int}(B_\sigma).$$
\end{lem}

\proof\
We only show i), so we assume $n$ is odd.
By the assumption, $\sigma\in \Sigma_n$ and
 $a_\sigma,b_\sigma\notin\RR_{n-1}$. Assume $B_\tau\in \BB_{n+1}$ and
  $B_\tau\subset B_\sigma$, then $\tau$ satisfies
\eqref{tau}.

 By Lemma \ref{tech-1} i) and Lemma \ref{endsymbol} iii),
 $$
 B_{\sigma0}\subset B_\sigma;\ \ B_{\sigma}\prec B_{\sigma1};\ \ a_{\sigma0}\notin\RR_n;\ \  b_{\sigma0}=z_\sigma\in \RR_n.
 $$
 So $B_{\sigma0}\in \BB_{n+1}$ has type $1_e$ and $B_{\sigma1}\not\subset B_\sigma$.

  By  Lemma \ref{tech-2} i),
  $$
  \begin{cases}
  B_{\sigma00}, B_{\sigma10}\subset {\rm int} (B_\sigma);\ \
   B_{\sigma00}\prec B_{\sigma0}<B_{\sigma10}\prec B_{\sigma1},\\
  B_{\sigma01}\subset B_{\sigma0};\ \  B_\sigma\prec B_{\sigma11}.
  \end{cases}$$
   So $B_{\sigma11}\not\subset B_\sigma$ and $B_{\sigma01}\notin \BB_{n+1}$;
    $B_{\sigma00}\not\subset B_{\sigma0}, B_{\sigma10}\not\subset B_{\sigma1}$.
   By Proposition \ref{BB-n}, $B_{\sigma00}, B_{\sigma10}\in \BB_{n+1}$.
    Since
   $
    B_{\sigma0}\subset B_\sigma
   $
   and $B_{\sigma00}\prec B_{\sigma0}$, by the definition, $B_{\sigma00}$ has type $3_{el}$.
   Since
   $
    B_{\sigma1}\not\subset B_\sigma
   $
   and $B_{\sigma}\prec B_{\sigma1}$, by the definition, $B_{\sigma10}$ has type $3_{er}$.
\hfill $\Box$

As two examples,
 $B_\emptyset$ has type $0_e$ and
 $B_0$ has type $0_o$,  see Figure \ref{fig-type} (a) and (b) for the order configurations of their son intervals.

Now, we study type $1_\ast$  bands.

\begin{lem}\label{type-1}
Assume $B_\sigma\in \BB_n$.

i) If $B_\sigma$ has type $1_o$, then
 $B_{\sigma}$ contains exact two bands in $\BB_{n+1}$ with the  order:
$$
 B_{\sigma0}<B_{\sigma10},
$$
and with the type $ 2_e, 3_{er}$, respectively. Moreover $B_{\sigma0}=[a_\sigma,z_\sigma]$ and  $B_{\sigma10}\subset {\rm int}(B_\sigma).$

ii)  If $B_\sigma$ has type $1_e$, then
$B_{\sigma}$ contains exact two bands in $\BB_{n+1}$ with the  order:
$$
B_{\sigma01}< B_{\sigma1},
$$
and with the type $3_{ol}, 2_o$, respectively. Moreover $B_{\sigma1}=[z_\sigma,b_\sigma]$ and  $B_{\sigma01}\subset {\rm int}(B_\sigma).$
\end{lem}

\proof\
We only show i), so  we assume $n$ is odd.
By the assumption, $\sigma\in \Sigma_n$ and
 $\#\{a_\sigma,b_\sigma\}\cap\RR_{n-1}=1$. Assume $B_\tau\in \BB_{n+1}$ and
  $B_\tau\subset B_\sigma$, then $\tau$ satisfies
\eqref{tau}.

At first, we claim that $a_\sigma\in \RR_{n-1}$ and $b_\sigma\notin\RR_{n-1}.$
   Indeed, since $\partial B_\sigma\cap \RR_{n-1}\ne\emptyset, $ by Corollary \ref{tech-cor} i), $B_\sigma\subset B_{\sigma|_{n-1}}.$  Write $\hat\sigma=\sigma|_{n-1}$, notice that $n-1$ is even.
    By applying Lemma \ref{tech-1} ii) to $\hat\sigma$, we have $b_{\hat\sigma0}\notin\RR_{n-1}$ and
    $a_{\hat\sigma1}=z_{\hat\sigma}\in \RR_{n-1}.$ Since $\#\{a_\sigma,b_\sigma\}\cap\RR_{n-1}=1$, we have
    $a_\sigma\in \RR_{n-1}$ and $b_\sigma\notin\RR_{n-1}$ in either case of $\sigma=\hat\sigma0$ or $\sigma=\hat\sigma1.$ So the claim holds.

  By Lemma \ref{tech-1} i), $B_{\sigma0}=[a_\sigma,z_\sigma]\subset B_\sigma$ and
  $ B_\sigma\prec B_{\sigma1}.$ Since
  $a_{\sigma}\in\RR_{n-1}$ and $z_\sigma\in \RR_n$,  $B_{\sigma0}\in \BB_{n+1}$ has type $2_e.$
  $B_{\sigma1}\not\subset B_\sigma.$

  By  Lemma \ref{tech-2} i),
  $$
B_{\sigma00}, B_{\sigma01}\subset B_{\sigma0};\ \
  B_{\sigma10}\subset {\rm int}(B_{\sigma}),\ \  B_{\sigma10}\prec B_{\sigma1};\ \
   B_{\sigma}\prec B_{\sigma11}, \ \ B_{\sigma0}<B_{\sigma10}.
  $$
  So $B_{\sigma11}\not\subset B_\sigma $ and
  $B_{\sigma00},B_{\sigma01}\notin \BB_{n+1}$. Since $B_{\sigma10}\not\subset B_{\sigma1}$,
   by Proposition \ref{BB-n}, $B_{\sigma10}\in \BB_{n+1}$.
  Since  $B_\sigma\prec B_{\sigma1}$,
   $B_{\sigma10}$ has type $3_{er}.$ So the result holds.
\hfill $\Box$

As two examples,
 $B_{00}$ has type $1_e$ and
 $B_1$ has type $1_o$,  see Figure \ref{fig-type} (c) and (d) for the order configurations of their son intervals.

At last, we study type $2_\ast$  bands.

\begin{lem}\label{type-2}
Assume $B_\sigma\in \BB_n$.

i) If $B_\sigma$ has type $2_o,$ then $B_{\sigma}$ contains exact three bands in $\BB_{n+1}$ with the order:
$$
 B_{\sigma0}<B_{\sigma10}\prec B_{\sigma1},
$$
and with the type $2_e, 3_{el}, 1_e$, respectively.

ii) If $B_\sigma$ has type $2_e,$ then $B_{\sigma}$ contains exact three bands in $\BB_{n+1}$ with the order:
$$
B_{\sigma0}\prec B_{\sigma01}< B_{\sigma1},
$$
and with the type $1_o, 3_{or}, 2_o$, respectively.
\end{lem}

\proof\
We only show i), so we assume $n$ is odd.
By the assumption, $\sigma\in \Sigma_n$ and
 $a_\sigma,b_\sigma\in\RR_{n-1}$. By Lemma \ref{endsymbol} ii), $b_{\sigma1}=b_\sigma.$ Assume $B_\tau\in \BB_{n+1}$ and  $B_\tau\subset B_\sigma$, then $\tau$ satisfies
\eqref{tau}.

 By Lemma \ref{tech-1} i), $B_{\sigma0}=[a_\sigma,z_\sigma],B_{\sigma1}\subset B_\sigma$,
  $a_{\sigma1}\notin\RR_n$ and $b_{\sigma1}=b_\sigma\in \RR_n$. So
  $B_{\sigma1}\in \BB_{n+1}$ has type $1_e.$ Since $a_\sigma\in \RR_{n-1}$ and $z_\sigma\in \RR_n$,
   $B_{\sigma0}\in \BB_{n+1}$ has type $2_e$.

   By  Lemma \ref{tech-2} i),
  $$
B_{\sigma00}, B_{\sigma01}\subset B_{\sigma0};\ \
  B_{\sigma10}\subset B_{\sigma};\ \
  B_{\sigma11}\subset B_{\sigma1}, \ \ B_{\sigma0}<B_{\sigma10}\prec B_{\sigma1}.
  $$
  So $B_{\sigma00},B_{\sigma01}, B_{\sigma11}\notin \BB_{n+1}$.
  Since $B_{\sigma10}\not\subset B_{\sigma1}$,
   by Proposition \ref{BB-n}, $B_{\sigma10}\in \BB_{n+1}$.
    Since
   $
    B_{\sigma1}\subset B_\sigma
   $
   and $B_{\sigma10}\prec B_{\sigma1}$, by the definition, $B_{\sigma10}$ has type $3_{el}$.
\hfill $\Box$

As two examples,
 $B_{10}$ has type $2_e$ and
 $B_{101}$ has type $2_o$,  see Figure \ref{fig-type} (e) and (f) for the order configurations of  their son intervals.

We summarize Lemmas \ref{type-3}-\ref{type-2} in a symbolic way as follows:

\begin{eqnarray}
\label{evo-3}&& 3_{ol}, 3_{or}\stackrel{\emptyset}{\longrightarrow} 0_e;\ \ \ 3_{el}, 3_{er}\stackrel{\emptyset}{\longrightarrow} 0_o\\
\label{evo-0}&&
\begin{cases}
  0_o\stackrel{00}{\longrightarrow} 3_{el}; \ 0_o\stackrel{0}{\longrightarrow} 1_{e};\  0_o\stackrel{10}{\longrightarrow} 3_{er}; \ \ \  3_{el}\prec 1_e<3_{er};\\
  0_e\stackrel{01}{\longrightarrow} 3_{ol}; \ 0_e\stackrel{1}{\longrightarrow} 1_{o};\  0_e\stackrel{11}{\longrightarrow} 3_{or}; \ \ \  3_{ol}<1_o\prec 3_{or}.
  \end{cases}\\
  \label{evo-1}&&
\begin{cases}
  1_o\stackrel{0}{\longrightarrow}  2_{e};\  1_o\stackrel{10}{\longrightarrow} 3_{er}; \ \ \   2_e<3_{er};\\
  1_e\stackrel{01}{\longrightarrow} 3_{ol}; \ 1_e\stackrel{1}{\longrightarrow} 2_{o}; \ \ \  3_{ol}<2_o.
  \end{cases}\\
\label{evo-2}
&&\begin{cases}
  2_o\stackrel{0}{\longrightarrow} 2_{e}; \ 2_o\stackrel{10}{\longrightarrow} 3_{el};\  2_o\stackrel{1}{\longrightarrow} 1_{e}; \ \ \  2_{e}<3_{el}\prec 1_{e};\\
  2_e\stackrel{0}{\longrightarrow} 1_{o}; \ 2_e\stackrel{01}{\longrightarrow} 3_{or};\  2_e\stackrel{1}{\longrightarrow} 2_{o}; \ \ \  1_o\prec3_{or}<2_{o}.
  \end{cases}
\end{eqnarray}
Here $\alpha\stackrel{\tau}{\longrightarrow} \beta$ means type $\alpha$ can evolve to type $\beta$ and the directed edge $\alpha\beta$ is labeled by $\tau$.

Finally, consider $\BB_0=\{B_0,B_\emptyset\}$. By \eqref{initial}, $B_0\prec B_\emptyset$. By Example \ref{exam-type}, $B_0$ has type $3_{el}$ and  $B_\emptyset $ has type $0_e$. So we complete the order of the types by
\begin{equation}\label{order-addi}
  3_{el}\prec 0_e.
\end{equation}

\subsection{$\{\BB_n:n\ge0\}$ as a NS}\

 At first we have

\begin{prop} \label{unique-father}
 For each $\tilde B\in \BB_{n+1}$, there is a unique band $B\in \BB_n$ such that $\tilde B\subset B.$
\end{prop}

\proof\
 Let $\sigma\in\Sigma_{n+1}\cup\Sigma_{n+2}$ be the unique word such that
  $\tilde B=B_\sigma.$ If $B\in\BB_n$ is such that
 $\tilde B\subset B$, then by Corollary \ref{tech-cor} iii), the only two possibilities are
 $B=B_{\sigma|_n}$ or $B_{\sigma|_{n+1}}.$

 If $\tilde B$ has type $3_\ast,$ then $\sigma\in \Sigma_{n+2}$. By Proposition \ref{BB-n},
  $B_\sigma\not\subset B_{\sigma|_{n+1}}.$ By Corollary \ref{tech-cor} ii),
   $B_\sigma\subset B_{\sigma|_n}\in \B_n\subset \BB_n.$

   If $\tilde B$ has type $0_*$, then $\sigma\in \Sigma_{n+1}$ and
   $\partial B_\sigma\cap \RR_n=\emptyset$.  By Corollary \ref{tech-cor} i),
  $B_\sigma\not\subset B_{\sigma|_{n}}.$ By Proposition \ref{BB-n},
  $B_\sigma\in \BB_n$ and $\tilde B= B_\sigma.$

   If $\tilde B$ has type $1_\ast$ or $2_\ast$, then $\sigma\in \Sigma_{n+1}$ and
   $\partial B_\sigma\cap \RR_n\ne\emptyset$.  By Corollary \ref{tech-cor} i),
  $B_\sigma\subset B_{\sigma|_{n}}\in \BB_n.$ By the definition of $\BB_n$,
  $B_\sigma\notin \BB_n$.

  In all the cases, there is a unique $B\in \BB_n$ such that $\tilde B\subset B.$
\hfill $\Box$

See Sec. \ref{sec-NS} for the definition of NS. Now we can show that

\begin{cor}\label{NS-BB}
$\BB:=\{\BB_n:n\ge0\}$ is a NS and $A(\BB)=\sigma_\lambda$.
\end{cor}

\proof
By Proposition \ref{gap-open-sigma-n}, $\B_n$ is a disjoint family. Combine with  the definition \eqref{def-BB-n} of $\BB_n$, $\BB$ is optimal.
By Proposition \ref{unique-father}, $\BB$ is nested.
By Lemmas \ref{type-3}-\ref{type-2},
each band in $\BB_n$ contains at least one band in $\BB_{n+1}$, so $\BB$ is minimal. By the definition of NS (see Sec. \ref{sec-NS}), $\BB$ is a NS.
By \eqref{sigma-lamb-n} and the definition of $\BB_n$,
$$
\bigcup_{B\in \BB_n}B=\sigma_{\lambda,n}\cup\sigma_{\lambda,n+1}.
$$
By \eqref{spec-covering}, $\sigma_\lambda$ is the limit set of $\BB$.
\hfill $\Box$

As a final complement, we show the following

\begin{lem}\label{0-momotone}
If $B\in \BB_n$ has type $0_e(0_o)$, then $h_n$ is increasing (decreasing) on $B$.
\end{lem}

\proof\
We only show the case that $B$ has type $0_e$. The proof of the other case is the same.

If $n=0$, then $B=B_\emptyset$ and $h_0(E)=E-\lambda$, so the result holds.

Now assume $n\ge1$ and $n$ is even. Assume $B=B_\sigma$ for some $\sigma\in \Sigma_n$. We claim that
$\sigma_n=1$. If otherwise $\sigma=\tau0$ for some $\tau\in \Sigma_{n-1}$. Then by Lemma \ref{tech-1}
i), we have $b_{\sigma}=z_{\tau}\in \RR_{n-1}$, which contradicts with the fact that
 $B$ has type $0_e$. Now by
Proposition \ref{basic-facts} i), $h_n$ is increasing on $B.$
\hfill $\Box$


\section{Symbolic space $\Omega_\infty$}\label{sec-symbolic}

In this section, we define two orders on  the symbolic space $\Omega_\infty$ and characterize the gaps of it by using the map $\Pi$.

We have defined the symbolic space $\Omega_\infty$ in Sec. \ref{sec-main-coding}. Through Sec. \ref{sec-type-evo-order},  we understand that $\A$, defined by \eqref{alphabet},  is just the set of all the possible types; \eqref{adm} is just a summary of  the evolution laws appeared in  \eqref{evo-3}-\eqref{evo-2}. These information tell us how to construct the directed graph $\mathbb G.$ The labels of edges of  $\mathbb G$ can also be read off from \eqref{evo-3}-\eqref{evo-2}. Now the order relations appeared in \eqref{evo-0}-\eqref{order-addi} suggest the following definitions of orders on $\Omega_\infty.$

\subsection{Two orders on $\Omega_\infty$}\

Define the set of {\it admissible words} as
 \begin{equation}\label{Omega-n-ast}
 \Omega_n:=\{w\in \A^{n+1}: w_0\in \{3_{el},0_e\}, a_{w_jw_{j+1}}=1\},
  \ n\ge 0\ \ \text{ and }\ \ \ \Omega_\ast:=\bigcup_{n\ge0} \Omega_n.
  \end{equation}
  For $\omega=\omega_0\omega_1\cdots\in \Omega_\ast\cup\Omega_\infty$, write $\omega|_n:=\omega_0\cdots \omega_n$. For $\omega,\omega'\in \Omega_\ast\cup\Omega_\infty$, write $\omega\wedge \omega'$ for their maximal common prefix.
  For $w\in \Omega_n$, define the cylinder as
$$
[w]:=\{\omega\in \Omega_\infty: \omega|_n=w\}.
$$

\subsubsection{Two orders on $\A$}\

Define $\preceq$ to be the smallest   partial order  on $\A$ such that the following holds:
\begin{equation}\label{order}
3_{ol}\prec 1_{o}\prec 3_{or}\prec 2_o;\ \  2_e\prec 3_{el}\prec 1_{e}\prec 3_{er};\ \  3_{el}\prec 0_e.
\end{equation}
Define $\le $ to be the smallest partial order  on $\A$ such that the following holds:
\begin{equation}\label{order-strong}
3_{ol}< 1_{o}< 2_o;\ \ 3_{ol}< 3_{or}< 2_o;\ \  2_e< 3_{el}< 3_{er};\ \
 2_e< 1_{e}< 3_{er}.
\end{equation}
(Here we view partial order as a subset of $\A\times \A$.)

By \eqref{order} and \eqref{order-strong}, it
 is direct to check that if $\alpha<\beta$, then $\alpha\prec\beta.$

 \begin{rem}
{\rm
The idea for  defining \eqref{order} and \eqref{order-strong} is as follows. At first, we impose
\begin{equation*}\label{order-idea}
3_{ol}<1_{o}\prec 3_{or}<2_o;\ \  2_e<3_{el}\prec 1_{e}<3_{er};\ \  3_{el}\prec 0_e,
\end{equation*}
which is a summary of  the order relations appeared in \eqref{evo-0}-\eqref{order-addi}.
Since $<$ is stronger than $\prec$, by transitivity, we get \eqref{order} and \eqref{order-strong}.
}
\end{rem}

\subsubsection{Induced  orders on $\Omega_\infty$}\label{sec-order-Omega}\

Define a relation $\preceq(\le) $ on $\Omega_\infty$ as follows. For any $\omega\in \Omega_\infty,$ let
$\omega\preceq \omega(\omega\le \omega).$ Given $\omega,\hat\omega\in \Omega_\infty$
 and $\omega\ne\hat\omega.$ Assume
$n$ is such that $\omega|_{n-1}=\hat\omega|_{n-1}$ and $\omega_n\ne \hat\omega_n$. Let
$\omega\prec \hat\omega(\omega<\hat\omega)$ if $\omega_n\prec \hat\omega_n(\omega_n<\hat\omega_n).$

\begin{lem}\label{two-orderings}
Both $\preceq$ and $\le$ are partial orders on $\Omega_\infty$.
Moreover, $\preceq$ is a total order. If $\omega<\hat\omega$, then $\omega\prec \hat\omega.$
\end{lem}

\proof\
It is routine to check that both $\preceq$ and $\le$ are partial orders.

Now we show that $\preceq$ is a total order.
Assume $\omega,\hat\omega\in \Omega_\infty$ and $\omega\ne\hat\omega.$
Let $n\ge0$ be the minimal integer such that $\omega_n\ne \hat\omega_n$.
If $n=0$, then $\omega_0\ne\hat \omega_0$. By \eqref{def-omega-infty}, $\{\omega_0,\hat\omega_0\}=\{3_{el},0_e\}$.
By \eqref{order}, we can compare $\omega_0$ and $\hat\omega_0$.
Then by the definition, we can compare $\omega$ and $\hat\omega.$
If $n\ge1$, then $\omega\wedge\hat\omega=\omega|_{n-1}=\hat\omega|_{n-1}$ and by \eqref{adm},
 $\omega_{n-1}=\hat\omega_{n-1}\in\{0_e,0_o,1_e,1_o,2_e,2_o\}$. Now by
 \eqref{adm} and \eqref{order}, we can compare $\omega_n$ and $\hat\omega_n$.
  So, we can compare
 $\omega$ and $\hat\omega.$

 Since  in $\A$, $\alpha<\beta$ implies  $\alpha\prec\beta$, by the definitions, the last statement holds.
\hfill $\Box$

Take $\omega,\hat\omega\in\Omega_\infty$. Assume  $\omega\prec \hat\omega$. Define the {\it open interval}
$(\omega,\hat\omega)$ as
$$
(\omega,\hat\omega):=\{\tilde \omega\in \Omega_\infty: \omega\prec \tilde\omega\prec
 \hat \omega\}.
$$
We call $(\omega,\hat\omega)$ a {\it gap} of $\Omega_\infty$
if $(\omega,\hat\omega)=\emptyset$. We also say that $\omega$ (
$\hat\omega$) is the {\it left (right) edge} of the gap.

\subsection{Characterization of the gaps of $\Omega_\infty$}\label{sec-gap-symbolic}\

We denote the set of the gaps of $\Omega_\infty$ by $\G.$ The following result is a symbolic version of Theorem \ref{main-gap-zero-inftyenergy}:

\begin{thm}\label{symbolic-gap}
For $\Omega_\infty$, we have

i) $\omega_\ast\preceq \omega \preceq\omega^\ast$ for any $\omega\in \Omega_\infty$.

ii) (Type-I gaps): There exists a bijection $ \ell^o: \mathcal E_{l}^o\to \mathcal E_{r}^o\setminus \{\omega_\ast\}$   such that $\omega<\ell^o(\omega)$ for any $\omega\in \mathcal E_{l}^o$ and
\begin{equation}\label{gap-Io}
  \G_I^o:=\{(\omega,\ell^o(\omega)): \omega\in \mathcal E_{l}^o\}\subset \G.
\end{equation}
There exists a bijection $\ell^e: \mathcal E_{r}^e\to \mathcal E_{l}^e\setminus \{\omega^\ast\}$   such that $\ell^e(\omega)<\omega$ for any $\omega\in \mathcal E_{r}^e$ and
\begin{equation}\label{gap-Ie}
  \G_I^e:=\{(\ell^e(\omega),\omega): \omega\in \mathcal E_{r}^e\}\subset \G.
  \end{equation}

  iii) (Type-II gaps): There exists a bijection $\ell: \widetilde {\mathcal E}_l\to \widetilde {\mathcal E}_r$   such that $\omega\prec\ell(\omega)$ for any $\omega\in\widetilde {\mathcal E}_l$ and
\begin{equation}\label{gap-II}
  \G_{II}:=\{(\omega,\ell(\omega)): \omega\in \widetilde {\mathcal E}_l\}\subset \G.
\end{equation}

 iv)  $\G=\G_I\cup\G_{II}$, where $\G_I:=\G_I^o\cup\G_{I}^e.$
\end{thm}

We will prove Theorem \ref{symbolic-gap} by using  the map $\Pi$ defined by \eqref{def-Pi}. But at first, we need to know how to define $\ell, \ell^o$ and $\ell^e.$

\begin{figure}[htbp]

\begin{minipage}[t]{0.3\linewidth}
\centering
\includegraphics[width=1.88in]{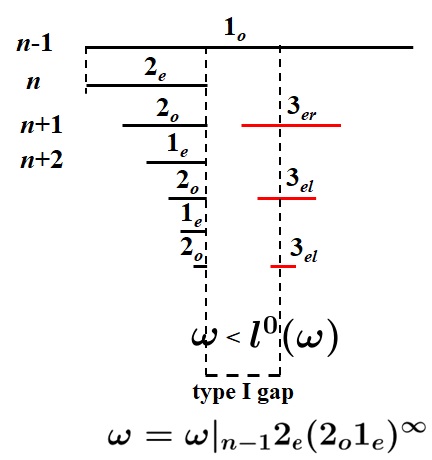}
\caption*{(a)}
\end{minipage}%
\hskip 3cm
\begin{minipage}[t]{0.3\linewidth}
\centering
\includegraphics[width=2in]{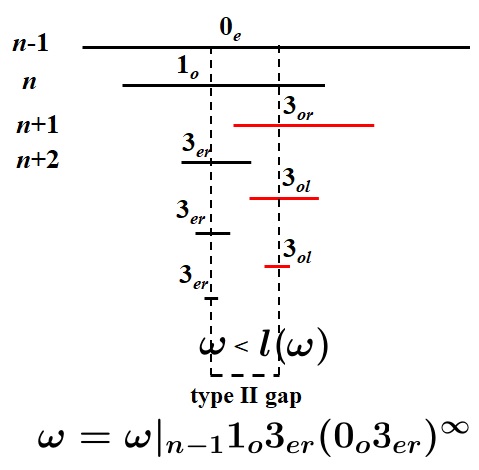}
\caption*{(b)}
\end{minipage}%
\caption{Two types of gaps}\label{2-type-gap}
\end{figure}

\subsubsection{Definitions and bijectivities  of $\ell, \ell^o, \ell^e$}\

At first, we define $\ell:\widetilde{\E_l}\to \widetilde{\E_r}$. Assume $\omega\in \widetilde{\E_l}$. Then there exists $n$ such that $\omega=\omega|_n3_{er}(0_o3_{er})^\infty$ with $\omega_n\ne 0_o$; or $\omega=\omega|_n(0_o3_{er})^\infty$ with $\omega_n\ne3_{er}.$  By tracing the graph $\mathbb G$ backward, in the former case, $ \omega_n=1_o$ and hence $n\ge1, \omega_{n-1}=0_e$ or $2_e$; in the latter case, $\omega_n=3_{el}$ and hence  $n\ge2, \omega_{n-1}=0_o$ or $2_o$, or $n=0$. Now define
\begin{equation}\label{def-elll}
  \ell(\omega):=
  \begin{cases}
   \omega|_{n-1}3_{or}(0_e3_{ol})^\infty,& \text{ if }\  \omega=\omega|_{n-1}1_o3_{er}(0_o3_{er})^\infty, n\ge1\\
   \omega|_{n-1}1_e3_{ol}(0_e3_{ol})^\infty,& \text{ if }\  \omega=\omega|_{n-1}3_{el}(0_o3_{er})^\infty, n\ge2\\
   (0_e3_{ol})^\infty,& \text{ if }\  \omega=3_{el}(0_o3_{er})^\infty.
   \end{cases}
\end{equation}
See Figure \ref{2-type-gap} (b) for an illustration  of the definition of $\ell$ and for some intuition of a gap of $\infty_{II}$-type.

Next, we define $ \ell^o: \mathcal E_{l}^o\to \mathcal E_{r}^o\setminus \{\omega_\ast\}$. Assume $\omega\in \E^o_l$. Then  there exists $n$ such that either $\omega=\omega|_n(2_o1_e)^\infty$ with $\omega_n\ne 1_e$; or $\omega=\omega|_n1_e(2_o1_e)^\infty$ with $\omega_n\ne 2_o$. By tracing the graph $\mathbb G$ backward, in the former case, $\omega_n=2_e$ and hence $n\ge2, \omega_{n-1}=1_o$ or $2_o$; in the latter case, $\omega_n=0_o$ and hence $n\ge1$. Now define
\begin{equation}\label{def-ell-o}
  \ell^o(\omega):=
  \begin{cases}
   \omega|_{n-2}1_o3_{er}(0_o3_{el})^\infty,& \text{ if }\  \omega=\omega|_{n-2}1_o2_e(2_o1_e)^\infty, n\ge2\\
   \omega|_{n-2}2_o3_{el}(0_o3_{el})^\infty,& \text{ if }\  \omega=\omega|_{n-2}2_o2_e(2_o1_e)^\infty,n\ge2\\
   \omega|_{n-1}0_o3_{er}(0_o3_{el})^\infty,& \text{ if }\  \omega=\omega|_{n-1}0_o1_e(2_o1_e)^\infty, n\ge1
   \end{cases}
\end{equation}
See Figure \ref{2-type-gap} (a) for an illustration  of the definition of $\ell^o$ and for some intuition of a gap of $\infty_{I}$-type.

Now, we define $ \ell^e: \mathcal E_{r}^e\to \mathcal E_{l}^e\setminus \{\omega^\ast\}$. Assume $\omega\in \E^e_r$. Then  there exists $n$ such that either $\omega=\omega|_n(2_e1_o)^\infty$ with $\omega_n\ne 1_o$; or $\omega=\omega|_n1_o(2_e1_o)^\infty$ with $\omega_n\ne 2_e$. By tracing the graph $\mathbb G$ backward, in the former case, $\omega_n=2_o$ and hence $n\ge3, \omega_{n-1}=1_e$ or $2_e$; in the latter case, $\omega_n=0_e$ and hence $n\ge0$. Now define
\begin{equation}\label{def-ell-e}
  \ell^e(\omega):=
  \begin{cases}
   \omega|_{n-2}1_e3_{ol}(0_e3_{or})^\infty,& \text{ if }\  \omega=\omega|_{n-2}1_e2_o(2_e1_o)^\infty, n\ge3\\
   \omega|_{n-2}2_e3_{or}(0_e3_{or})^\infty,& \text{ if }\  \omega=\omega|_{n-2}2_e2_o(2_e1_o)^\infty,n\ge3\\
   \omega|_{n-1}0_e3_{ol}(0_e3_{or})^\infty,& \text{ if }\  \omega=\omega|_{n-1}0_e1_o(2_e1_o)^\infty, n\ge0
   \end{cases}
\end{equation}

\begin{lem}\label{bij-ell-o-e}
i) $\ell:\widetilde{\E_l}\to \widetilde{\E_r}$ is a bijection and $\omega\prec\ell(\omega)$ for any $\omega\in\widetilde {\mathcal E}_l$.

ii) $ \ell^o: \mathcal E_{l}^o\to \mathcal E_{r}^o\setminus \{\omega_\ast\}$ is a bijection and $\omega<\ell^o(\omega)$ for any $\omega\in \mathcal E_{l}^o$.

iii) $\ell^e: \mathcal E_{r}^e\to \mathcal E_{l}^e\setminus \{\omega^\ast\}$  is a bijection and  $\ell^e(\omega)<\omega$ for any $\omega\in \mathcal E_{r}^e$.
\end{lem}

\proof\
i)
By \eqref{def-elll} and \eqref{order}, we have $\omega\prec\ell(\omega)$.

To show the bijection of $\ell$, we simply reverse the definition of $\ell$ and  construct another map $\eta:\widetilde{\E_r}\to \widetilde{\E_l}$, then check $\eta$ is the inverse of $\ell$. Assume $\omega\in \widetilde{\E_r}$. Then   either there exists $n\ge1$ such that  $\omega=\omega|_n(0_e3_{ol})^\infty$ with $\omega_n\ne3_{ol}$; or there exists $n\ge0$ such that $\omega=\omega|_n3_{ol}(0_e3_{ol})^\infty$ with $\omega_n\ne 0_e$; or $\omega=(0_e3_{ol})^\infty.$  By tracing the graph $\mathbb G$ backward,  in the first case, $\omega_n=3_{or}$ and hence  $\omega_{n-1}=0_e$ or $2_e$. In the second case, $ \omega_n=1_e$ and hence $n\ge2, \omega_{n-1}=0_o$ or $2_o$, Now define
\begin{equation*}\label{def-ell}
  \eta(\omega):=
  \begin{cases}
    \omega|_{n-1}1_o3_{er}(0_o3_{er})^\infty, & \text{ if }\  \omega=\omega|_{n-1}3_{or}(0_e3_{ol})^\infty, n\ge1\\
   \omega|_{n-1}3_{el}(0_o3_{er})^\infty,& \text{ if }\  \omega=\omega|_{n-1}1_e3_{ol}(0_e3_{ol})^\infty, n\ge2\\
   3_{el}(0_o3_{er})^\infty,&  \text{ if }\ \omega=(0_e3_{ol})^\infty.
   \end{cases}
\end{equation*}
One  check directly that $\eta\circ\ell={\rm Id}_{\widetilde{\E_l}}$ and $\ell\circ\eta={\rm Id}_{\widetilde{\E_r}}$. So $\ell$ is a bijection.

ii) and iii) can be proven similarly, we omit the proof.
\hfill $\Box$
\subsubsection{The property of the  map $\Pi: \Omega_\infty\to \Sigma_\infty$}\label{sec-pro-Pi}\

 Recall that $\Pi$ is defined by \eqref{def-Pi}. We treat the map $\Pi$ more precisely as follows.

Define the {\it label-assigning map} $\LL: \{\alpha\beta: \alpha\to\beta\}\to \{\emptyset,0,1,00,01,10,11\}$ as follows: $\LL(\alpha\beta)$ is the label of the edge $\alpha\beta$.

 By \eqref{evo-0}-\eqref{order-addi}, the following holds:
\begin{equation}\label{order-label}
 \alpha\to \beta,\beta' \text{ and }\  \beta<\beta'\Rightarrow \LL(\alpha\beta)=0\ast; \ \ \LL(\alpha\beta')=1\ast.
\end{equation}

Define a map $\Pi_\ast:\Omega_\ast\to \Sigma_\ast$ as
\begin{equation}\label{Pi-star}
  \Pi_\ast(w):=i(w)\LL(w_0w_1)\LL(w_1w_2)\cdots \LL(w_{n-1}w_n), \ \  w\in \Omega_n,
\end{equation}
where $i(w)=\emptyset$ if $w_0=0_e$ and $i(w)=0$ if $w_0=3_{el}$.

By \eqref{evo-3}-\eqref{evo-2} and
\eqref{def-Pi}, it is direct to check that
\begin{equation}\label{Pi-L}
 \Pi_\ast(\omega|_n) \lhd \Pi(\omega);\ \ |\Pi_\ast(\omega|_n)|=n \text{ or } n+1;\ \ \Pi(\omega)=\lim_{n\to\infty}\Pi_\ast(\omega|_n).
\end{equation}

Define
\begin{equation}\label{Simga-ast-oe}
  \Sigma_\ast^o:=\bigcup_{n\ge0}\Sigma_{2n+1}; \ \ \Sigma_\ast^e:=\bigcup_{n\ge0}\Sigma_{2n};\ \ \Sigma_\infty^{(2)}:=\{\sigma(01)^\infty: \sigma\in \Sigma_\ast\}.
\end{equation}

 We endow $\Sigma_\infty$ with the lexicographical order $\le$.

\begin{prop}\label{Pi-basic}

i) For any $\omega\in \widetilde{\E_l}$, we have $\Pi(\omega)=\Pi(\ell(\omega))$. Consequently \begin{equation}\label{equal-E-l-r}
\Pi(\widetilde{\E_l})=\Pi(\widetilde{\E_r})\subset \Sigma_\infty^{(2)}.
\end{equation}

ii) Assume $\sigma\in \Sigma_\infty.$ If $\#\Pi^{-1}(\{\sigma\})\ge2,$  then there exists $\omega\in \widetilde{\E_l}$ such that
$$\Pi^{-1}(\{\sigma\})=\{\omega,\ell(\omega)\}.$$
Consequently, $\#\Pi^{-1}(\{\sigma\})\ge2$  if and only if $\sigma\in \Pi(\widetilde{\E_l}).$

iii) $\Pi$ is surjective and order-preserving.
\begin{equation}\label{2-coding-bdy}
  \Pi(\omega_\ast)=0^\infty;\ \ \Pi(\omega^\ast)=1^\infty.
\end{equation}
Moreover, $\Pi: \Omega_\infty\setminus  \widetilde{\E_r}\to \Sigma_\infty$ is bijective.

iv)
If $\omega\in \E^o_l$, then there exists a unique $\sigma\in \Sigma_\ast^o$ such that
\begin{equation}\label{2-coding-Io}
  \Pi(\omega)=\sigma01^\infty;\ \ \Pi(\ell^o(\omega))=\sigma10^\infty.
\end{equation}

If $\omega\in \E^e_r$, then there exists a unique $\sigma\in \Sigma_\ast^e$ such that
\begin{equation}\label{2-coding-Ie}
  \Pi(\omega)=\sigma10^\infty;\ \ \Pi(\ell^e(\omega))=\sigma01^\infty.
\end{equation}

The following restrictions of $\Pi$ are bijections:
\begin{eqnarray}\label{bij-I-type}
\begin{cases}
  \Pi_{ol}:\E^o_l\to \{\sigma01^\infty: \sigma\in \Sigma^o_\ast \};\ \ \ \Pi_{or}:\E^o_r\setminus \{\omega_\ast\}\to \{\sigma10^\infty: \sigma\in \Sigma^o_\ast \}\\
  \Pi_{er}:\E^e_r\to \{\sigma10^\infty: \sigma\in \Sigma^e_\ast \};\ \ \ \Pi_{el}:\E^e_l\setminus \{\omega^\ast\}\to \{\sigma01^\infty: \sigma\in \Sigma^e_\ast \}
  \end{cases}
\end{eqnarray}

v)  $\Pi: \widetilde {\mathcal E}_{l}\cup \mathcal F\to \Sigma_\infty^{(2)}$ is bijective. Consequently
\begin{equation}\label{01-infty}
  \widetilde {\mathcal E}_{l} \cup \widetilde {\mathcal E}_{r}\cup \mathcal F
  =\Pi^{-1}(\Sigma_\infty^{(2)}).
\end{equation}

\end{prop}

\proof\
i) Assume $\omega\in \widetilde{\E_l}$. By \eqref{def-elll}, there are three cases.

 If  $\omega=\omega|_{n-1}1_o3_{er}(0_o3_{er})^\infty$, then $\omega_{n-1}=0_e$ or $2_e$ and $\ell(\omega)=\omega|_{n-1}3_{or}(0_e3_{ol})^\infty$. By tracing on the graph $\mathbb G$, for $m\ge n$,
\begin{equation}\label{gap-1}
  \begin{cases}
  \Pi_\ast(\omega|_m)=\Pi_\ast(\omega|_{n-1})\delta (10)^{[(m+1-n)/2]},\\
  \Pi_\ast(\ell(\omega)|_m)=\Pi_\ast(\omega|_{n-1})\delta1 (01)^{[(m-n)/2]}.
  \end{cases}
\end{equation}
Here $\delta=1$ if $\omega_{n-1}=0_e$ and $\delta=0$ if $\omega_{n-1}=2_e.$

 If  $\omega=\omega|_{n-1}3_{el}(0_o3_{er})^\infty$, then $\omega_{n-1}=0_o$ or $2_o$ and $\ell(\omega)=\omega|_{n-1}1_e3_{ol}(0_e3_{ol})^\infty$. By tracing on the graph $\mathbb G$, for $m\ge n$,
\begin{equation}\label{gap-2}
  \begin{cases}
  \Pi_\ast(\omega|_m)=\Pi_\ast(\omega|_{n-1})\delta 0(10)^{[(m-n)/2]},\\
  \Pi_\ast(\ell(\omega)|_m)=\Pi_\ast(\omega|_{n-1})\delta (01)^{[(m+1-n)/2]}.
  \end{cases}
\end{equation}
Here $\delta=0$ if $\omega_{n-1}=0_o$ and $\delta=1$ if $\omega_{n-1}=2_o.$

If  $\omega=3_{el}(0_o3_{er})^\infty$, then  $\ell(\omega)=(0_e3_{ol})^\infty$. Notice that, in this case $i(\omega)=0$ and $i(\ell(\omega))=\emptyset.$ By tracing on the graph $\mathbb G$, for $m\ge 0$,
\begin{equation}\label{gap-3}
  \begin{cases}
  \Pi_\ast(\omega|_m)= 0(10)^{[m/2]},\\
  \Pi_\ast(\ell(\omega)|_m)= (01)^{[(m+1)/2]}.
  \end{cases}
\end{equation}

In all cases, by \eqref{Pi-L} we have
$$
\Pi(\omega)=\lim_{m\to\infty}\Pi_\ast(\omega|_m)=\lim_{m\to\infty}\Pi_\ast(\ell(\omega)|_m)=\Pi(\ell(\omega)).
$$
Combine with Lemma \ref{bij-ell-o-e} i), \eqref{equal-E-l-r} holds.

ii) Assume $\#\Pi^{-1}(\{\sigma\})\ge2.$ Take any $\omega,\hat\omega\in \Pi^{-1}(\{\sigma\})$. WLOG,
 we assume $\omega\prec\hat\omega$.  Let $w:=\omega\wedge\hat\omega$, and assume $w\in \Omega_{n-1}$. Then $\omega_n\prec\hat\omega_n$. By \eqref{order-label}, $\omega_n\not<\hat\omega_n.$ So either  $n\ge1$ and
$$
(w_{n-1}, \omega_n, \hat\omega_n)\in \{(0_o,3_{el},1_e),(2_o,3_{el},1_e),(0_e,1_o,3_{or}),(2_e,1_o,3_{or})\};
$$
or $n=0$, $w=\emptyset$ and $\omega_0=3_{el}, \hat\omega_0=0_e.$

At first assume $n\ge1$ and  $(w_{n-1}, \omega_n, \hat\omega_n)=(0_o,3_{el},1_e)$. Then
$$
\Pi_\ast(\omega|_n)=\Pi_\ast(w3_{el})=\tau00;\ \ \Pi_\ast(\hat \omega|_n)=\Pi_\ast(w1_e)=\tau0, \text{ where } \tau=\Pi_\ast(w).
$$
Since $\Pi(\omega)=\sigma=\Pi(\hat\omega)$, we have  $\LL(\hat\omega_n\hat\omega_{n+1})=\LL(1_e\hat\omega_{n+1})=0\ast$. By checking $\mathbb G$, the only possibility is  $\hat\omega_{n+1}=3_{ol}$. So $\LL(\hat\omega_n\hat\omega_{n+1})=01$ and $\Pi_\ast(\hat\omega|_{n+1})=\tau001$. Again by $\Pi(\omega)=\Pi(\hat\omega)$, we have $$\LL(\omega_n\omega_{n+1})\LL(\omega_{n+1}\omega_{n+2})=\LL(3_{el}0_o)\LL(0_0\omega_{n+2})
=\LL(0_0\omega_{n+2})=1\ast.$$ By checking $\mathbb G$, the only possibility is  $\omega_{n+1}=3_{er}$. As a consequence $\Pi_\ast(\omega|_{n+2})=\tau0010.$ This process can continue infinitely. By induction, we get
$$
\omega=w3_{el}(0_o3_{er})^\infty\in \widetilde{\E}_l ;\ \ \ \hat\omega=w1_e3_{ol}(0_e3_{ol})^\infty\in \widetilde{\E}_r.
$$
By \eqref{def-elll}, $\hat \omega=\ell(\omega)$.

For the other four cases, the same proof shows that
\begin{equation*}\label{pairing}
\omega\in \widetilde{\E}_l ;\ \ \ \hat\omega\in \widetilde{\E}_r \ \ \ \text{ and }\ \ \ \hat \omega=\ell(\omega).
\end{equation*}

This implies that
$$
\Pi^{-1}(\{\sigma\})\subset \widetilde{\E_l}\cup \widetilde{\E_r},\ \ \#\widetilde{\E_l}\cap \Pi^{-1}(\{\sigma\}), \#\widetilde{\E_r}\cap \Pi^{-1}(\{\sigma\})\le1.
$$
So the first statement  follows.

If $\#\Pi^{-1}(\{\sigma\})\ge2$, then $\Pi^{-1}(\{\sigma\})=\{\omega,\ell(\omega)\}$ for some $\omega\in \widetilde{\E_l}$, hence $\sigma\in \Pi(\widetilde{\E_l})$. If $\sigma\in \Pi(\widetilde{\E_l})$, then by i), $\sigma\in \Pi(\widetilde{\E_r})$, so $\#\Pi^{-1}(\{\sigma\})\ge2$. Thus the second statement holds.

iii) We will prove the Surjectivity of $\Pi$ in Sec. \ref{sec-Pi-surj} by using the coding of $\BB_n$.

Now we show that $\Pi$ preserves the order. Given $\omega\prec \hat\omega$.  Then there exists $n\ge0$ such that  $\omega|_{n-1}=\hat\omega|_{n-1}=:w$ and $\omega|_{n}\prec \hat\omega|_{n}$. If $\omega|_{n}<\hat\omega|_{n}$, then by \eqref{order-label},
$$
\mathcal L(w_{n-1}\omega_{n})=0\ast;\ \ \ \mathcal L(w_{n-1}\hat\omega_{n})=1\ast.
$$
Hence $\Pi(\omega)<\Pi(\hat\omega)$.

Now assume $\omega|_{n}\not<\hat\omega|_{n}$. Then
$$
(\omega_n,\hat\omega_n)\in \{(3_{el},1_e),(1_o,3_{or}),(3_{el},0_e)\}.
$$

At first assume $(\omega_n,\hat\omega_n)=(3_{el},1_e)$. So either $\omega=\omega|_{n-1}3_{el}(0_o3_{er})^\infty:=\tilde \omega$ or there exists $m\ge0$ such that $\omega=\omega|_{n-1}3_{el}(0_o3_{er})^m 0_o\alpha\cdots$ with $\alpha=1_e$ or $3_{el}$. Since
$$
\LL(0_o1_e)=0;\ \ \LL(0_o3_{el})=00;\ \ \LL(0_o3_{er})=10,
$$
We conclude  that $\Pi(\omega)\le \Pi(\tilde \omega).$

The same argument shows that $\Pi(\ell(\tilde \omega))\le\Pi(\hat\omega).$ So by i) we have
$$
\Pi(\omega)\le \Pi(\tilde\omega)=\Pi(\ell(\tilde \omega))\le \Pi(\hat\omega).
$$

For the other two cases, the same argument shows that $\Pi(\omega)\le \Pi(\hat\omega)$. Hence $\Pi$ is order preserving.

By \eqref{two-boundary} and the definition of $\Pi$, one  get \eqref{2-coding-bdy} by tracing  the graph $\mathbb G$.

By ii) $\Pi^{-1}(\Pi(\widetilde{\E_l}))=\widetilde{\E_l}\cup \widetilde{\E_r}$ and $\Pi:\widetilde{\E_l}\to \Pi(\widetilde{\E_l})$ is bijective. Since $\Pi$ is surjective,
the restriction
$\Pi: \Omega_\infty\setminus \left(\widetilde{\E_l}\cup \widetilde{\E_r}\right)\to \Sigma_\infty\setminus \Pi(\widetilde{\E_l})$ is also surjective. Again by ii), this restriction is injective. So it is bijective.
Combine these two bijections, we conclude that $\Pi: \Omega_\infty\setminus  \widetilde{\E_r}\to \Sigma_\infty$ is a bijection.

iv) \eqref{2-coding-Io} and \eqref{2-coding-Ie} follow from \eqref{def-ell-o} and \eqref{def-ell-e} by using the definition of $\Pi$ and tracing  on $\mathbb G.$ Now we show that all the restriction maps in \eqref{bij-I-type} are bijections. By iii), they are all injections. Now we show that they are  surjections. We take the first one as an example. Since $\Pi$ is surjective,  for each $\sigma\in \Sigma_\ast^o$, there exists some $\omega\in \Omega_\infty$ such that $\Pi(\omega)=\sigma01^\infty$. We only need to show that $\omega\in \E_l^o.$ We  delete all the edges of $\mathbb G$ with label containing $0$. Then only  two connected subgraphs are left:
$$
0_e\stackrel{11}{\longrightarrow} 3_{or} \stackrel{\emptyset}{\longrightarrow}0_e;\ \ \ 1_e\stackrel{1}{\longrightarrow} 2_o \stackrel{1}{\longrightarrow}1_e.
$$
Each of these subgraphs can generates $1^\infty$ and $1^\infty$ can only be generated by these two graphs. Thus $\omega\in \E_l^o$ or $\E_l^e$. However, by \eqref{2-coding-bdy} and \eqref{2-coding-Ie}, if $\omega\in \E_l^e$, then $\Pi(\omega)=1^\infty$ or $\sigma01^\infty$ with $\sigma\in \Sigma_\ast^e$.
 So we conclude  that $\omega\in \E_l^o.$ That is, $\Pi_{ol}$ is surjective. The other three cases can be proven by exactly the same way.

v)
By \eqref{gap-edge} and tracing on $\mathbb G,$ it is seen that $\Pi(\mathcal F)\subset \Sigma_\infty^{(2)}$. Combine with \eqref{equal-E-l-r}, we have
$
\Pi( \widetilde {\mathcal E}_{l} \cup \mathcal F)
  \subset\Sigma_\infty^{(2)}.
$
By iii), The restriction of $\Pi$ on $\widetilde {\mathcal E}_{l} \cup \mathcal F$ is injective. Thus we only need to show that $\Sigma_\infty^{(2)}\subset \Pi( \widetilde {\mathcal E}_{l} \cup \mathcal F).$

Fix any $\sigma\in\Sigma_\infty^{(2)}.$ Since $\Pi$ is surjective, there exists $\omega\in \Omega_\infty$ such that $\Pi(\omega)=\sigma.$

At first, we claim that $1_e$ does not appear in $\omega$ infinitely often (i.o.). Indeed if otherwise, by tracing backward  on $\mathbb G$, one of the following words must appear in $\omega$ i.o.:
$$
1_o3_{er}0_o1_e;\ \ 0_o3_{er}0_o1_e;\ \ 2_o3_{el}0_o1_e;\ \ 0_o3_{el}0_o1_e;\ \  2_e2_o1_e;\ \ 1_e2_o1_e.
$$
By applying $\LL$ to these words, we see that for the first four cases, $00$ will appear and for the last two cases, $11$ will appear. So $00$ or $11$ will appears in $\Pi(\omega)$ i.o., a contradiction.

 The same proof shows that $1_o$ also only appear in $\omega$ finite times. Now we go back to the graph $\mathbb G$ and delete $1_o$ and $1_e$ and all the edges connected to $1_o$ and $1_e$. We also delete $0_e3_{or}$ and $0_o3_{el}$ since they have labels $11$ and $00$, respectively. Now the remaining part of the graph has only three connected components:
 $$
0_e\stackrel{01}{\longrightarrow} 3_{ol} \stackrel{\emptyset}{\longrightarrow}0_e;\ \ \ 0_o\stackrel{10}{\longrightarrow} 3_{er} \stackrel{\emptyset}{\longrightarrow}0_o;\ \ \ 2_o\stackrel{0}{\longrightarrow} 2_e \stackrel{1}{\longrightarrow}2_o,
$$
and $(01)^\infty$ can be realized by them and only by them. So $\omega$ can only be eventually $(ab)^\infty$, with
$$ab\in \{0_o3_{er}, 0_e3_{ol}, 2_e2_o\}.$$ That is, $\omega\in \widetilde {\mathcal E}_{l} \cup \widetilde {\mathcal E}_{r}\cup \mathcal F$. Now combine with i), $\Pi(\omega)\in \Pi(\widetilde {\mathcal E}_{l}\cup \mathcal F)$.
Hence, $\Sigma_\infty^{(2)}\subset \Pi(\widetilde {\mathcal E}_{l}\cup \mathcal F).$ Then the restriction of $\Pi$ on $\widetilde {\mathcal E}_{l}\cup \mathcal F$ is bijective.

Now by i), \eqref{01-infty} follows.
\hfill $\Box$

\begin{rem}\label{Pi-ast-order}
{\rm
Recall that we have defined an order $\preceq$ on $\Sigma_\ast$ by \eqref{order-Sigma-ast}.  If $(\omega,\hat\omega)\in \G_{II}$, by  \eqref{gap-1}-\eqref{gap-3}, we conclude that for any $m,m'> |\omega\wedge\hat\omega|$,
$$
\Pi_\ast(\omega|_m)\prec \Pi_\ast(\hat\omega|_{m'}).
$$

}
\end{rem}

\subsubsection{Proof of Theorem \ref{symbolic-gap}}\

Before giving the formal proof, let us explain the mechanism to form a gap in $\Omega_\infty$. At first, if $(\sigma,\hat\sigma)$ is gap of $\Sigma_\infty$, then it can be lifted to a gap of $\Omega_\infty$ via $\Pi.$ It is seen that  the gap of $\Sigma_\infty$ has the form $(\tau01^\infty,\tau10^\infty)$. Secondly, if $\Pi^{-1}(\sigma)$ has more than one elements, then by Proposition \ref{Pi-basic} ii), it has exactly two elements, and the elements form a gap. These are the only two ways to form a gap.

\noindent {\bf Proof of Theorem \ref{symbolic-gap}.}\
i) Fix any $\omega\ne\omega_\ast$. By \eqref{2-coding-bdy}, $\Pi(\omega_\ast)=0^\infty\notin \Sigma_\infty^{(2)}$. By Proposition  \ref{Pi-basic} ii), $\Pi(\omega)\ne 0^\infty$.  We claim that $\omega_\ast\prec \omega$. Indeed, if otherwise,
$\omega\preceq \omega_\ast,$ then $\Pi(\omega)\le \Pi(\omega_\ast)=0^\infty$.  Since $0^\infty$ is the minimum of $\Sigma_\infty$, we have $\Pi(\omega)=0^\infty$, a contradiction. So $\omega_\ast\prec \omega.$

The same argument shows that if $\omega\ne\omega^\ast$, then $\omega\prec \omega^\ast$. So the result follows.

ii) Combine with Lemma \ref{bij-ell-o-e}, we only need  to show \eqref{gap-Io} and \eqref{gap-Ie}. Fix $\omega\in \E_l^o$. By \eqref{2-coding-Io}, there exists $\sigma\in \Sigma_\ast$ such that $\Pi(\omega)=\sigma01^\infty=:\tau$ and $\Pi(\ell^o(\omega))=\sigma10^\infty=:\hat\tau.$ It is well-known that $(\tau,\hat\tau)$ is a gap of $\Sigma_\infty.$ We claim that $(\omega,\ell^o(\omega))$ is a gap of $\Omega_\infty$. Indeed if otherwise, there exists $\tilde \omega$ such that $\omega\prec \tilde \omega\prec \ell^o(\omega)$, so
$$
\tau\le \Pi(\tilde\omega)\le\hat\tau.
$$
So either $\Pi(\tilde \omega)=\tau$ or $\hat\tau$. In the former case, $\Pi^{-1}(\{\tau\})\supset\{\omega,\tilde\omega\}$. So by Proposition \ref{Pi-basic} ii), $\tau\in \Sigma_\infty^{(2)}$, which is a contradiction. In the latter case, we get a contradiction by the same reasoning. So $(\omega,\ell^o(\omega))$ is a gap of $\Omega_\infty$ and \eqref{gap-Io} holds.

The same proof shows that \eqref{gap-Ie} holds.

iii) Combine with Lemma \ref{bij-ell-o-e}, we only need to show \eqref{gap-II}. Take $\omega\in \widetilde{\E_l}$, we claim that $(\omega,\ell(\omega)) $ is a gap. Indeed if otherwise, there exists $\tilde \omega$ such that $\omega\prec \tilde \omega\prec \ell(\omega)$, so
$$
\Pi(\omega)\le \Pi(\tilde\omega)\le\Pi(\ell(\omega)).
$$
By Proposition \ref{Pi-basic} i), $\Pi(\omega)=\Pi(\ell(\omega))=:\sigma$. So $\Pi^{-1}(\{\sigma\})\supset\{\omega,\tilde\omega,\ell(\omega)\}$, which contradicts with Proposition \ref{Pi-basic} ii). So \eqref{gap-II} holds.

iv) Now assume $(\omega,\hat\omega)\in \G.$ Then either $\Pi(\omega)=\Pi(\hat\omega)=:\sigma$, or $\Pi(\omega)<\Pi(\hat\omega)$.

In the former case, by Proposition \ref{Pi-basic} ii), $\omega\in \widetilde{\E_l}$ and $\hat\omega=\ell(\omega)$. So $(\omega,\hat\omega)\in \G_{II}$.

In the latter case, we claim that  $(\Pi(\omega),\Pi(\hat\omega))$ is a gap of $\Sigma_\infty.$ Indeed, if otherwise, there exists $\sigma\in\Sigma_\infty$ such that $\Pi(\omega)<\sigma<\Pi(\hat\omega).$ Assume $\Pi(\tilde \omega)=\sigma$, then $\tilde\omega\ne \omega,\hat\omega$. Since $\Pi$ preserve the order, we must have $\omega\prec \tilde \omega\prec \hat\omega$, a contradiction.

 So there exists $\sigma\in \Sigma_\ast$ such that $\Pi(\omega)=\sigma01^\infty$ and $\Pi(\hat\omega)=\sigma10^\infty.$ By Proposition \ref{Pi-basic} iv) and ii), either $\omega\in \E_l^o, \hat\omega=\ell^o(\omega)$, or $\hat\omega\in \E_r^e, \omega=\ell^e(\hat\omega)$. So $(\omega,\hat\omega)\in \G_I.$

 As a result, $\G\subset \G_I\cup\G_{II}$. Combine with ii) and iii), we get the equality.
\hfill $\Box$

We end with the following lemma, which will be used to show that $\pi$ is order-preserving.

\begin{lem}\label{prec-not-<}
Assume $\omega\prec \hat\omega$  but $\omega\not<\hat\omega.$ Then there exists $(\tau,\hat\tau)\in \G_{II}$ such that $\omega\preceq \tau$ and $\hat\tau \preceq \hat\omega$. Moreover if $\omega\ne\tau (\hat \tau\ne \hat \omega)$, then $\omega<\tau (\hat\tau<\hat\omega)$.
\end{lem}

\proof\
Assume $n\ge0$ is  such that  $\omega|_{n-1}=\hat\omega|_{n-1}=:w$ and $\omega|_{n}\prec \hat\omega|_{n}$.
Since $\omega\not<\hat\omega,$
$$
(\omega_n,\hat\omega_n)\in \{(3_{el},1_e),(1_o,3_{or}),(3_{el},0_e)\}.
$$

At first assume $(\omega_n,\hat\omega_n)=(3_{el},1_e)$. So either $\omega=\omega|_{n-1}3_{el}(0_o3_{er})^\infty=:\tau\in\widetilde {\mathcal E}_l$ or there exists $m\ge0$ such that $\omega=\omega|_{n-1}3_{el}(0_o3_{er})^m 0_o\alpha\cdots$ with $\alpha=1_e$ or $3_{el}$. In the latter case, by \eqref{order-strong}, $\omega< \tau.$ The same argument shows that $\ell(\tau)\preceq\hat\omega$ and if $\hat\omega\ne \ell(\tau)$, then $\ell(\tau)<\hat\omega$.

For the case that $(\omega_n,\hat\omega_n)=(1_o,3_{or})$ or $(3_{el},0_e)$, the proof is the same.
\hfill $\Box$


\section{Coding and the Hausdorff  dimension of the spectrum  }\label{sec-coding-Hausdorff}

In this section, we prove a weak version of Theorem \ref{main-coding-spectrum}, which is enough for estimating the dimension of the spectrum. Then we prove Theorem \ref{main-low-dim}.

\subsection{Coding of the spectrum}

\subsubsection{Coding of $\BB_n$}\label{Sec-code-BB-n}\

We start with coding  $\BB_n$ by $\Omega_n$. We summarize Lemmas \ref{type-3}-\ref{type-2} in the following Corollary, the proof of which is by direct checking.

\begin{cor}\label{admissible-explain}
Assume $B\in \BB_n$ has type $\alpha$. Then $B$ contains a
 band $\tilde B\in \BB_{n+1}$ of type $\beta$
if and only if $\alpha\to \beta.$ Moreover, if $\alpha\to\beta$, then $B$ contains exact one
band $\tilde B\in \BB_{n+1}$ of type $\beta$.
\end{cor}

Now we can code $\BB_n$ by $\Omega_n$ inductively as follows.

 Recall that $\BB_0=\{B_0, B_\emptyset\}$; $B_0$ has type
$3_{el}$ and $B_\emptyset$ has type $0_e$ and by \eqref{initial},
 $B_0\prec B_\emptyset.$ We code $\BB_0$ by
$\Omega_0=\{3_{el}, 0_e\}$ as
$$
I_{3_{el}}:=B_0;\ \ I_{0_e}:=B_\emptyset.
$$
 Assume we have coded $\BB_{n}$ by $\Omega_n$:
 $$\BB_n=\{I_w: w\in \Omega_n\},\ \ I_w \text{ has type } w_n \ \ \text{ and }\ \  I_v\ne I_w \ \text{if}\ v\ne w.$$
 Take any  $u=u_0\cdots u_{n+1}\in \Omega_{n+1}$, then $u|_n\in \Omega_n$ and $u_n\to u_{n+1}.$
  by Corollary \ref{admissible-explain}, there exists a unique
 $\hat B\in \BB_{n+1}$ with type $u_{n+1}$ such that $\hat B\subset I_{u|_n}$.
 Define  $I_u:=\hat B$. Then
 $$
 \{I_w: w\in \Omega_{n+1}\}\subset \BB_{n+1}.
 $$
 Now  take any $\tilde B\in \BB_{n+1}$, assume it has type $\beta$.
 By Proposition \ref{unique-father}, there exists a unique
    $B\in \BB_n$ such that $\tilde B\subset B$. By induction hypothesis, there is a unique
     $w\in \Omega_n$
     such that  $ B=I_{w}$. Again by
     Corollary \ref{admissible-explain}, $w_n\to\beta$ and $\tilde B=I_{w\beta}$ and no
      other $u\in \Sigma_{n+1}$ satisfies $\tilde B=I_u$.
     Thus
 $$
 \BB_{n+1}=\{I_w: w\in \Omega_{n+1}\},\ \ I_w \text{ has type } w_{n+1} \ \ \text{ and }\ \  I_v\ne I_w \ \text{if}\ v\ne w.
 $$
 By induction, we finish the coding process.

 Combine with Lemmas \ref{type-3}-\ref{type-2}, the above coding process also implies that for $w\in \Omega_\ast$,
 \begin{equation}\label{I-and-B}
   I_w=B_{\Pi_\ast(w)}.
 \end{equation}

 By Corollary \ref{admissible-explain}, we also have the following: Assume $v\in \Omega_n$ and $w\in \Omega_{n-1}$, then
 \begin{equation}\label{inclusion-prefix}
 I_{v}\subset I_w \Leftrightarrow w\lhd v.
 \end{equation}

 With this coding, we can explain two orders \eqref{order} and \eqref{order-strong} as follows:

 \begin{cor}\label{order-explain}
 Assume $w\in \Omega_n$ and $w_n\to \alpha, \beta$. Then
 $$
 I_{w\alpha}\prec I_{w\beta}\Leftrightarrow
 \alpha\prec \beta; \ \ \ I_{w\alpha}< I_{w\beta}\Leftrightarrow\alpha< \beta.
 $$
 \end{cor}

 The proof is again  by  applying   Lemmas \ref{type-0}-\ref{type-2} and checking directly.

\subsubsection{Coding of $\sigma_\lambda$}\label{sec-coding-sigma-lambda}\


For any $\omega\in \Omega_\infty$, we have $\omega|_n\in \Omega_n$. Hence $I_{\omega|_n}\in \BB_n$.
 By the construction, we have $I_{\omega|_{n+1}}\subset I_{\omega|_n}.$ By \eqref{band-length},
 $|I_{\omega|_n}|\to0$. Hence $\bigcap_{n\ge 0} I_{\omega|_n}$ is a singleton. By Corollary \ref{NS-BB}, this point is in $\sigma_\lambda.$

 We define the coding map $\pi: \Omega_\infty \to \sigma_\lambda$ as
 \begin{equation*}\label{coding-map}
 \pi(\omega):=\bigcap_{n\ge 0} I_{\omega|_n}.
 \end{equation*}

We have the following weak version of Theorem \ref{main-coding-spectrum}:

\begin{prop}\label{main-coding-spectrum-weak}

  $\pi: (\Omega_\infty,\preceq)\to (\sigma_\lambda,\le )$ is continuous, surjective and   preserves the orders.
\end{prop}

To prove it, we need two lemmas.

 \begin{lem}\label{order-<-pre}
 Given $\omega,\hat\omega\in \Omega_\infty$. If $\omega<\hat\omega$, then $\pi(\omega)<\pi(\hat\omega).$
 \end{lem}

 \begin{proof}
 Assume $n$ is such that
  $\omega|_{n-1}=\hat\omega|_{n-1}$ and $\omega_{n}<\hat\omega_{n}$. By Corollary \ref{order-explain},
   $I_{\omega|_{n}}<I_{\hat\omega|_n}.$
    Since $\pi(\omega)\in I_{\omega|_n}$ and $\pi(\hat \omega)\in I_{\hat\omega|_n}$,
    we conclude that $\pi(\omega)<\pi(\hat\omega)$.
 \hfill $\Box$
 \end{proof}

 \begin{lem}\label{order-gap-II-weak}
  If $(\omega,\hat\omega)\in \G_{II}$, then
 $
 \pi(\omega)\le \pi(\hat \omega).
$
 \end{lem}

 \begin{proof}  Write
 $$
 \sigma^{(n)}=\Pi_\ast(\omega|_n);\ \ \ \hat\sigma^{(n)}=\Pi_\ast(\hat\omega|_n).
 $$
 By Remark \ref{Pi-ast-order}, if  $n>|\omega\wedge\hat\omega|$, then $\sigma^{(n)}\prec\hat\sigma^{(n)}. $
 By Proposition \ref{zero-order}, $z_{\sigma^{(n)}}<z_{\hat\sigma^{(n)}}$.

 On the other hand, by \eqref{I-and-B},
 $$
 \pi(\omega)\in I_{\omega|_n}=B_{\sigma^{(n)}}; \ \ \ \pi(\hat\omega)\in I_{\hat\omega|_n}=B_{\hat\sigma^{(n)}}.
 $$
By \eqref{band-length}, the length of the bands  $B_{\sigma^{(n)}}$ and $B_{\hat\sigma^{(n)}}$
  tends to zero as $n$ tends to infinity. Since
  $
 z_{\sigma^{(n)}}\in B_{\sigma^{(n)}}
 $
 and $z_{\hat\sigma^{(n)}}\in B_{\hat\sigma^{(n)}}$, we have
 \begin{equation*}\label{lim-E12}
 \pi(\omega)=\lim_n z_{\sigma^{(n)}}\le \lim_n z_{\hat\sigma^{(n)}}
 =\pi(\hat\omega)
\end{equation*}
So the result follows.
\hfill $\Box$
 \end{proof}

\begin{rem}
{\rm In Lemma \ref{order-gap-II},  we will show that the strict inequality holds.  Here we just mention that the proof is highly nontrivial. We need to study the $\infty_{II}$ energies to finally reach the strict inequality.
}
\end{rem}

\noindent {\bf Proof of Proposition \ref{main-coding-spectrum-weak}.}\
 By Corollary \ref{NS-BB}, $\pi $ is surjective.

 Assume $|\omega\wedge\hat\omega|=n$, then $d(\omega,\hat\omega)=2^{-n}$. We have
 $$
 |\pi(\omega)-\pi(\hat\omega)|\le |I_{\omega\wedge\hat\omega}|.
 $$
 Combine with \eqref{band-length},  $\pi$ is continuous.

 Assume  $\omega\prec \hat\omega.$ If $\omega<\hat\omega,$  by Lemma \ref{order-<-pre}, $\pi(\omega)<\pi(\hat\omega)$. If $\omega\not<\hat\omega,$ by Lemma \ref{prec-not-<}, there exists $(\tau,\hat\tau)\in \G_{II}$ such that $\omega\preceq \tau$ and $\hat\tau \preceq \hat\omega$, moreover if $\omega\ne\tau (\hat \tau\ne \hat \omega)$, then $\omega<\tau (\hat\tau<\hat\omega)$. By Lemmas \ref{order-<-pre} and \ref{order-gap-II-weak},
 $$
 \pi(\omega)\le \pi(\tau)\le\pi(\hat \tau)\le \pi(\hat\omega).
 $$
So $\pi$ is order-preserving.
\hfill $\Box$

\subsubsection{Surjectivity of $\Pi$}\label{sec-Pi-surj}\

As an application of the  coding of $\BB_n$, we can show that $\Pi$ is surjective.

 Given $\sigma\in \Sigma_\infty.$ For each $n\ge0$, we always have $B_{\sigma|_n}\in \BB_n$.
By
 the coding of $\BB_n$, there exists a
unique $w^{(n)}\in\Omega_n$ such that $I_{w^{(n)}}=B_{\sigma|_n}$.

\noindent {\bf Claim 1:  } For each $n\ge2$,
   $
   w^{(n-1)}\lhd w^{(n)}$ or $ w^{(n-2)}\lhd w^{(n)}.
   $

\noindent $\lhd $ If $B_{\sigma|_n}\notin\BB_{n-1},$ then $B_{\sigma|_n}\subset \tilde B$ for some $\tilde B\in\B_{n-1}$.
by Corollary \ref{tech-cor} iii), $\tilde B=B_{\sigma|_{n-1}}.$  That is,
 $I_{w^{(n)}}\subset I_{w^{(n-1)}}$. By \eqref{inclusion-prefix}, $w^{(n-1)}\lhd w^{(n)}$.

 If $B_{\sigma|_n}\in\BB_{n-1},$ then $B_{\sigma|_n}\not\subset B_{\sigma|_{n-1}}$.  By  Corollary \ref{tech-cor} ii),
 $B_{\sigma|_n}\subset B_{\sigma|_{n-2}}$. Assume $\hat w^{(n-1)}\in \Omega_{n-1}$ is such that $B_{\sigma|_n}=I_{\hat w^{(n-1)}}$.
   By \eqref{inclusion-prefix}, we have $w^{(n-2)}\lhd \hat w^{(n-1)}$. We also have
   $I_{ w^{(n)}}=I_{\hat w^{(n-1)}}$, again by \eqref{inclusion-prefix}, we have $\hat w^{(n-1)}\lhd w^{(n)}$. So, $w^{(n-2)}\lhd w^{(n)}$.
   \hfill $\rhd$

 By the claim, we conclude that
 \begin{equation}\label{prefix}
\text{ for any }\  n>k\ge0, \text{ either}\  w^{(k)}\lhd w^{(n)}\ \text{ or }\  w^{(k+1)}\lhd w^{(n)}.
 \end{equation}

Define  inductively an integer  sequence  $(m_k)_{k\ge1}$ and a decreasing  sequence $(\Xi_k)_{k\ge1}$ with $\Xi_{k}\subset \N$ and $\#\Xi_k=\infty$  as follows. Write $\Xi_0:=\N.$
By \eqref{prefix}, at least one of $\{n\in \Xi_0: w^{(0)}\lhd w^{(n)}\}$ and $\{n\in \Xi_0: w^{(1)}\lhd w^{(n)}\}$ is infinite.
If $\{n\in \Xi_0: w^{(0)}\lhd w^{(n)}\}$ is infinite, define $m_1:=0$. If otherwise, then $\{n\in \Xi_0: w^{(1)}\lhd w^{(n)}\}$ is infinite, define $m_1:=1$.  Define $\Xi_1:=\{n\in \Xi_0: w^{(m_1)}\lhd w^{(n)}\}.$ Then $\#\Xi_1=\infty.$

Assume $m_k$ and $\Xi_k$ has been defined with the desired property. By \eqref{prefix}, at least one of $\{n\in \Xi_k: w^{(m_k+1)}\lhd w^{(n)}\}$ and $\{n\in \Xi_k: w^{(m_k+2)}\lhd w^{(n)}\}$ is infinite. If $\{n\in \Xi_k: w^{(m_k+1)}\lhd w^{(n)}\}$ is infinite, define $m_{k+1}:=m_{k}+1$. If otherwise, then $\{n\in \Xi_k: w^{(m_k+2)}\lhd w^{(n)}\}$ is infinite, define $m_{k+1}:=m_k+2$.  Define $\Xi_{k+1}:=\{n\in \Xi_k: w^{(m_{k+1})}\lhd w^{(n)}\}.$ Hence $\Xi_{k+1}\subset \Xi_k$ and is infinite.

By induction, we finish the definition.

Now we claim that for any $k$, $w^{(m_k)}\lhd w^{(m_{k+1})}$. Indeed, by the defining process, we have $w^{(m_k)}\lhd w^{(n)}$ for any $n\in \Xi_k$. Choose some $n\in \Xi_{k+1}$, then $ n\in \Xi_k$ since $\Xi_{k+1}\subset \Xi_k$. Hence  we have  $w^{(m_k)},w^{(m_{k+1})}\lhd w^{(n)}$. Since $|w^{(m_k)}|< |w^{(m_{k+1})}|$, we conclude that $w^{(m_k)}\lhd w^{(m_{k+1})}$.

Let $\omega=\lim_{k\to\infty} w^{(m_k)}$. Then we have $\Pi(\omega)=\sigma.$
\hfill $\Box$

It is an  interesting question  to find a purely combinatorial proof for the surjectivity of $\Pi.$

\subsection{Lower bound for the Hausdorff dimension of $\sigma_\lambda$}\

In this subsection, at first we present a sufficient condition for estimating from below the Hausdorff dimension of the limit set of a SNS. Then we apply it to a sub-SNS of $\{\BB_n:n\ge0\}$ to obtain a lower bound for the Hausdorff dimension of $\sigma_\lambda$.

\subsubsection{Dimension estimation for limit set of SNS}\

We need the following general result on the dimension of limit set of SNS, which
 is also useful in its own right.

\begin{prop}\label{SNS-lower-dim}
Let $\I=\{\I_n:n\ge0\}$ be a SNS. Assume $\I$ satisfies:

i) There exist $\lambda\in (0,1)$ and $C>0$ such that
$$
|I|\ge C\lambda^n,\ \ \forall n\ge0, \forall I\in \I_n.
$$

ii) There exists $C'>0$ such that for any $n,k\ge0$, and  $I,I'\in \I_n$,
$$
\frac{\#\{J\in \I_{n+k}: J\subset I\}}{\#\{J\in \I_{n+k}: J\subset I'\}}\le C'.
$$

Then the limit set satisfies
\begin{equation*}
  \dim_H A(\I)\ge \liminf_{n\to\infty}\frac{\log \# \I_n}{-n\log \lambda}.
\end{equation*}
\end{prop}

\proof
For any $n\ge0$, let
$$\kappa_n=\# \I_n,\ \alpha=\liminf_{n\to\infty}\sqrt[n]{ \kappa_n},\ t=-\log\alpha/\log\lambda.$$
If $\alpha=1$, we are done. Suppose $\alpha>1$.
Take any $0<\delta<t$.
There exists $0<\varepsilon<\alpha$ such that, for any large enough $n$,
$$t-\delta<\frac{-n\log(\alpha-\varepsilon)}{(n+1)\log\lambda+\log C}.$$
Note that, for any large enough $n$,
$\kappa_n>(\alpha-\varepsilon)^n.$

For any $n\ge0$, define a probability measure $\mu_n$ on $\R$ as follows: $\mu_n$ is supported on
$\bigcup_{B\in \I_n}B$.  For any $B\in \I_n$,
$$\mu_n(B)=\kappa_{n}^{-1},$$
and $\mu_n$ is uniform on the interval $B$.

Take any $k,n>0$. Let $B\in \I_n$. Suppose there are $m$ intervals in $\I_{n+k}$ that is contained in $B$.
Then $\mu_{n+k}(B)=m/\kappa_{n+k}$. By ii),
$$m/C'\le \frac{\kappa_{n+k}}{\kappa_n}\le C'm.$$
So we  have
\begin{equation}\label{mu-n+k}
\kappa_n^{-1}/C'\le \mu_{n+k}(B)\le C' \kappa_n^{-1}.
\end{equation}
Take a weak-star limit of $\{\mu_n\}_{n\ge1}$, say $\mu$.
Since $\I$ is a SNS, $\mu$ is a probability measure supported on $A(\I)$.
Moreover, by \eqref{mu-n+k}, for any $n>0$ and $B\in \I_n$,
$$\kappa_n^{-1}/C'\le \mu(B)\le C' \kappa_n^{-1}.$$

For any open interval $U$ with $|U|$ small enough, there exists $k\in\N$  such that
$$C \lambda^{k+1}<|U|\le C \lambda^{k}.$$
At most two intervals in $\I_k$   intersect  $U$. Then
$$\mu(U)\le 2C' \kappa_k^{-1}\le 2C'(\alpha-\varepsilon)^{-k}\le
2C'|U|^{\frac{-k\log (\alpha-\varepsilon)}{(k+1)\log\lambda+\log C}}
\le 2C'|U|^{t-\delta}.$$
By mass distribution principle,
$$\dim_H A(\I)\ge t-\delta.$$
Since $\delta>0$ can be arbitrarily small, the result follows.
\hfill $\Box$

\subsubsection{Lower bound for $\dim_H\sigma_\lambda$}\label{sec-lower-dim}\

Now we construct a sub-NS of $\{\BB_n:n\ge0\}$ such that it is a SNS and satisfies the conditions in Proposition \ref{SNS-lower-dim}.

 Define a sub-alphabet $\widetilde A=\{1_e,1_o,2_e,2_o\}$  of $\A$ and a sub directed graph $\widetilde {\mathbb G}$ with  restricted admissible relation:
$$
1_e\to 2_o;\ \ 1_o\to 2_e;\ \ 2_e\to 1_o, 2_o;\ \  2_o\to 1_e, 2_e.
$$
Consider the sub-NS $\widetilde \BB=\{\widetilde \BB_n: n\ge1\}$ defined by
$$
\widetilde \BB_n:=\{I_{0_e1_ow}\in \BB_{n+1}: w=w_1\cdots w_n, w_i\in \widetilde A\}.
$$

\begin{prop}\label{SNS}
$\widetilde \BB$ is a SNS and its limit set $A(\widetilde \BB)\subset \sigma_\lambda\cap B_\infty.$ Moreover There exists constant $C>0$  such that for  all $n,k\ge0$ and  $I, I'\in \widetilde \BB_n,$
$$
|I|\ge C\ 4^{-n} \ \ \text{ and }\ \ \frac{\#\{J\in \widetilde \BB_{n+k}: J\subset I\}}{\#\{J\in \widetilde \BB_{n+k}: J\subset I'\}}\le 2.
$$
\end{prop}

\proof
One check directly that $\widetilde \BB$ is a sub-NS of $\BB$, so $A(\widetilde \BB)\subset A(\BB)=\sigma_\lambda$. Moreover, since $2_e<1_e$ and $1_o<2_o$, combine with Corollary \ref{order-explain}, one conclude that $\widetilde \BB$ is a SNS.

Take any $n>0$ and any $I\in \widetilde \BB_n$. We claim that for any $1\le  k\le n+1$ and any $x\in I$,
$$|h_k(x)|\le2.$$
Indeed, assume $I=I_{0_e1_ow}$ with $w\in \widetilde \A^n$, then by \eqref{I-and-B} and tracing on $\mathbb G,$  for any $k=0,\cdots,n$ we have
$$
I^{(k)}:=I_{0_e1_ow_1\cdots w_k}=B_{1\sigma_1\cdots\sigma_k}, \text{ where }\ \ \sigma_i\in \{0,1\}.
$$
Thus $I^{(k)}\in \B_{k+1}$ and $I=I^{(n)}\subset I^{(n-1)}\subset \cdots\subset  I^{(0)}=B_1$. Since $|h_k|\le 2$ on any $B\in \B_k$, the claim follows.

The claim implies immediately that $A(\widetilde \BB)\subset B_\infty.$

Recall that
$$h_{k+1}(x)=h_{k}(x)(h_{k-1}^2(x)-2)-2,$$
we have
$$h'_{k+1}(x)= (h^2_{k-1}(x)-2)h'_k(x)+2h_{k-1}(x)h_k(x)h'_{k-1}(x).$$
Then, for $x\in I$ and $2\le k\le n$,
$$|h'_{k+1}(x)|\le 2|h'_k(x)|+8|h'_{k-1}(x)|.$$
This implies there exists $c>0$ depending only on $B_1$ such that, for $x\in I$ and $1\le k\le n+1$,
$$|h'_{k}(x)|\le c 4^k.$$
By Floquet theory, $h_{n+1}$ is monotone on $I\in \B_{n+1}$ and
$h_{n+1}(I)=[-2,2]$. Writing $I=[x_0,x_1]$, we have
$$4=\int_{x_0}^{x_1}|h'_{n+1}(x)|dx\le c(x_1-x_0)4^{n+1}.$$
Thus
$|I|=x_1-x_0\ge c^{-1} 4^{-n}.$

It is seen that $\widetilde \BB$ is corresponding to a sub-shift of finite type with alphabet $\{1_o,2_o,1_e,2_e\}$ and incidence  matrix
$$\tilde A=
\left[\begin{array}{cccc}
0&0&0&1\\
0&0&1&1\\
0&1&0&0\\
1&1&0&0
\end{array}\right].
$$
By this, one can show easily that
$$\frac{\#\{J\in \widetilde \BB_{n+k}: J\subset I\}}{\#\{J\in \widetilde \BB_{n+k}: J\subset I'\}}\le \frac{F_k}{F_{k-1}},$$
 where $\{F_k:k\ge0\}$ is the Fibonacci sequence defined by $F_0=1, F_1=2$ and $F_{n+1}=F_n+F_{n-1}$. Then the result follows easily.
\hfill $\Box$

\noindent {\bf Proof of Theorem \ref{main-low-dim}.} By the definition of $\widetilde \BB,$ it is seen that $\# \widetilde \BB_n=F_n$.  It is well known that $F_n\sim \alpha^n$. By Propositions \ref{SNS} and \ref{SNS-lower-dim},
\begin{equation}\label{dim-lim-set}
\dim_H \sigma_\lambda\ge \dim_H A(\widetilde \BB)\ge \liminf_{n\to\infty}\frac{\log\# \widetilde \BB_n}{n\log 4}=\frac{\log \alpha}{\log 4}.
\end{equation}
\hfill $\Box$

\begin{rem}\label{rem-dim}
{\rm

1)
Indeed, we can start from any band $B$ with type in $\widetilde \A$ and do the same construction by following the sub graph $\widetilde{\mathbb G}$. Then we construct a Cantor subset of $\sigma_\lambda\cap B$ with the same dimension estimation. This means that the ``local" dimension of $\sigma_\lambda$ is also uniformly bounded from below.

2) If we view the SNS $\widetilde \BB$ as a dynamical system coded by the subshift $\Omega_{\tilde A}$, then  $\log \alpha$ is just the topological entropy of the system and $\log 4$ can be viewed as an upper bound of the lyapunov exponents.
}\end{rem}


\section{ $\infty$-type energies, gaps of the spectrum and zeros of trace polynomials }\label{sec-infty-energy}

In this section, at first, we establish  the existence of  $\infty$-type energies;  then we  prove Theorems \ref{main-infty-energy}, \ref{main-coding-spectrum} and  \ref{main-gap-zero-inftyenergy}.

\subsection{$\infty_I$-type energies}\

We will show that  $E\in \pi(\mathcal E_{l}^e\cup \mathcal E_{r}^o)$ is an $\infty$-type energy. Recall that in Sec. \ref{subsec-gap-zero-infty} we call it  $\infty_I$-type  energy. Indeed, we   can be more precise:

\begin{prop}\label{infty-energy-I}

i)
$\pi(\mathcal E^e_{l})$ is  dense in $\sigma_\lambda.$ If $E\in \pi(\mathcal E_{l}^e),$ then
$$
\lim_{n\to\infty} |h_{2n+1}(E)|\to\infty\ \ \text{ and }\ \ \ \lim_{n\to\infty} |h_{2n}(E)|\to\sqrt{2}.
$$

ii) $\pi(\mathcal E^o_{r})$  is  dense in $\sigma_\lambda.$ If $E\in \pi(\mathcal E_{r}^o),$ then
$$
\lim_{n\to\infty} |h_{2n}(E)|\to\infty\ \ \text{ and }\ \ \ \lim_{n\to\infty} |h_{2n+1}(E)|\to\sqrt{2}.
$$

\end{prop}

We  only prove i), since the proof of ii) is the same.

\begin{lem}\label{evenendpoints}

Suppose $\sigma\in\Sigma_n$ and $B_\sigma$ has type $0_e$. Then
\begin{equation}\label{cascade}
b_{\sigma1}<b_{\sigma111}<b_{\sigma11}<b_\sigma,  
\end{equation}
and
\begin{eqnarray}
 \label{bd-h-n+1} h_{n+1}(E) \ge 2  \ \ \ \text{ for } \ \ E\in [b_{\sigma1},b_\sigma], \\
 \label{bd-h-n} h_n(E) \ge\sqrt{2} \ \ \ \ \text{ for } \ \ E\in [a_{\sigma11},b_\sigma].
\end{eqnarray}
\end{lem}

\proof\
Since $B_\sigma$ has type $0_e$,  $n$ is even and $b_\sigma\notin\RR_{n-1}$.
 By Lemma \ref{endsymbol} i), $a_{\sigma1}=z_\sigma$. By Lemma \ref{endsymbol} iii),
  $b_{\sigma1^k}\notin\RR_{n+k-1}$ for any $k\ge1.$
  Applying Lemma \ref{tech-2} i) to
 $B_{\sigma1}$, we have
 $  B_{\sigma1}\prec B_{\sigma111}.$ So $b_{\sigma1}<b_{\sigma111}$.
 Applying Lemma \ref{tech-1} ii) to
 $B_{\sigma11}$, we have
 $  B_{\sigma111}\subset {\rm int} (B_{\sigma11}).$ So $b_{\sigma111}<b_{\sigma11}$.
By Lemma
\ref{tech-2} ii), $a_{\sigma11}\notin\RR_{n+1}$ and
$
 B_{\sigma1}\prec B_{\sigma11}\subset {\rm int}(B_\sigma).
$
So $z_\sigma=a_{\sigma1}<a_{\sigma11}$ and   $b_{\sigma11}<b_\sigma.$ Thus \eqref{cascade} follows.

If $\sigma=1^n$, then $b_{\sigma1}=b_{1^{n+1}}$ is the maximal root of $h_{n+1}(E)=2$. Hence $h_{n+1}(E)\ge2$ for any $E\in [b_{\sigma1},\infty)$. Now assume $\sigma<1^n$.
Since $b_\sigma\notin \RR_{n-1}$, by Lemma \ref{endsymbol} iv), $a_{\sigma^+}\in \RR_{n-1}$. By Lemma \ref{endsymbol} ii),
  $a_{\sigma^+0}=a_{\sigma^+}$.
 By Proposition \ref{basic-facts} i), $h_{n+1}$ is increasing on $B_{\sigma1}$ and decreasing on $B_{\sigma^+0}$, hence,
 $$h_{n+1}(b_{\sigma1})=2;\ \ h_{n+1}(a_{\sigma^+0})=2;\ \ h_n([b_{\sigma1},a_{\sigma^+0}])\subset[2,\infty).$$
  Since
  $b_{\sigma1}<b_\sigma<a_{\sigma^+}=a_{\sigma^+0}$, \eqref{bd-h-n+1} holds.

  By Proposition \ref{basic-facts} i), we have
 $h_{n+2}(a_{\sigma11})=-2.$ So
 $$-2=h_{n+2}(a_{\sigma11})=h_{n+1}(a_{\sigma11})(h_{n}(a_{\sigma11})^2-2)-2.$$
 Since $a_{\sigma11}\notin\RR_{n+1}$, $h_{n+1}(a_{\sigma11})\ne0.$ So we must have
 $h_{n}(a_{\sigma11})=\pm \sqrt{2}$. By Lemma \ref{0-momotone},
 $h_n$ is increasing on $B_\sigma$. Since $z_\sigma<a_{\sigma11}<b_{\sigma11}<b_\sigma,$ we must have
 $h_{n}(a_{\sigma11})=\sqrt{2}$. Consequently  \eqref{bd-h-n} holds.
 \hfill $\Box$

\noindent{\bf Proof of Proposition \ref{infty-energy-I}.}\
We only show i). At first we show that $\pi(\mathcal E_{l}^e)$ is dense in $\sigma_\lambda$.
Since $\{I_w: w\in \Omega_n\}$ is a covering
of the spectrum and $\max\{|I_w|: w\in \Omega_n\}\to0$ as $n\to\infty,$ we only need to show that  $I_w\cap \pi(\mathcal E_{l}^e) \ne\emptyset$ for any
$w\in \Omega_\ast$.
 Since $\mathbb G$ is connected, for any $\alpha\in \A,$
there exists an admissible path $u_\alpha:=\alpha \alpha_1\cdots \alpha_{k-1}0_e$. Assume $w\in \Omega_m$,
define $v:=w|_{m-1}u_{w_m}$,
 then $I_v$ has type $0_e$,
$I_{v}\subset I_w$ and $\pi(v3_{or}(0_e3_{or})^\infty)\in I_v\cap \pi(\mathcal E_{l}^e).$ Hence $I_w\cap \pi(\mathcal E_{l}^e) \ne\emptyset.$

Assume $E\in \pi(\mathcal E_{l}^e)$, then there exist
 some even $n$  and $w\in\Omega_n$
with $w_n=0_e$ such that $E=\pi(w(3_{or}0_e)^\infty)$. write $\sigma=\Pi_\ast(w)$, then $\sigma\in \Sigma_n$. For any $k\ge0$, by \eqref{Pi-star} and tracing on $\mathbb G$, we have $\Pi_\ast(w(3_{or}0_e)^k)=\sigma 1^{2k}.$
By \eqref{I-and-B},
$$
B_{\sigma 1^{2k}}=I_{w(3_{or}0_e)^k}.
$$
Write $\sigma^{(k)}:=\sigma1^{2k}$. Then
 $B_{\sigma^{(k)}}$ has type $0_e$. By the definition of $\pi$, we have
 $E\in B_{\sigma^{(k)}}$. We claim that $E\ge b_{\sigma 1^{2k+1}}$. Indeed,
 by applying \eqref{cascade} to every $\sigma^{(k)}$, we have
 $$
 b_{\sigma1}<b_{\sigma1^3}<\cdots <b_{\sigma1^{2k+1}}<\cdots<b_{\sigma1^{2k}}<\cdots<b_{\sigma1^2}<b_{\sigma}.
 $$
Since $E\in B_{\sigma^{(k)}}$ and the length of $ B_{\sigma^{(k)}}$ tend to 0, we conclude that
$$
E=\inf_{k\ge 1}b_{\sigma1^{2k}}\ge \sup_{k\ge1} b_{\sigma1^{2k+1}}.
$$

So $E\in [b_{\sigma^{(k)}1},b_{\sigma^{(k)}}]$. By \eqref{bd-h-n+1},
\begin{equation*}\label{n+2k+1}
  h_{n+2k+1}(E)\ge 2.
\end{equation*}

 Since $E\in B_{\sigma^{(k+1)}}=[a_{\sigma^{(k)}11}, b_{\sigma^{(k)}11}]\subset B_{\sigma^{(k)}}=[a_{\sigma^{(k)}}, b_{\sigma^{(k)}}],$ by \eqref{bd-h-n},
\begin{equation*}\label{n+2k}
  h_{n+2k}(E)\ge \sqrt{2}.
\end{equation*}

By Lemma \ref{type-0} ii), $E\in B_{\sigma1^4}\subset {\rm int}(B_{\sigma1^2})$. By Lemma \ref{basic-zeros} ii),  $E\notin\RR_{n+1}$.

  Now for any $N\ge1$, since $E\in B_{\sigma1^{2N}}$,
we have $|h_{n+2N}(E)|\le 2$.
Then by Lemma \ref{fundamental},
\begin{eqnarray*}
  |h_{n+2N+1}(E)| &\ge& 2\lambda \prod_{j=0}^{n+2N}|h_j(E)|-|h_{n+2N}^2(E)-2| \\
   &\ge&  2^{3N/2+1}\lambda\prod_{j=0}^n|h_j(E)| -6.
\end{eqnarray*}
Hence $|h_{n+2N+1}(E)|\to\infty$.
Now by \eqref{recurrence},
$$
|h_{n+2N}^2(E)-2|=\frac{|h_{n+2N+2}(E)+2|}{|h_{n+2N+1}(E)|}\le \frac{4}{|h_{n+2N+1}(E)|}\to0.
$$
So $|h_{n+2N}(E)|\to\sqrt{2}$.
\hfill $\Box$

\subsection{$\infty_{II}$-type energies}\

We will show that  $E\in \pi(\widetilde{\mathcal E}_l\cup \widetilde{\mathcal E}_r)$ is also an $\infty$-type energy. Recall that in Sec. \ref{subsec-gap-zero-infty} we call it  $\infty_{II}$-type  energy.

\begin{prop}\label{infty-energy-II}
i) Both $ \pi(\widetilde {\mathcal E}_{l})$
  and $
  \pi(\widetilde{\mathcal E}_{r})$ are dense in $\sigma_\lambda.$

ii) Assume $(\omega,\hat\omega)\in \G_{II}$, then
\begin{equation*}\label{homeo-1}
  \lim_{n\to\infty}|h_{2n}(\pi(\omega))|=\lim_{n\to\infty}|h_{2n+1}(
  \pi(\hat\omega))|=\infty.
\end{equation*}
\end{prop}


  Although  this proposition looks similar with
 Proposition \ref{infty-energy-I}, the proof of which is much more involved. We need to study
 a dynamical system induced by the recurrence  relation of trace polynomials, which is inspired by
 our previous work \cite{LQY}.

 \begin{figure}[htbp]

\begin{minipage}[t]{0.3\linewidth}
\centering
\includegraphics[width=1.85in]{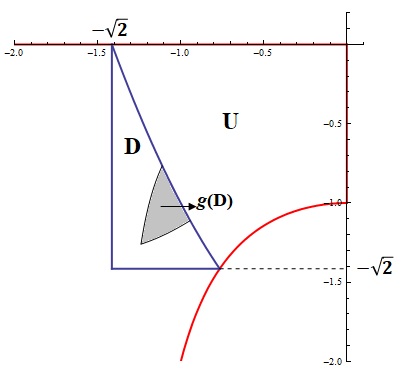}
\caption*{(a)}
\end{minipage}%
\hskip 1.5cm
\begin{minipage}[t]{0.3\linewidth}
\centering
\includegraphics[width=2in]{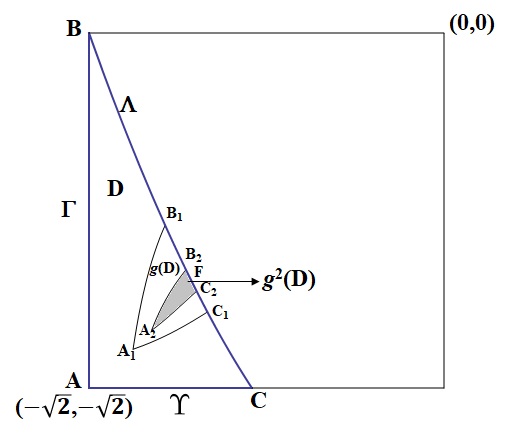}
\caption*{(b)}
\end{minipage}%
\caption{Dynamics of $f$}\label{dyn-f}
\end{figure}

 Define the {\it period-doubling trace  map} $f:\R^2\to\R^2$ as
\begin{equation}\label{f}
  f(x,y)=\left(y(x^2-2)-2, [y(x^2-2)-2](y^2-2)-2\right).
\end{equation}
By \eqref{recurrence}, for any $n\ge0$ and $E\in \R$, we have
\begin{equation*}\label{iteration-f}
  \left(h_{n+2}(E),h_{n+3}(E)\right)=f\left(h_{n}(E),h_{n+1}(E)\right).
\end{equation*}
We have
\begin{equation*}\label{Df}
Df(x,y)=
\begin{bmatrix}
2xy&x^2-2\\
2xy(y^2-2)&(x^2-2)(3y^2-2)-4y
\end{bmatrix}.
\end{equation*}
So we have
\begin{equation*}\label{Jf}
\det Df(x,y)=4xy^2(y(x^2-2)-2).
\end{equation*}
Define a simply connected domain $U$ (enclosed by the red curve in Figure \ref{dyn-f} (a)) as
$$
U:=\{(x,y):x,y<0, y(x^2-2)-2<0\}.
$$
Define a compact set  $D$ (enclosed by the blue curve in Figure \ref{dyn-f} (a)) as
\begin{equation*}\label{Domain-D}
  D:=\{(x,y):  -\sqrt{2}\le x,y\le0; y\le x^2-2\}.
\end{equation*}
Notice that $(-1,-1)\in \partial D.$ Write
$$A:=(-\sqrt{2}, -\sqrt{2}),\  B:=(-\sqrt{2},0),\  C:=(-\sqrt{2-\sqrt{2}},-\sqrt{2}), \ \ F:=(-1,-1).$$
Then $A,B,C$ are the vertices of the ``triangle" $D$. Define
$$
\begin{cases}
\Gamma:=\{(-\sqrt{2},y): -\sqrt{2}\le y\le 0\},\\
\Upsilon:=\{(x,-\sqrt{2}): -\sqrt{2}\le x\le -\sqrt{2-\sqrt{2}}\},\\
\Lambda:=\{(x,x^2-2):-\sqrt{2}\le x\le -\sqrt{2-\sqrt{2}}\}.
\end{cases}
$$
Then $B,C, F\in \Lambda$ and   $\partial D=\Gamma\cup\Upsilon\cup\Lambda$, see Figure \ref{dyn-f} (b).
\subsubsection{Local dynamical property of $f$}\label{sec-loc-pro-f}\


We summarize the needed properties of $f$ in the following proposition.

\begin{prop}\label{f-basic}
i) $f:U\to f(U)$ is a diffeomorphism and $D\subset f(U)$. $f$ has a unique fixed point $F$ in $U.$

ii) Let $g: f(U)\to U$ be the inverse of $f$, then $g$ has a unique fixed point $F$ in $f(U)$. Moreover,  $g(D)\subset D\cap U$.

iii) $g(\Lambda)\subset \Lambda$ and $g^n(\Lambda)\downarrow \{F\}$.

iv) Write $g=(g_1,g_2)$. For any $(x,y)\in D$ and $i=1,2$,
\begin{equation}\label{g-monotone}
  \frac{\partial g_i}{\partial x}(x,y),\ \  \frac{\partial g_i}{\partial y}(x,y)>0.
\end{equation}
\end{prop}

\proof\
i)  Write $(x_1,y_1):=f(x,y)$, we have
$$y_1=x_1\frac{(x_1+2)^2}{(x^2-2)^2}-2x_1-2.$$
Since $x_1<0$ for $(x,y)\in U$, we have $y_1<-2x_1-2$. Moreover, for $-2\le x<0$, $y_1<x_1^3/4+x_1^2-x_1-2$.
Now by computing the image of the fiber $U_a:=U\cap \{y=a\}$ under $f$ for all $a<0$, one can verify that
$$
f(U)=\left\{(x,y): x<0;\ y<-2x-2;\ y<\frac{x^3}{4}+x^2-x-2\ \mbox{ if } -2\le x <0 \right\}.
$$
On check directly that $D\subset f(U)$.
Define a function $g$ as
\begin{equation*}\label{inverse-f}
  g(x,y):=\left(-\sqrt{2-\frac{2+x}{\sqrt{2+\frac{2+y}{x}}}},-\sqrt{2+\frac{2+y}{x}}\right).
\end{equation*}
One  check directly  that  $g$ can be defined on $f(U)$ and is smooth, and
$$g\circ f|_U={\rm Id}_U;\ \ \  f\circ g|_{f(U)}={\rm Id}_{f(U)}.$$
So $f:U\to f(U)$ is a diffeomorphism and $g$ is the inverse of $f$.

By direct computation, the fixed points of $f$ are
$$
(-1,-1),\ \ (2,2),\ \ (-\alpha,\alpha-1),\ \ (\alpha^{-1},-(1+\alpha^{-1})),
$$
where $\alpha=(\sqrt{5}+1)/2$. So $F$ is the only fixed point of $f$ in $U.$

ii) The first assertion is a restatement of i). Fix any $(x, y)\in D$. By direct computation,
  $$
   -\sqrt{2}<-1.3\leq  g_1(x, y)\leq -1.2<-\sqrt{2-\sqrt{2}};\ \ \ -\sqrt{2}<-1.3\leq  g_2(x, y)\leq -0.7.
  $$
   Moreover,
$$g_2(x,y)-g_1^2(x,y)+2=\frac{y-x^2+2}{x g_2(x,y)}\le0.$$
By the definitions of $D$ and $U$, we conclude that $g(D)\subset D\cap U$.

iii) Define $\phi:\R\to\R$ as $\phi(x)=x^2-2.$ Write $a:=-\sqrt{2}, b:=-\sqrt{2-\sqrt{2}}$. Then
$$
\Lambda=\{(x,\phi(x)): x\in [a,b]\}.
$$
By the definition of $f$, for any $x\in\R,$
$$
f(x,\phi(x))=(\phi^2(x),\phi^3(x)),
$$
where $\phi^n$ is the $n$-time iteration of $\phi$.
By solving $\phi^2(x)=a$ and $b$, there exist $a_1, b_1$ such that $$a<a_1=-1.11\cdots<-1<b_{1}<b=-0.76\cdots;\ \  \phi^2(a_1)=a;\ \ \phi^2(b_1)=b.$$
Since $[\phi^2(x)]'=4x(x^2-2)$, for any $x\in[a_1,b_1]$,
\begin{equation}\label{bd-deriv-phi-2}
[\phi^2(x)]'\ge 3.
\end{equation}
Define $\Lambda_1:=\{(x,\phi(x)): x\in [a_1,b_1]\}$, then $\Lambda_1\subset \Lambda$ and $f( \Lambda_1)=\Lambda.$ Since $\Lambda\subset D\subset f(U)$ and $g:f(U)\to U$ is a diffeomorphism, we conclude that
$$
g(\Lambda)= \Lambda_1\subset \Lambda.
$$
Notice that $\phi^2(-1)=-1$.
By \eqref{bd-deriv-phi-2} and contraction principle, for any $n\ge2$ one can find $a_n,b_n$ such that $a_n> a_1, b_n<b_1$
and $a_n\uparrow -1, b_n\downarrow -1$ and $\phi^{2n}(a_n)=a, \phi^{2n}(b_n)=b$. By repeating the above proof, we get
$$g^n(\Lambda)=\Lambda_n:=\{(x,\phi(x)): x\in [a_n,b_n]\}\downarrow \{F\}.$$

iv) By ii), we have $-1.3\le g_1(x, y), g_2(x, y)\le -0.7$.  Thus
\begin{equation*}
\begin{split}
\frac{\partial g_1}{\partial x}&=\frac{1}{2g_1g_2}+\frac{(2-g_{1}^{2})(2-g_{2}^2)}{4xg_1g_{2}^2}>0, \quad
\frac{\partial g_1}{\partial y}=\frac{2-g_{1}^{2}}{4xg_1g_{2}^2}>0;\\
\frac{\partial g_2}{\partial x}&=\frac{2-g_{2}^2}{2xg_2}>0, \quad \frac{\partial g_2}{\partial y}=\frac{1}{2xg_2}>0.
\end{split}
\end{equation*}
So, \eqref{g-monotone} holds.
\hfill $\Box$

Define a partial order on $\R^2$ as follows. Assume
$p=(x,y)$ and $q=(\tilde x,\tilde y)$, we say $p\le q$ if $x\le\tilde x, y\le\tilde y$. Assume $\mathcal C$
is a continuous curve in $\R^2$ with a parametrization $\gamma: [a,b]\to \mathcal C$ such that
 $\gamma(t)\le \gamma(s)$ whenever $t<s$, then we call $\mathcal C$  an {\it increasing curve} in $\R^2.$

\begin{lem}\label{monotone}
If $p,q\in D$ and $p\le q$, then $g(p)\le g(q)$.
\end{lem}

\proof\
Notice that, if $p,q\in D$ and $p\le q$, then for any $\hat p\in\R^2$ with $p\le\hat p\le q$, we have $\hat p\in D$.
 By \eqref{g-monotone},
 $g_1(x,\cdot),$ $ g_1(\cdot,y),$  $g_2(x,\cdot),$ $ g_2(\cdot,y)$ are all strictly increasing.
Now assume $p=(x,y), q=(\tilde x,\tilde y)$ and $p\le q$, then
\begin{eqnarray*}
  g_1(p) &=& g_1(x,y)\le g_1(\tilde x,y)\le g_1(\tilde x,\tilde y)=g_1(q), \\
 g_2(p) &=& g_2(x,y)\le g_2(\tilde x,y)\le g_2(\tilde x,\tilde y)=g_2(q).
\end{eqnarray*}
Hence $g(p)\le g(q)$.
\hfill $\Box$

Now we have the following consequence.

\begin{prop}\label{contraction-g}
Assume $p\in D$. Then
 $f^n(p)\in D$ for any $n\in \N$ if and only if $p=F$.
\end{prop}

\proof\
 Since $f(F)=F$ and $F\in D,$
  we have  $f^n(F)\in  D$ for any $n$. So the if part holds.

  To prove the only if part, we need the following claim:

  \noindent{\bf Claim:} We have
\begin{equation}\label{intersection}
 g^n(D)\downarrow\{(-1,-1)\}.
\end{equation}

\noindent $\lhd$
since $g(F)=F$, we have $F\in g^n(D)$ for any $n\ge0.$
 By Proposition \ref{f-basic} ii), $g(D)\subset D$. So
 $\{g^{n}(D):n\ge 0\}$ is a decreasing  sequence of compact sets. To prove \eqref{intersection}, we only
need to show that the diameter of $g^{n}(D)$ tends to 0.

Recall that, $A,B,C$ are the vertices of the ``triangle" $D$. Define
$$A_n:=g^n(A),\ \  B_n:=g^n(B),\ \  C_n:=g^n(C).$$
Since $B,C\in \Lambda,$ By Proposition \ref{f-basic} iii),  $B_n, C_n\to F.$ We also have
$$
A\le A_1=\left(-\sqrt{2-\frac{2-\sqrt{2}}{\sqrt{3-\sqrt{2}}}},-\sqrt{3-\sqrt{2}}\right),
\ \ \ A, A_1\in D.
$$
 We have $A_n\in D$ since $g(D)\subset D.$  By Lemma \ref{monotone},
$$
A\le A_1\le A_2\le\cdots
$$
Hence $A_\infty:=\lim_n A_n$ exists and $A_\infty\in D$, since $g$ is continuous on $D$,  $g(A_\infty)=A_\infty$. By Proposition \ref{f-basic} ii),  $g$ has a unique fixed
point $F$ in $D$,  so $A_\infty=F.$

 Recall that $\Gamma\cup \Upsilon\cup \Lambda$ is the boundary of $D$. By Proposition \ref{f-basic} i) and ii), $g: f(U)\to U$ is a diffeomorphism
and $D\subset f(U)$.  Since $g^n(D)\subset D,$ we conclude that
 $g^n(\Gamma)\cup g^n(\Upsilon)\cup g^n(\Lambda)$ is the boundary of $g^n(D)$. Notice that
 $\Gamma$ is the vertical
 interval $AB$ which is an increasing interval. By Lemma \ref{monotone},
 $g^n(\Gamma)$ is also increasing and has endpoints $A_n, B_n$. Similarly,
 $g^n(\Upsilon)$ is  increasing and has endpoints $A_n, C_n$. By Proposition \ref{f-basic} iii),
 $g^n(\Lambda)\subset \Lambda$ is the segment of  $\Lambda$ with endpoints $B_n,C_n$.
 So we conclude that
 $$
 g^n(D)\subset [x_{A_n}, x_{C_n}]\times [y_{A_n}, y_{B_n}].
 $$
 Since $A_n, B_n, C_n\to F$, we conclude that ${\rm diam}(g^n(D))\to0$.
 \hfill $\rhd$

 By the claim, we can write $ D\setminus \{F\}$ as the following disjoint union:
 $$
 D\setminus \{F\}=\bigcup_{n\ge1}(g^{n-1}(D)\setminus g^{n}(D)).
 $$
 Now if $p\in D\setminus \{F\}$, then there exists  $m\ge1$ such that
 $p\in g^{m-1}(D)\setminus g^{m}(D).$ So we have
 $$
 f^{m}(p)\in f(f^{m-1}(g^{m-1}(D)\setminus g^{m}(D)))=f(D\setminus g(D)).
 $$
 By Proposition \ref{f-basic} ii), $g(D)\subset D\cap U$. Notice that $B,C\notin U$ and
  $D\setminus \{B,C\}\subset U$, so  $g(D)\subset D\setminus\{B,C\}\subset U. $ Since $f$ is injective on $U$, we conclude that
  $$
  f((D\setminus\{B,C\})\setminus g(D))=f(D\setminus \{B,C\})\setminus D.
  $$
  Hence $f((D\setminus\{B,C\})\setminus g(D))\cap D=\emptyset.$ On the other hand, by direct computation,
  $$
  f(B)=(-2,2), f(C)=(0,-2)\notin D.
  $$
  So we have
  $$
  D\cap f(D\setminus g(D))=D\cap \left(f((D\setminus\{B,C\})\setminus g(D))\cup\{f(B),f(C)\}\right)=\emptyset.
  $$
That is, $f^m(p)\notin D$. So  the only if part follows.
\hfill $\Box$

\subsubsection{Proof of Proposition \ref{infty-energy-II}}\

Assume $(\omega,\hat\omega)\in \G_{II}$. Write $E_1=\pi(\omega)$ and
$ E_2=\pi(\hat\omega)$. By Lemma \ref{order-gap-II-weak}, $E_1\le E_2.$

\begin{lem}\label{bndy}
 If $E_1=E_2=:E$, then there exists $m\in\N$ such that
 for  any $k\ge m$,
  $$(h_k(E),h_{k+1}(E))\in {\rm int}(D).$$

 If $E_1<E_2$, then as $k\to \infty$,
 $$
 |h_{2k}(E_1)|, |h_{2k+1}(E_2)|\to\infty.
 $$
\end{lem}

\proof\
By \eqref{def-elll}, $\omega$ takes one of the following two forms:
 $$
 w1_o3_{er}(0_o3_{er})^\infty;\ \ \  w3_{el}(0_o3_{er})^\infty.
 $$

We prove the lemma for $\omega=w1_o3_{er}(0_o3_{er})^\infty$, and leave the other case to the reader.

 Assume $w\in\Omega_n$. In this case,  $n$ is even, $w_{n}=0_e$ or $2_e$ and $\hat\omega=\ell(\omega)=w3_{or}(0_e3_{ol})^\infty$. For $k\ge1$, by \eqref{gap-1},
\begin{equation*}\label{}
  \begin{cases}
  \Pi_\ast(\omega|_{n+2k+1})=\Pi_\ast(w)\delta (10)^k=:\sigma^{(k)},\\
  \Pi_\ast(\hat\omega|_{n+2k+2})=\Pi_\ast(w)\delta1 (01)^{k}=:\hat\sigma^{(k)}.
  \end{cases}
\end{equation*}
Here $\delta=1$ if $w_n=0_e$ and $\delta=0$ if $w_n=2_e.$
 Notice that $\hat\sigma^{(k)}=\sigma^{(k)}1$.
By \eqref{I-and-B},
 $$I_{\omega|_{n+2k+1}}=I_{w1_o(3_{er}0_o)^k}=B_{\sigma^{(k)}},\ \ \  I_{\hat\omega|_{n+2k+2}}=I_{w3_{or}0_e(3_{ol}0_e)^k}=B_{\hat\sigma^{(k)}}.$$
%
%
We have $E_1\in B_{\sigma^{(k)}}$ and $E_2\in B_{\hat\sigma^{(k)}}$. Notice that $B_{\sigma^{(k)}}$ has type $0_o$ and $B_{\hat\sigma^{(k)}}$ has type $0_e$.
%
%
By
applying Lemma \ref{type-0} i) to $B_{\sigma^{(k)}}$,  we have
\begin{equation}\label{E-1}
z_{\sigma^{(k)}}=b_{\sigma^{(k)}0}<a_{\sigma^{(k+1)}}<E_1<b_{\sigma^{(k+1)}}<b_{\sigma^{(k)}}.
\end{equation}
In particular, $E_1\in {\rm int}(B_{\sigma^{(k)}})$.
By Proposition \ref{basic-facts} i), $h_{n+2(k+1)+1}$ is decreasing on $B_{\sigma^{(k+1)}}$,
so $h_{n+2(k+1)+1}(b_{\sigma^{(k+1)}})=-2.$ Hence by \eqref{recurrence},
$$
-2=h_{n+2(k+1)}(b_{\sigma^{(k+1)}})(h_{n+2k+1}^2(b_{\sigma^{(k+1)}})-2)-2.
$$
Since $B_{\sigma^{(k+1)}}$ has type $0_o$, $b_{\sigma^{(k+1)}}\notin\RR_{n+2(k+1)}$. Thus $h_{n+2k+1}(b_{\sigma^{(k+1)}})=\pm\sqrt{2}.$ Since
$h_{n+2k+1}$ is decreasing on $B_{\sigma^{(k)}}$, by \eqref{E-1},
\begin{equation}\label{72}
-\sqrt{2}=h_{n+2k+1}(b_{\sigma^{(k+1)}})<h_{n+2k+1}(E_1)<h_{n+2k+1}(z_{\sigma^{(k)}})=0.
\end{equation}

On the other hand, since $E_1\in {\rm int}(B_{\sigma^{(k)}})$
 with $|\sigma^{(k)}|=n+2k+1$ and $\sigma^{(k)}$ ends with $0$, by Lemma \ref{fundamental} and Corollary
  \ref{sgn-prod},
 \begin{equation}\label{below-para}
h_{n+2k+1}(E_1)-(h_{n+2k}^2(E_1)-2)=(-1)^{n+2k}2\lambda\prod_{j=0}^{n+2k} h_j(E_1)<0.
\end{equation}

By a symmetric argument, one can show that
\begin{equation}\label{73}
-\sqrt{2}<h_{n+2k+2}(E_2)<0 \ \ \text{ and }\ \  h_{n+2k+2}(E_2)-(h_{n+2k+1}^2(E_2)-2)<0.
\end{equation}

Now assume  $E_1=E_2=:E$. By \eqref{72}, \eqref{below-para} and \eqref{73}, for any $k\ge n+3=:m$,
$$-\sqrt{2}< h_k(E)<0;\ \ \ \  h_{k+1}(E)<h_k^2(E)-2.$$
By the definition of $D$, we have $(h_k(E),h_{k+1}(E))\in {\rm int}(D).$

Next assume $E_1<E_2$. By \eqref{E-1}, we have $z_{\sigma^{(k)}}<E_1$. By a symmetric argument, we have $E_2<z_{\hat \sigma^{(k)}}.$ So
\begin{equation*}\label{E-1<E-2}
z_{\sigma^{(k)}}<E_1<E_2<z_{\hat \sigma^{(k)}}.
\end{equation*}
Since the length of $B_{\sigma^{(k)}}$ and $B_{\hat\sigma^{(k)}}$ tend to $0$, There exists $K\ge n+3$
such that for any $k\ge K$, $B_{\sigma^{(k)}}< B_{\hat \sigma^{(k)}}.$
Recall that $\hat\sigma^{(k)}=\sigma^{(k)}1$.  So $b_{\sigma^{(k)}}<a_{\hat \sigma^{(k)}}=a_{\sigma^{(k)}1}$. By Proposition \ref{basic-facts} i), $h_{n+2k+2}$ is decreasing on $B_{\sigma^{(k)}0}$ and increasing on $B_{\sigma^{(k)}1}$. By Floquet theory,   for any
$E\in (b_{\sigma^{(k)}0},a_{\sigma^{(k)}1})$, we have $h_{n+2k+2}(E)<-2.$  By \eqref{E-1}, $b_{\sigma^{(k)}0}<E_1<b_{\sigma^{(k)}}<a_{\sigma^{(k)}1}$, so
\begin{equation}\label{E_1}
  h_{n+2k+2}(E_1)<-2.
\end{equation}

We claim that
\begin{equation}\label{E_1-1}
  h_{n+2k+1}(E_1)\le -1.
\end{equation}
If otherwise, by \eqref{72}, we have $-1<h_{n+2k+1}(E_1)<0.$ Then
$$
h_{n+2k+3}(E_1)=h_{n+2k+2}(E_1)(h_{n+2k+1}^2(E_1)-2)-2\ge0,
$$
which contradicts with \eqref{72}.

By \eqref{below-para}, $E_1\notin\ZZ$. Since $E_1\in B_{\sigma^{(k)}}$, we have $|h_{n+2k+1}(E_1)|\le2.$  Now for any $N> K$, by
Lemma \ref{fundamental} and  \eqref{E_1}, \eqref{E_1-1},
\begin{eqnarray*}
  |h_{n+2N+2}(E_1)| &\ge& 2\lambda \prod_{j=0}^{n+2N+1}|h_j(E_1)|-|h_{n+2N+1}^2(E_1)-2| \\
   &\ge& \lambda 2^{(N-K)/2}\prod_{j=0}^{n+2K}|h_j(E_1)| -2.
\end{eqnarray*}
Hence $|h_{n+2N+2}(E_1)|\to\infty$ as $N\to\infty$.

By exactly the same proof, we can show that $|h_{n+2N+1}(E_2)|\to\infty$ as $N\to\infty$.
So the proof is finished.
\hfill $\Box$

As an application, we can improve Lemma \ref{order-gap-II-weak} as follows:

 \begin{lem}\label{order-gap-II}
 If $(\omega,\hat\omega)\in \G_{II}$, then
 $
 \pi(\omega)< \pi(\hat \omega).
$
 \end{lem}

\proof\
 We prove it by contradiction.
Assume $\pi(\omega)=\pi(\hat \omega)=:E$. By Lemma \ref{bndy}, there exists  $m\in \N$ such that  for any $k\ge m$,
 $(h_k(E),h_{k+1}(E))\in {\rm int}(D).$ Recall that
$f(h_k(E),h_{k+1}(E))=(h_{k+2}(E),h_{k+3}(E))$, thus for any $k\in \N$, we have
$$
f^k(h_m(E),h_{m+1}(E))=(h_{m+2k}(E),h_{m+2k+1}(E))\in {\rm int}(D).
$$
By Proposition \ref{contraction-g}, $(h_m(E),h_{m+1}(E))=F$. However,  $F\in \partial D$, a contradiction.
So  $\pi(\omega)<\pi(\hat \omega)$.
\hfill $\Box$

\noindent {\bf Proof  of Proposition \ref{infty-energy-II}}
i) The proof is the same as Proposition \ref{infty-energy-I}.

ii)
It is a direct consequence of Lemma \ref{order-gap-II} and Lemma \ref{bndy}.
\hfill $\Box$

\subsection{Proof of Theorems \ref{main-infty-energy}, \ref{main-coding-spectrum} and \ref{main-gap-zero-inftyenergy}
}\

\noindent {\bf Proof of Theorem \ref{main-infty-energy}.}\
At first we show  $B_\infty\subset \sigma_\lambda.$ Fix $E\in B_\infty$, we claim that for any $n\ge0$, $|h_n(E)|\le 2$ or $|h_{n+1}(E)|\le2$. If otherwise, there exists $n_0\ge 0$ and $\delta>0$ such that
$|h_{n_0}(E)|,|h_{n_0+1}(E)|\ge2+\delta$. By \eqref{recurrence},
$$
|h_{n_0+2}(E)|\ge |h_{n_0+1}(E)|(h_{n_0}^2(E)-2)-2\ge 2+10\delta.
$$
Using \eqref{recurrence} again, we have $|h_{n_0+3}(E)|\ge 2+10\delta$. Now by induction, it is seen that
$$
|h_{n_0+k}(E)|\ge 2+10^{[k/2]}\delta.
$$
Thus $|h_n(E)|\to \infty$, which  contradicts with $E\in B_\infty$. So the claim holds. Now by \eqref{spec-covering}, we conclude that $E\in \sigma_\lambda.$ Thus $B_\infty\subset \sigma_\lambda$.

By Lemma \ref{basic-zeros} iii), $\ZZ\subset B_\infty$ and $B_\infty$ is dense in $\sigma_\lambda.$ By Proposition \ref{SNS} and \eqref{dim-lim-set},
$$
\dim_H B_\infty\ge \dim_H A(\widetilde \BB)\ge \log \alpha/\log 4>0,
$$
so $B_\infty$ is uncountable.

By the definition of $\mathscr{E}_\infty$, $\mathscr{E}_\infty=\sigma_\lambda\setminus B_\infty$ always holds. Now we show that  $\mathscr{E}_\infty$ is dense in $\sigma_\lambda$ and uncountable.

By  Proposition \ref{infty-energy-I}, $\pi(\mathcal E^e_{l})\subset \mathscr{E}_\infty$ and is dense in $\sigma_\lambda.$  So  $\mathscr{E}_\infty$ is dense in $\sigma_\lambda.$

Next we show that $\mathscr{E}_\infty$ is uncountable.
We prove it by contradiction. Otherwise, $\mathscr{E}_\infty$ is countable.
 Let
$\mathscr{E}_\infty=\{\gamma_i\}_{i=1}^{\infty}.
$
Since $\sigma_\lambda$ has no isolated point, we conclude that $\{\gamma_i\}$ is nowhere dense in $\sigma_\lambda$ for each $i$.
Denote
\begin{equation*}
T:=\sigma_\lambda\setminus \mathscr{E}_\infty=\{E\in \sigma_\lambda:
\mathscr O(E) \text{ is bounded } \}.
\end{equation*}
We can write $T$ as a countable union of closed subsets,
\begin{equation*}
T=\bigcup_{k=1}^{\infty}A_k,
\end{equation*}
where
\begin{equation*}
A_k=\bigcap_{n=0}^{\infty}\{E\in\sigma_\lambda: |h_{n}(E)|\leq k\}.
\end{equation*}
Since $\mathscr{E}_\infty$ is dense in $\sigma_\lambda$ and $\mathscr{E}_\infty \cap A_k=\emptyset$,
$A_k$ is nowhere dense for each $k\in \mathbb{N}$. Thus
\begin{equation*}
\sigma_\lambda=\bigcup_{i=1}^{\infty}\{\gamma_i\}\cup\bigcup_{k=1}^{\infty}A_k,
\end{equation*}
i.e., $\sigma_\lambda$ is a countable union of  nowhere dense closed sets.
This contradicts with  Baire's category theorem.
\hfill $\Box$

We remark that, once we know that $\mathscr{E}_\infty$ is dense in $\sigma_\lambda$, the rest  proof of the theorem  is purely topological.

\noindent {\bf Proof of Theorem \ref{main-coding-spectrum}.}\
By Proposition \ref{main-coding-spectrum-weak}, we only need to show that $\pi$ is a homeomorphism. Since $\Omega_\infty$ is compact and $\sigma_\lambda$ is Hausdorff, it suffice to show that $\pi$ is injective.

 Assume $\omega\ne \hat\omega$. WLOG, we assume $\omega\prec \hat\omega.$ If $\omega<\hat\omega,$ then by Lemma \ref{order-<-pre}, $\pi(\omega)<\pi(\hat\omega)$. If $\omega\not<\hat\omega,$ by Lemma \ref{prec-not-<}, there exists $(\tau,\hat\tau)\in \G_{II}$ such that $\omega\preceq \tau$ and $\hat\tau \preceq \hat\omega$. By Proposition \ref{main-coding-spectrum-weak} and Lemma \ref{order-gap-II},
 $$
 \pi(\omega)\le \pi(\tau)<\pi(\hat \tau)\le \pi(\hat\omega).
 $$
 Hence $\pi$ is injective.
\hfill $\Box$

\noindent {\bf Proof of Theorem \ref{main-gap-zero-inftyenergy}.}\
i)-iv)
By Theorem \ref{main-coding-spectrum}, $(a,b)$ is a gap of $\sigma_\lambda$ if and only if $(\pi^{-1}(a),\pi^{-1}(b))$ is a gap of $\Omega_\infty.$  So i)-iv)
  follows directly from  i)-iv) of Theorem \ref{symbolic-gap} and Theorem \ref{main-coding-spectrum}.

v) At first, we prove $\pi(\E_l^o)=\ZZ^o$.
Recall that $\ZZ^o=\{z_\sigma: \sigma\in \Sigma_\ast^o\}$. Fix any $\sigma\in \Sigma_\ast^o$, by \eqref{bij-I-type}, there exists a unique $\omega_\sigma\in \E_l^o$ such that $\Pi(\omega)=\sigma01^\infty$. Moreover
$$\{\omega_\sigma:\sigma\in \Sigma_\ast^o\}=\E_l^o.$$
 We claim that $\pi(\omega_\sigma)=z_\sigma.$ Once we show the claim, we immediately get $\pi(\E_l^o)=\ZZ^o.$ Now we show the claim. At first, by Lemma \ref{endsymbol} i), we have $z_\sigma=b_{\sigma01^k}$ for any $k\ge0$. On the other hand, by the definition of $\Pi$, there exists $N\in \N$ such that when $m\ge N$, $\Pi_\ast(\omega_\sigma|_m)=\sigma01^{l_m}$. Note the $\{I_{\omega|_m}: m\ge1\}$ is a decreasing set sequence, so we have
 $$
\pi(\omega_\sigma)=\bigcap_{m\ge0}I_{\omega_\sigma|_m}=\bigcap_{m\ge N}I_{\omega_\sigma|_m}
=\bigcap_{m\ge N}B_{\sigma01^{l_m}}.
$$
Consequently $z_\sigma\in \bigcap_{m\ge N}B_{\sigma01^{l_m}}$. Thus $\pi(\omega_\sigma)=z_\sigma.$

The proof of $\pi(\E_r^e)=\ZZ^e$ is the same. 

vi) It is a direct consequence of Propositions \ref{infty-energy-I} and \ref{infty-energy-II}.
\hfill $\Box$


\section{IDS and gap labelling }\label{sec-ids-gap}

In this section, we prove Theorems \ref{main-comp-IDS} and \ref{main-gap-label}.

\subsection{Computation of IDS}\

\subsubsection{The IDS for $\ZZ$}\

 Recall that   $\ZZ=\{z_\sigma:\sigma\in \Sigma_\ast\}$; $\varpi_n$ is defined by \eqref{varpi-n}; an order $\preceq$ on $\Sigma_\ast$ is defined by \eqref{order-Sigma-ast}. We have

\begin{lem}\label{IDS-R}
For  $\sigma\in \Sigma_n$, we have
\begin{equation}\label{IDS-zero}
N(z_\sigma)=\sum_{i=1}^{n}\frac{\sigma_i}{2^{i}}+\frac{1}{2^{n+1}}.
\end{equation}
Consequently,
\begin{equation}\label{dyadic}
 N(\ZZ)= \{N(z_\sigma): \sigma\in \Sigma_\ast\}
 =\left\{\frac{j}{2^n}: n\ge1; 1\le j<2^n\right\}=\mathscr D\cap (0,1).
\end{equation}
\end{lem}

\begin{proof}
By the property of IDS and Proposition \ref{zero-order}, we have
\begin{eqnarray*}
N(z_\sigma)&=&\lim_{m\to\infty}\frac{\#\{z\in \ZZ_m: z\le z_\sigma\}}{2^m}
=\lim_{m\to\infty}\frac{\#\{\alpha\in\Sigma_m: \alpha\preceq \sigma\}}{2^m}\\
&=&\lim_{m\to\infty}\frac{\#\{\alpha\in\Sigma_m: \alpha\preceq \sigma01^{m-n-1}\}}{2^m}
=\lim_{m\to\infty}\frac{\varpi_m(\sigma01^{m-n-1})+1}{2^m}\\
&=&\lim_{m\to\infty}\left(\sum_{j=1}^n\frac{\sigma_j}{2_j}+\sum_{j=n+2}^m 2^{-j}+2^{-m}\right)
=\sum_{j=1}^n\frac{\sigma_j}{2_j}+\frac{1}{2^{n+1}}.
\end{eqnarray*}
Hence we have
$$\{N(z_\sigma): \sigma\in \Sigma_0\}=\{1/2\};
\ \ \{N(z_\sigma): \sigma\in \Sigma_1\}=\{1/4,3/4\};\ \ \cdots $$
Then \eqref{dyadic} holds.
\hfill $\Box$
\end{proof}

\subsubsection{The IDS for $\sigma_\lambda$}\

Now we compute the IDS for any energy in the spectrum.

\noindent{\bf Proof of Theorem \ref{main-comp-IDS}.}\
 By Proposition \ref{Pi-basic} iii), $\Pi$ is surjective.

Now fix any $\omega\in \Omega_\infty.$ Write $E=\pi(\omega)$ and $\sigma=\Pi(\omega).$
By \eqref{Pi-L}, $\Pi_\ast(\omega|_{n})\lhd \sigma$ for any $n$ and
 $k_n:=|\Pi_\ast(\omega|_{n})|=n$
or $n+1$. So
 $$E=\pi(\omega)=\bigcap_{n\ge0} I_{\omega|_n}=\bigcap_{n\ge 0}B_{\sigma|_{k_n}}.$$
 By \eqref{band-length}, $|B_{\sigma|_{k_n}}|\to 0$ as $n\to\infty$.
 Since $z_{\sigma|_{k_n}}\in B_{\sigma|_{k_n}}$, we have $z_{\sigma|_{k_n}}\to E$.
 Since $N$ is continuous,  by \eqref{IDS-zero}, we have
 $$
 N(E)=\lim_{n\to\infty} N(z_{\sigma|_{k_n}})=
 \lim_{n\to\infty}\sum_{i=1}^{k_n}\frac{\sigma_i}{2^{i}}+\frac{1}{2^{k_n+1}}
 =\sum_{i=1}^{\infty}\frac{\sigma_i}{2^{i}}=\varepsilon(\sigma).
 $$
 That is, $N(\pi(\omega))=\varepsilon(\Pi(\omega)).$
 \hfill $\Box$

\subsection{Gap labelling of $\sigma_\lambda$}\

\noindent{\bf Proof of Theorem \ref{main-gap-label}.}\
By  Theorem \ref{main-gap-zero-inftyenergy} ii), v) and \eqref{dyadic}, we have
\begin{eqnarray*}
  \{N(G): G\in \GG_I\} =\{N(\pi(\omega)): \omega\in \E_l^o\cup\E_r^e\}=\{N(z): z\in \ZZ\}=\mathscr D\cap(0,1).
\end{eqnarray*}
It is well known that the restriction of $\varepsilon$ on $\Sigma_\infty^{(2)}$ is injective and
\begin{equation}\label{D/3}
  (\mathscr D/{3}\setminus \mathscr D)\cap(0,1)=\varepsilon(\Sigma_\infty^{(2)}).
\end{equation}
By  Theorem \ref{main-gap-zero-inftyenergy} iii), Theorem \ref{main-comp-IDS}, Proposition \ref{Pi-basic}  v), injectivity of $\varepsilon$ on $\Sigma_\infty^{(2)}$ and \eqref{D/3}, we have
\begin{eqnarray*}
  \{N(G): G\in \GG_{II}\} &=&\{N(\pi(\omega)): \omega\in \widetilde{\E}_l\}
  =\{\varepsilon(\Pi(\omega)): \omega\in \widetilde{\E}_l\}=\varepsilon(\Pi(\widetilde{\E}_l))\\
  &=&\varepsilon(\Sigma_\infty^{(2)}\setminus\Pi(\mathcal F))
  =\varepsilon(\Sigma_\infty^{(2)})\setminus\varepsilon(\Pi(\mathcal F))\\
  &=&[(\mathscr D/{3}\setminus \mathscr D)\cap(0,1)]\setminus\varepsilon(\Pi(\mathcal F)).
\end{eqnarray*}
So the result follows.
\hfill $\Box$

\section{Proof of technical lemmas and related results }\label{sec-pf-tech-lemma}

In this section, we give the proofs of the technical lemmas and related results in Sec. \ref{sec-optimal-covering} and Sec. \ref{sec-type-evo-order}.

\subsection{The properties of $\ZZ$}\

\subsubsection{The order of $\ZZ$}\label{sec-order-Z}\

The following lemma tells us how the zeros  in $\RR_n$ are
ordered.

\begin{lem}\label{zeros}
$\#\RR_n=2^{n+1}-1$. Moreover $\ZZ_n$ interlaces $\RR_{n-1}$. That is,
if we list the elements of
$\RR_{n-1}$ as
$$
r_1^{(n-1)}<r_2^{(n-1)}<\cdots <r_{2^n-1}^{(n-1)},
$$
then $\RR_n$ is ordered as
\begin{equation}\label{order-R-n}
z_{\varpi_n^{-1}(0)}<r_1^{(n-1)}<z_{\varpi_n^{-1}(1)}<r_2^{(n-1)}<\cdots<
z_{\varpi_n^{-1}(2^n-2)}<r_{2^n-1}^{(n-1)}<z_{\varpi_n^{-1}(2^n-1)}^{(n-1)}.
\end{equation}
In particular, for any $\sigma\in \Sigma_n$ we have
\begin{equation}\label{B_sigma}
B_{\sigma}\subset [r_{\varpi_n(\sigma)}^{(n-1)},r_{\varpi_n(\sigma)+1}^{(n-1)}].
\end{equation}
Here we use the convention:
 $r_0^{(n-1)}=-\infty$ and $r_{2^n}^{(n-1)}=\infty$.

\end{lem}

\begin{proof}
Notice that $\#\mathcal{Z}_n=2^n$.
 By Lemma \ref{basic-zeros} i),
 $\mathcal{Z}_n\cap \mathcal{Z}_m=\emptyset$ for $n\ne m.$ Hence,
$$
\#\mathcal{R}_n=\sum_{j=0}^n\#\mathcal{Z}_j=2^{n+1}-1.
$$

We prove \eqref{order-R-n} by induction. By \eqref{recurrence},
$$z_{\emptyset}=\lambda;\ \ \ z_{0}=-\sqrt{2+\lambda^2};\ \ \  z_{1}=\sqrt{2+\lambda^2}.
$$
So
$
\RR_0=\mathcal{Z}_0=\{z_{\emptyset}\} $ and $ \mathcal{Z}_1=\{z_{0}, z_{1}\}.
$
Since  $z_{0}<z_{\emptyset}<z_{1},$
\eqref{order-R-n} holds for $n=1.$

Now assume the result holds for $n<k$. We consider $\mathcal{R}_k.$
By Lemma \ref{basic-zeros} i) and  induction hypothesis,
\begin{equation}\label{odd-even}
\begin{cases}
\{r_2^{(k-1)},r_4^{(k-1)},\cdots,r_{2^k-2}^{(k-1)}\}=&{\mathcal  R}_{k-2}\subset \{h_k=2\};  \\ \{r_1^{(k-1)},r_3^{(k-1)},\cdots,r_{2^k-1}^{(k-1)}\}=&\ZZ_{k-1}\subset \{h_k=-2\}.
\end{cases}
\end{equation}
Thus in each interval $(r_i^{(k-1)},r_{i+1}^{(k-1)})$ there is an element of $\ZZ_k.$
Moreover since $h_k(x)\to+\infty$ when $|x|\to\infty$, there is one
element of $\ZZ_k$ in $(r_0^{(k-1)},r_1^{(k-1)})$ and one element of $\ZZ_k$ in $(r_{2^k-1}^{(k-1)},r_{2^k}^{(k-1)})$.
Since  $\#\ZZ_k=2^k$, we get all elements of $\ZZ_k.$ Thus the result holds for $n=k$.

By induction, \eqref{order-R-n} holds.

By \eqref{order-R-n},   $z_\sigma\in [r_{\varpi_n(\sigma)}^{(n-1)},r_{\varpi_n(\sigma)+1}^{(n-1)}]$. If $\sigma\ne0^n$ or $1^n$, by
\eqref{odd-even}, $h_n(r_{\varpi_n(\sigma)}^{(n-1)})=\pm2$ and $r_{\varpi_n(\sigma)+1}^{(n-1)}=\mp2.$  If $\sigma=0^n,$  $h_n(r_0^{(n-1)})=\infty$ and $h_n(r_1^{(n-1)})=-2$. If $\sigma=1^n,$ $h_n(r_{2^n}^{(n-1)})=\infty$ and $h_n(r_{2^n-1}^{(n-1)})=-2$. By Proposition \ref{basic-facts} i), \eqref{B_sigma} holds.
\hfill $\Box$
\end{proof}

\begin{cor}\label{order-zero}
 For any $n,t\ge0$ and $\sigma\in \Sigma_n$, we have
\begin{equation}\label{sigma10t}
z_{\sigma01^t}<z_\sigma<z_{\sigma10^t}.
\end{equation}
\end{cor}

\proof\
Fix  any $\sigma\in\Sigma_n$.

 \noindent {\bf Claim:} $z_{\sigma0}<z_\sigma<z_{\sigma1}$ are three consecutive numbers in $\RR_{n+1}$.

 \noindent $\lhd$ By \eqref{order-R-n}, $z_\sigma$ is the $2\varpi_n(\sigma)+1$-th number in $\RR_n.$ So
\begin{equation}\label{z-sigma}
z_\sigma=r^{(n)}_{2\varpi_n(\sigma)+1}=r^{(n)}_{\varpi_{n+1}(\sigma0)+1}=r^{(n)}_{\varpi_{n+1}(\sigma1)}.
\end{equation}
By applying \eqref{order-R-n} for $\RR_{n+1}$, we conclude that
$$z_{\sigma0}<r^{(n)}_{\varpi_{n+1}(\sigma1)}=z_\sigma<z_{\sigma1}$$
are three consecutive numbers in $\RR_{n+1}$.
 \hfill $\rhd$

 Applying the claim to $\sigma0$, we know that $z_{\sigma0}<z_{\sigma01}$ are two consecutive numbers in $\RR_{n+2}$. Since $z_\sigma\in \RR_{n+2}$ and $z_{\sigma0}<z_\sigma$, we have $z_{\sigma0}<z_{\sigma01}<z_\sigma$. A symmetric argument shows that $z_\sigma<z_{\sigma10}$. We can continue this process. By an inductive  argument, we get \eqref{sigma10t}.
\hfill $\Box$

\noindent {\bf Proof of Proposition \ref{zero-order}.}\ Since the map  $\sigma\to z_\sigma$ is a bijection between $\Sigma_\ast$ and $\ZZ$, we  only need to show that if $\sigma\prec \tau$, then $z_\sigma<z_\tau.$

Assume $\sigma\prec \tau$ and $|\sigma|=m, |\tau|=n.$ Write $\theta=\sigma\wedge\tau.$

If $|\theta|<m,n$, then $\sigma=\theta0\ast\le \theta01^{m-1-|\theta|}\in \Sigma_m$ and $\tau=\theta1\ast\ge\theta10^{n-1-|\theta|}\in \Sigma_n$. By \eqref{coding-Z-n} and  Corollary \ref{order-zero},
$$
z_\sigma\le z_{\theta01^{m-1-|\theta|}}<z_\theta<z_{\theta10^{n-1-|\theta|}}\le z_\tau.
$$
If $|\theta|=m,$ then $\sigma\lhd\tau$. Since $\sigma\prec \tau$, we have $\tau=\sigma1\ast\ge \sigma10^{n-1-m}$. Again by
\eqref{coding-Z-n} and  Corollary \ref{order-zero},
$$
z_\sigma<z_{\sigma10^{n-1-m}}\le z_\tau.
$$
If $|\theta|=n$, the same proof shows that $z_\sigma<z_\tau.$
\hfill $\Box$

\noindent {\bf Proof of Corollary \ref{sgn-prod}.}\
By \eqref{B_sigma},
$${\rm int}(B_\sigma)\subset (r_{\varpi_n(\sigma)}^{(n-1)},
 r_{\varpi_n(\sigma)+1}^{(n-1)})=:J_\sigma.$$
 So we only need to show the following claim.

\noindent{\bf Claim}: If $\sigma_n=0(=1)$, then $\prod_{j=0}^{n-1}h_j(E)<0(>0)$ for any $E\in J_\sigma.$

\noindent $\lhd$
We show it by induction. If $n=1$, then $J_0=(-\infty,\lambda) $ and $J_1=(\lambda,\infty)$. Since $h_0(E)=E-\lambda$,
$h_0(E)<0$ for $E\in J_0$ and $h_0(E)>0$ for $E\in J_1$.

Assume the result holds for $n=k$. Now fix any $\sigma\in \Sigma_k$. By \eqref{order-R-n},  for any $j$ we have
$$r^{(k-1)}_j=r^{(k)}_{2j}.$$
Together with \eqref{z-sigma}, we have
\begin{eqnarray*}
J_{\sigma0}&=&(r^{(k)}_{\varpi_{k+1}(\sigma0)},r^{(k)}_{\varpi_{k+1}(\sigma0)+1})
=(r_{\varpi_k(\sigma)}^{(k-1)},z_{\sigma}),\\
J_{\sigma1}&=&(r^{(k)}_{\varpi_{k+1}(\sigma1)},r^{(k)}_{\varpi_{k+1}(\sigma1)+1})
=(z_{\sigma},r_{\varpi_k(\sigma)+1}^{(k-1)}).
\end{eqnarray*}
Take $E\in {\rm int}(J_{\sigma0})$.
If $\sigma_k=0$, then by the induction hypothesis, $\prod_{j=0}^{k-1}h_j(E)<0$.
By Proposition \ref{basic-facts} i), $h_k$ is decreasing in a neighbourhood of $z_\sigma$. Since $z_{\sigma^-}<r_{\varpi_k(\sigma)}^{(k-1)}<z_\sigma$, we conclude that
$h_k(E)>0$.
 If $\sigma_k=1$, then by the induction hypothesis, $\prod_{j=0}^{k-1}h_j(E)>0$.
By a symmetric argument, $h_k(E)<0$. So we always have $\prod_{j=0}^{k}h_j(E)<0$.

Take $E\in {\rm int}(J_{\sigma1})$. By the same argument, we can show that $\prod_{j=0}^{k}h_j(E)>0$.

By induction, the claim holds.
\hfill $\rhd$
\hfill $\Box$

\subsubsection{Local behaviors of trace polynomials at zeros}\

Next we study the  monotone properties of the trace polynomials around a zero.

\begin{lem}\label{deriv-trpoly}

 Assume $n\ge0$ and  $\sigma\in\Sigma_n$.  Then
\begin{equation}\label{sign}
{\rm sgn}(h_{n+1}'(z_\sigma))=(-1)^n;\  \ \ {\rm sgn}(h_{n+k}'(z_\sigma))=(-1)^{n+1},\ (k\ge2).
\end{equation}

\end{lem}

\proof\
 If $n=0$, then $z_\sigma=z_\emptyset=\lambda$ and $h_{n+1}'(z_\sigma)=h_1'(\lambda)=2\lambda>0$.

  For $n\ge1$, by
 the recurrence relation \eqref{recurrence},
 we have,
\begin{equation}\label{0}
h_{n+1}^\prime(E)=h_{n}^\prime(E)(h_{n-1}^2(E)-2)+2h_{n}(E)h_{n-1}(E)h_{n-1}^\prime(E).
\end{equation}
Take $E=z_\sigma$. Since  $h_{n}(z_\sigma)=0$,  we get
\begin{equation}\label{1}
h_{n+1}^\prime( z_\sigma)=h_{n}^\prime(z_\sigma)(h_{n-1}^2(z_\sigma)-2).
\end{equation}
By Lemma \ref{fundamental}, we have
\begin{equation}\label{2}
-(h_{n-1}^2( z_\sigma)-2)=h_{n}(z_\sigma)-(h_{n-1}^2 (z_\sigma) -2)=(-1)^{n-1}2\lambda\prod_{c\in{\mathcal
R}_{n-1}}(z_\sigma-c).
\end{equation}
By Proposition \ref{basic-facts} i),
$${\rm sgn}(h_n'(z_\sigma))=(-1)^{\sigma_n+1}.$$
By  \eqref{order-R-n} and \eqref{varpi-n},
$$
{\rm sgn}\left(\prod_{c\in{\mathcal
R}_{n-1}}(z_\sigma-c)\right)=(-1)^{\#\RR_{n-1}-\varpi_n(\sigma)}=(-1)^{\sigma_n+1}.
$$
 Combining \eqref{1} and \eqref{2}, the sign of $h_{n+1}'(z_\sigma)$ is $(-1)^n$.
So the first equation of \eqref{sign} holds.

By Lemma \ref{basic-zeros} i), $h_{n+1}(z_\sigma)=-2, h_{n+k}(z_\sigma)=2$ for $k\ge2.$ Then by
 \eqref{0},  for any $n\ge0$,
\begin{equation*}
\begin{cases}
h_{n+2}^\prime(z_\sigma)=-2h_{n+1}^\prime(z_\sigma);\\
h_{n+3}'(z_\sigma)=2h_{n+2}'(z_\sigma)-8h_{n+1}^\prime(z_\sigma);\\
h_{n+k}'(z_\sigma)=2h_{n+k-1}'(z_\sigma)+8h_{n+k-2}^\prime(z_\sigma), \ \ (k\ge4).
\end{cases}
\end{equation*}
From this, we conclude that, there exists an increasing positive integer sequence $\{\tau_k:k\ge 2\}$ such that $\tau_2=2, \tau_3=12,\tau_{k+2}=2\tau_{k+1}+8\tau_{k}, (k\ge2)$ and
\begin{equation*}
h_{n+k}'(z_\sigma)=-\tau_k  h_{n+1}^\prime(z_\sigma), \ \  k\ge 2.
\end{equation*}
So the second equation of \eqref{sign} holds.
\hfill $\Box$

\subsubsection{The relation between $\ZZ$ and band edges}\label{sec-Z-band-edge}\

\noindent {\bf Proof of Lemma \ref{endsymbol}.}\
i)
By Lemma \ref{basic-zeros} i),
\begin{equation}\label{-2-and-2}
h_{n+1}(z_\sigma)=-2;\ \  h_{n+1+t}(z_\sigma)=2,\ \ \forall t\ge1.
\end{equation}
By Proposition \ref{basic-facts} iii), $z_\sigma$
must be an endpoint of some band in $\B_{n+1+t}$. By Corollary \ref{order-zero}, $z_{\sigma01^t}<z_\sigma<z_{\sigma10^t}$. Since $(\sigma01^t)^+=\sigma10^t$, $B_{\sigma01^t}$ and $B_{\sigma10^t}$ are two consecutive bands in $\B_{n+1+t}$.  Hence
$z_\sigma$ can only be $b_{\sigma01^t}$ or $a_{\sigma10^t}$.

If $n$ is odd, by Lemma \ref{deriv-trpoly},
$$h_{n+1}'(z_\sigma)<0 \ \ \text{ and }\ \ h_{n+1+t}'(z_\sigma)>0 (t\ge1).$$
By \eqref{-2-and-2} and  Proposition \ref{basic-facts} iii),  $z_\sigma=b_{\sigma01^t}$.

If $n$ is even, by Lemma \ref{deriv-trpoly},
$$h_{n+1}'(z_\sigma)>0 \ \ \text{ and }\ \ h_{n+1+t}'(z_\sigma)<0 (t\ge1).$$
By \eqref{-2-and-2} and Proposition \ref{basic-facts} iii),  $z_\sigma=a_{\sigma10^t}$.

In both cases, by Proposition \ref{basic-facts} ii), $b_{\sigma01^t}<a_{(\sigma01^t)^+}=a_{\sigma10^t}$.

ii)
If $a_\sigma\in \RR_{n-1},$ then by i),
$\sigma=\hat \sigma10^s$ for some $\hat \sigma\in\Sigma_*$, $s\ge0$ and
$
z_{\hat \sigma}=a_{\hat \sigma 10^{t}}
$
for any $t\ge0$.
This implies $a_\sigma=a_{\sigma0^t}$ for any $t>0$.

If $b_\sigma\in \RR_{n-1},$ then by i),
$\sigma=\hat \sigma01^s$ for some $\hat \sigma\in\Sigma_*$, $s\ge0$ and
$
z_{\hat \sigma}=b_{\hat \sigma 01^{t}}
$
for any $t\ge0$.
This implies $b_\sigma=b_{\sigma1^t}$ for any $t>0$.

iii) If $a_\sigma\in \RR_{n-1}$, then by ii), $a_{\sigma0}=a_\sigma\in \RR_{n-1}\subset \RR_n.$ If $a_{\sigma0}\in \RR_{n},$ then by i),
$\sigma0=\hat \sigma10^s$ for some $\hat \sigma\in\Sigma_*$, $s\ge0$ and
$
z_{\hat \sigma}=a_{\hat \sigma 10^{t}}
$
for any $t\ge0$. In this case, we must have $s\ge1$. So we have
$a_{\sigma}=a_{\hat\sigma 10^{s-1}}=z_{\hat\sigma}\in\RR_{n-1}$. Thus the first assertion holds. The same proof shows that the second assertion holds.

iv) By Lemma \ref{zeros}, $\# (z_\sigma,z_{\sigma^+})\cap \RR_{n-1}=1$. By Lemma \ref{basic-zeros} ii), $$\left((z_\sigma,b_\sigma)\cup (a_{\sigma^+},z_{\sigma^+})\right)\cap \RR_{n-1}=\emptyset.$$
By Floquet theory and Lemma \ref{basic-zeros} i),
$$h_n((b_\sigma,a_{\sigma^+}))\subset (-\infty,-2)\cup(2,\infty)\ \ \text{ and }\ \ h_n(\RR_{n-1})\subset\{\pm2\}.$$
So $(b_\sigma,a_{\sigma^+})\cap \RR_{n-1}=\emptyset.$ Combine these facts, the result follows.
\hfill $\Box$

\noindent {\bf Proof of Proposition  \ref{gap-open-sigma-n}.}\
The possible gap of $\sigma_{\lambda,n}$ has the form: $(b_\sigma,a_{\sigma^+})$ with $\sigma\ne 1^n.$
By Lemma \ref{endsymbol} iv),
  $\{b_\sigma,a_{\sigma^+}\}\cap \RR_{n-1}=\{z_\ast\}$. By Lemma \ref{deriv-trpoly}, $h_n'(z_\ast)\ne0$.
   By Proposition \ref{basic-facts} ii), the gap $(b_w,a_{w^+})$ is open.
\hfill $\Box$

\subsection{Proof of Lemma \ref{tech-1}}\label{sec-pf-Lem2.3}\

At first we need the following lemma:

\begin{lem}\label{h-(n+1)}
 Assume $n\ge1$ and $\sigma\in \Sigma_n.$

 i) If $n$ is odd, then
 $$
 \begin{cases}
 h_{n+1}(a_\sigma)>2,& \text{ if } a_\sigma\notin \RR_{n-1};\\
  h_{n+1}(b_\sigma)<2,& \text{ if } b_\sigma\notin \RR_{n-1}.
 \end{cases}
 $$

 ii) If $n$ is even, then
 $$
 \begin{cases}
 h_{n+1}(b_\sigma)>2,& \text{ if } b_\sigma\notin \RR_{n-1};\\
  h_{n+1}(a_\sigma)<2,& \text{ if } a_\sigma\notin \RR_{n-1}.
 \end{cases}
 $$
 \end{lem}

 \begin{proof}
By \eqref{order-R-n},
$$\#(\mathcal{R}_n\cap(-\infty,z_\sigma))=2\varpi_n(\sigma).$$

i)\
At first assume $a_\sigma\notin\RR_{n-1}. $ Since $h_n(a_\sigma)=\pm2,$ and $\RR_n=\RR_{n-1}\cup \ZZ_n$, $a_\sigma\notin \RR_n.$ By Lemma \ref{fundamental}, we have
$$
h_{n+1}(a_\sigma)=h_n^2(a_\sigma)-2+(-1)^n 2\lambda \prod_{c\in\mathcal{R}_n}(a_\sigma-c)=2-2\lambda \prod_{c\in\mathcal{R}_n}(a_\sigma-c).
$$
By Lemma \ref{basic-zeros} ii), $(a_\sigma,z_\sigma)\cap\mathcal{R}_{n}=\emptyset$.
By Lemma \ref{zeros}, we have
\begin{eqnarray*}\#(\mathcal{R}_n\cap(a_\sigma,+\infty))&=&2^{n+1}-1-\#(\mathcal{R}_n\cap(-\infty,a_\sigma])\\
&=&2^{n+1}-1-\#(\mathcal{R}_n\cap(-\infty,z_\sigma))
=2^{n+1}-1-2\varpi_n(\sigma).
\end{eqnarray*}
Thus $\prod_{c\in\mathcal{R}_n}(a_\sigma-c)< 0,$ and consequently $h_{n+1}(a_\sigma)> 2.$

Now assume $b_\sigma\notin\RR_{n-1}$. Then by the same argument, $b_\sigma\notin\RR_n$.  Still by Lemma \ref{fundamental},
$$
h_{n+1}(b_\sigma)=h_n^2(b_\sigma)-2+(-1)^n 2\lambda \prod_{c\in\mathcal{R}_n}(b_\sigma-c)=2-2\lambda \prod_{c\in\mathcal{R}_n}(b_\sigma-c).
$$
By Lemma \ref{basic-zeros} ii), $[z_\sigma,b_\sigma]\cap\mathcal{R}_{n}=[z_\sigma,b_\sigma)\cap\mathcal{R}_{n}=\{z_\sigma\}$.
By Lemma \ref{zeros}, we have
\begin{eqnarray*}
\#(\mathcal{R}_n\cap(b_\sigma,+\infty))&=&2^{n+1}-1-\#(\mathcal{R}_n\cap(-\infty,b_\sigma])\\
&=&2^{n+1}-2-\#(\mathcal{R}_n\cap(-\infty,z_\sigma))
=2^{n+1}-2-2\varpi_n(\sigma).
\end{eqnarray*}
Thus  $\prod_{c\in\mathcal{R}_n}(b_\sigma-c)>0,$
and consequently $h_{n+1}(b_\sigma)<2.$

  ii)\  The proof is the same as i).
 \hfill $\Box$
 \end{proof}


\noindent {\bf Proof of Lemma \ref{tech-1}.}\
i)
By Lemma \ref{endsymbol} i),  $b_{\sigma0}=z_\sigma\in\RR_n$.
  By Lemma \ref{endsymbol} iv),
  $$a_{\sigma1}=a_{(\sigma0)^+}\notin \RR_n.$$

  If $a_\sigma\in \RR_{n-1}$, by Lemma \ref{endsymbol} ii),  $a_{\sigma}=a_{\sigma0}.$
  Hence $B_{\sigma0}=[a_{\sigma0},b_{\sigma0}]=[a_{\sigma},z_\sigma]$.

  If $a_\sigma\notin \RR_{n-1}$, by Lemma \ref{h-(n+1)} i),
   $h_{n+1}(a_\sigma)>2$. By Proposition \ref{basic-facts} i), $h_{n+1}(B_{\sigma0})=[-2,2].$  Since $a_\sigma<z_{\sigma}=b_{\sigma0}$, we conclude that
     $a_{\sigma}<a_{\sigma0}$.
   Hence $B_{\sigma0}\subset (a_{\sigma},z_\sigma]$.

    In both cases, we have $B_{\sigma0}\subset B_\sigma.$

    By Proposition \ref{gap-open-sigma-n}, $b_{\sigma0}<a_{(\sigma0)^+}=a_{\sigma1}$. So
  $z_\sigma=b_{\sigma0}<a_{\sigma1}$.

  If $b_\sigma\in \RR_{n-1}$, by Lemma \ref{endsymbol} ii),  $b_{\sigma}=b_{\sigma1}.$ So  we have $B_{\sigma1}\subset(z_\sigma,b_\sigma].$

If $b_\sigma\notin \RR_{n-1}$, By Lemma \ref{h-(n+1)} i),  $h_{n+1}(b_\sigma)<2$. If $\sigma=1^n$, then $b_{\sigma1}=b_{1^{n+1}}$ is the largest root of $h_{n+1}(x)=2.$ Since $h_{n+1}(b_\sigma)<2$ and $h_{n+1}(x)\to\infty$ when $x\to\infty$, we conclude that $b_\sigma<b_{\sigma1}$. Now assume $\sigma<1^n$.
By Lemma \ref{endsymbol} iv), $a_{\sigma^+}\in \RR_{n-1}.$ By Lemma \ref{endsymbol} ii), $a_{\sigma^+}=a_{\sigma^+0}$. By Proposition \ref{gap-open-sigma-n},
$$
b_\sigma<a_{\sigma^+};\ \ \ b_{\sigma1}<a_{\sigma^+0}=a_{\sigma^+}.
$$
By  Proposition \ref{basic-facts} i), $h_{n+1}$ is increasing on $B_{\sigma1}$ and decreasing on $B_{\sigma^+0}$, so
$$h_{n+1}([b_{\sigma1},a_{\sigma^+}])=h_{n+1}([b_{\sigma1},a_{\sigma^+0}])\subset [2,\infty).$$
Since $h_{n+1}(b_\sigma)<2$, we conclude that
\begin{equation}\label{sigma-1-+}
b_\sigma<b_{\sigma1}<a_{\sigma^+}.
\end{equation}
Together with $z_\sigma<a_{\sigma1}$,  we always have $[z_\sigma,b_\sigma]\prec B_{\sigma1}$.

ii) The proof is the same as i). In particular, we have that if $\sigma> 0^n$, then
\begin{equation}\label{sigma-0--}
b_{\sigma^-}<a_{\sigma0}<a_{\sigma}.
\end{equation}

iii) By Proposition \ref{gap-open-sigma-n}, $B_{\sigma0}<B_{\sigma1}$. If $n$ is odd and  $\sigma> 0^n$, then $B_{\sigma^-}<B_\sigma$ and $B_{\sigma0}\subset B_\sigma$. Hence $B_{\sigma^-}<B_{\sigma0}$.
  If $n$ is even and  $\sigma> 0^n$, then by \eqref{sigma-0--}, $B_{\sigma^-}<B_{\sigma0}$. If $n$ is odd and  $\sigma<1^n$, then by \eqref{sigma-1-+}, $B_{\sigma1}<B_{\sigma^+}$.  If $n$ is even and  $\sigma<1^n$, then
   $B_{\sigma}<B_{\sigma^+}$ and $B_{\sigma1}\subset B_\sigma$. Hence $B_{\sigma1}<B_{\sigma^+}$.
   \hfill $\Box$

\subsection{Proof of Lemma \ref{tech-2} and Corollary \ref{tech-cor}}\label{sec-pf-lem3.1-cor3.2}\

\noindent{\bf Proof of Lemma \ref{tech-2}. }\
i)  Applying Lemma \ref{tech-1} i) to $\sigma$ and  ii) to $\hat \sigma=\sigma0$, we have
$$a_{\sigma1}\notin \RR_n,\ \ B_{\sigma0}\subset B_\sigma \  \ \text{and }\ \  B_{\sigma01}\subset B_{\sigma0}.$$

Applying Lemma \ref{endsymbol} iii) and  Lemma \ref{tech-1} ii), iii) to $\hat \sigma=\sigma1$, we have
\begin{equation}\label{sig0}
a_{\sigma10},b_{\sigma10}\notin\RR_{n+1},\
B_{\sigma0}=B_{(\sigma1)^-}<B_{\sigma10}\prec B_{\sigma1}.
\end{equation}

If $a_\sigma\in \RR_{n-1}$, by Lemma \ref{endsymbol} ii), $a_{\sigma0}=a_\sigma\in \RR_{n-1}\subset\RR_n.$
Applying Lemma \ref{tech-1} ii) to $\hat \sigma=\sigma0$,
   we have
  \begin{equation}\label{sigma00-1}
  B_{\sigma00}\subset B_{\sigma0}\subset B_{\sigma}.
  \end{equation}

If $a_\sigma\notin\RR_{n-1}$, by Lemma \ref{endsymbol} iii),
 $a_{\sigma0}\notin \RR_n$ and hence  $a_{\sigma00}\notin\RR_{n+1}$.
Applying Lemma \ref{tech-1} ii) to $\hat \sigma=\sigma0$, we have
 $$b_{\sigma00}\notin\RR_{n+1},\ \ \ B_{\sigma00}\prec B_{\sigma0}.$$
  Now we show $a_\sigma<a_{\sigma00}$. By Lemma \ref{h-(n+1)} i), $h_{n+1}(a_\sigma)>2$. Hence
  $$
  h_{n+2}(a_\sigma)=h_{n+1}(a_\sigma)(h_n^2(a_\sigma)-2)-2>2.
  $$
  Since $a_{\sigma01}\in B_{\sigma01}
  \subset B_{\sigma0}\subset B_{\sigma}$, we have $a_{\sigma}\le a_{\sigma01}$. By Proposition \ref{basic-facts} i),
 $$h_{n+2}(a_{\sigma01})=-2.$$
 Since $a_{\sigma00}$ is the maximal solution of $h_{n+2}(E)=2$
 which are less than $a_{\sigma01}$, we conclude that $a_\sigma<a_{\sigma00}$. Since
 $b_{\sigma00}<b_{\sigma01}$ and by Lemma \ref{endsymbol} i),
 $b_{\sigma01}=b_{\sigma0}=z_\sigma$, we have
 \begin{equation}\label{sigma00-2}
 B_{\sigma00}=[a_{\sigma00},b_{\sigma00}]\subset (a_\sigma,z_\sigma)\subset {\rm int}( B_\sigma).
 \end{equation}

  Applying Lemma \ref{tech-1} ii) to $\hat \sigma=\sigma1$, we always have
 $$
 B_{\sigma11}\subset B_{\sigma1}.
 $$

  If $b_\sigma\in \RR_{n-1}$,  by Lemma
  \ref{tech-1} i),
 $
 B_{\sigma1}\subset B_{\sigma}.
 $
 Hence
 $$
 B_{\sigma11}\subset B_{\sigma1}\subset B_\sigma.
 $$
By \eqref{sig0}, $B_{\sigma0}<B_{\sigma10}\prec B_{\sigma1}$. Since $B_{\sigma0}, B_{\sigma1}\subset B_\sigma,$  we have
\begin{equation}\label{sigma10-1}
B_{\sigma10}=[a_{\sigma10},b_{\sigma10}]\subset (a_{\sigma0},b_{\sigma1})\subset {\rm int}(B_\sigma).
\end{equation}

If $b_\sigma\notin\RR_{n-1}$, we claim that
$$b_{\sigma10}<b_\sigma<b_{\sigma11}.$$

In fact, by Lemma \ref{endsymbol} iv), either $\sigma=1^n$, or  $a_{\sigma^+}\in \RR_{n-1}$.
 By Lemma \ref{h-(n+1)} i), $h_{n+1}(b_\sigma)<2$. Hence
  $$
  h_{n+2}(b_\sigma)=h_{n+1}(b_\sigma)(h_n^2(b_\sigma)-2)-2<2.
  $$
  If $\sigma=1^n$, then since $h_{n+2}(E)\to\infty$ as $E\to\infty$, there is a root of $h_{n+2}(E)=2$ in $(b_\sigma,\infty).$ Since $b_{\sigma11}=b_{1^{n+2}}$ is the largest root of $h_{n+2}(E)=2$, we conclude that $b_\sigma<b_{\sigma11}.$ Now assume $\sigma\ne1^n$, then $a_{\sigma^+}\in \RR_{n-1}$. By Lemma \ref{endsymbol} ii), $a_{\sigma^+0^t}=a_{\sigma^+}$. By Proposition \ref{basic-facts} i), $h_{n+2}$ is increasing on $B_{\sigma11}$ and decreasing on $B_{(\sigma11)^+}=B_{\sigma^+00}$. Hence by Floquet theory,
  $$
  h_{n+2}([b_{\sigma11},a_{\sigma^+}])=h_{n+2}([b_{\sigma11},a_{\sigma^+00}])\subset [2,\infty).
  $$
  Since $b_\sigma<a_{\sigma^+} $ and $h_{n+2}(b_\sigma)<2$, we conclude that $b_{\sigma}<b_{\sigma11}<a_{\sigma^+}$.
As a result, we always have $b_\sigma<b_{\sigma11}$.

On the other hand, we have
$$-2=h_{n+2}(b_{\sigma10})=h_{n+1}(b_{\sigma10})(h_n^2(b_{\sigma10})-2)-2.$$
So $h_{n+1}(b_{\sigma10})=0$ or $|h_n(b_{\sigma10})|= \sqrt{2}$.
By \eqref{sig0}, $b_{\sigma10}\notin\RR_{n+1}$. So $|h_n(b_{\sigma10})|=\sqrt{2}$.
If $\sigma=1^n$, then since $h_n(E)\ge 2$ for any $E\ge b_\sigma$, we have $b_{\sigma10}<b_\sigma.$
Now assume $\sigma\ne 1^n$.
By Floquet theory, $h_n([b_\sigma,a_{\sigma^+}])\cap (-2,2)=\emptyset.$
Since $b_{\sigma10}<b_{\sigma11}<a_{\sigma^+00}=a_{\sigma^+}$ and $h_{n}(b_{\sigma10})\in (-2,2)$,
we conclude that $b_{\sigma10}<b_\sigma.$ As a result, we always have $b_{\sigma10}<b_\sigma$. Hence the claim holds.

By \eqref{sig0}, $B_{\sigma0}<B_{\sigma10}$, so we have
$a_\sigma<z_\sigma=b_{\sigma0}<a_{\sigma10}<a_{\sigma11}.$
Together with the claim above, we have
$$B_{\sigma}\prec B_{\sigma11}.$$

By $B_{\sigma0}<B_{\sigma10}$ and the claim above, we also conclude that
\begin{equation}\label{sigma10-2}
B_{\sigma10}=[a_{\sigma10},b_{\sigma10}]\subset (b_{\sigma0}, b_\sigma)=(z_\sigma,b_\sigma)\subset {\rm int}(B_\sigma).
\end{equation}

By \eqref{sigma00-1} and \eqref{sigma00-2}, $B_{\sigma00}\subset B_\sigma$ always holds.

By \eqref{sigma10-1} and \eqref{sigma10-2}, $B_{\sigma10}\subset {\rm int}(B_\sigma)$ always holds.

ii) The proof is the same as i).
\hfill $\Box$

\noindent {\bf Proof of Corollary \ref{tech-cor}.}\
i) By Lemma \ref{basic-zeros} ii), ${\rm int}(B_\sigma)\cap \RR_{n-1}=\emptyset.$ So $B_\sigma\cap \RR_{n-1}\ne \emptyset$ if and only if $\partial B_\sigma\cap \RR_{n-1}\ne \emptyset$.
Write  $\hat\sigma:=\sigma|_{n-1}$.

At first assume $\partial B_\sigma\cap \RR_{n-1}\ne \emptyset$. If $a_\sigma\in \RR_{n-1}$, then by Lemma \ref{endsymbol} i), there exist $m<n$ even and  $\tau\in \Sigma_m$ such that $\sigma=\tau10^{n-m-1}$ and $z_\tau=a_{\tau10^{n-m-1}}=a_{\sigma}.$  If $m=n-1$, applying Lemma \ref{tech-1} ii) to $\tau$, we have $B_\sigma=B_{\tau1}\subset B_\tau=B_{\hat \sigma}$. If $m\le n-2$, then $a_{\hat\sigma}=z_\tau\in \RR_{n-2}$. By Lemma \ref{tech-1}, we always have
$B_{\sigma}=B_{\hat\sigma0}\subset B_{\hat\sigma}$. If $b_\sigma\in \RR_{n-1}$, the same proof shows that $B_{\sigma}\subset B_{\hat\sigma}.$

Next assume $B_\sigma\subset B_{\hat\sigma}$. If $n$ is even, by applying Lemma \ref{tech-1} i) to $\hat\sigma$, either $\sigma=\hat\sigma0$, hence $b_\sigma=z_{\hat\sigma}\in \RR_{n-1}$; or  $\sigma=\hat\sigma1$, hence $b_{\hat\sigma}\in \RR_{n-2}$. In the latter case, by Lemma \ref{endsymbol} ii), $b_{\sigma}=b_{\hat\sigma1}=b_{\hat\sigma}\in \RR_{n-1}$. Thus we always have $\partial B_\sigma\cap \RR_{n-1}\ne \emptyset.$ If $n$ is odd, the same proof shows that $\partial B_\sigma\cap \RR_{n-1}\ne \emptyset.$

ii) Write $\hat\sigma:=\sigma|_{n-1}$ and $\tilde \sigma:=\sigma|_{n-2}$. Assume $B_\sigma\not\subset B_{\hat\sigma}$. If $n=1$, by the convention, $B_{\hat\sigma}=B_\emptyset$ and $B_{\tilde\sigma}=\R$. So $B_\sigma\subset B_{\tilde\sigma}.$ Now assume $n\ge2$. At first, assume $n$ is odd. By applying Lemma \ref{tech-1} ii) to $\hat\sigma,$ the only possibility is that $a_{\hat\sigma}\notin \RR_{n-2}$ and $\sigma=\hat\sigma0$. So $\sigma=\tilde\sigma\sigma_{n-1}0$. Now by applying Lemma \ref{tech-2} i) to $\tilde\sigma$, we conclude that $B_{\sigma}=B_{\tilde\sigma\sigma_{n-1}0}\subset B_{\tilde\sigma}$. If $n$ is even, the same proof show that $B_{\sigma}\subset B_{\tilde\sigma}$. So ii) follows.

iii) At first assume  $\tau\in \Sigma_{n}$. By Proposition \ref{gap-open-sigma-n}, we must have $\sigma=\tau.$

  Next assume  $\tau\in \Sigma_{n+1}$. By Lemma \ref{tech-1} iii), if $\sigma\ne \tau|_n,$ then $B_\tau\cap B_\sigma=\emptyset.$ So we must have $\sigma=\tau|_n$.

Now assume  $\tau\in \Sigma_{n+2}$. By ii), $B_\tau\subset B_{\tau|_{n+1}}$ or $B_\tau\subset B_{\tau|_{n}}$. If $B_\tau\subset B_{\tau|_{n+1}}$, again by Lemma \ref{tech-1} iii), if $\sigma\ne \tau|_n,$ then $B_{\tau|_{n+1}}\cap B_\sigma=\emptyset.$ Consequently $B_\tau\cap B_\sigma=\emptyset$. So we must have $\sigma=\tau|_n$. If $B_\tau\subset B_{\tau|_{n}}$, by Proposition \ref{gap-open-sigma-n}, we must have $\sigma=\tau|_n$.

As a result, $\sigma=\tau|_n\lhd\tau.$
\hfill $\Box$

\bigskip

\noindent{\bf Acknowledgement}. 
 Liu was  supported by the National Natural Science Foundation of China, No. 11871098.  Qu was supported by the National Natural Science Foundation of China, No. 11790273 and No. 11871098. Yao was supported by the National Natural Science Foundation of China, No.11901311 and key technologies research and development program, 2020YFA0713300.


\section*{Appendix}

For the reader's convenience, we include  two tables of indexes  in this appendix. One is for the notations used in this paper, another is for the various orders used in this paper.

\smallskip


\begin{center}
\bf Table I: Index of notations
\end{center}
\begin{longtable}{ll}
\hline
$\xi, V_\xi$ & P-D sequence and potential, see Sec. \ref{sec-def}\\
\hline
$H_\lambda, \sigma_{\lambda}$ & PDH, the spectrum, see  \eqref{H-lambda}\\
\hline
$h_n(E)$ & trace polynomial, see \eqref{trace-poly-n} \\
\hline
$\mathscr O(E)$ & trace orbit of  $E$, see \eqref{O-E}\\
\hline
$\mathscr E_\infty$ & the set of $\infty$-type energies, see \eqref{E-lambda}\\
\hline
$\A,\mathbb G, A=[a_{\alpha \beta}]$ & alphabet, graph, adjacency matrix, see Sect. \ref{sec-main-coding} \\
\hline
$\Omega_A, \Omega_\infty$ & symbolic spaces, see Sect. \ref{sec-main-coding}\\
\hline
$\Omega_n,\Omega_\ast$& the set of admissible words, see \eqref{Omega-n-ast} \\
\hline
$\ZZ^o,\ZZ^e, \ZZ,\ZZ_n, \RR_n$ & zeros of $\{h_n:n\ge0\}$, see \eqref{Zero}, \eqref{Z-n-R-n}\\
\hline
$\E^o_l,\E^o_r,\E^e_l,\E^e_r,\widetilde {\mathcal E}_{l},\widetilde {\mathcal E}_{r},\mathcal F$ & eventually 2-periodic codings, see \eqref{gap-edge}\\
\hline
$\omega_\ast,\omega^\ast$ & minimum and maximum of $\Omega_\infty$, see \eqref{two-boundary}\\
\hline
$\GG, \GG_I,\GG_{II}, \GG_I^o, \GG_I^e$ & the set of  gaps of $\sigma_{\lambda}$, see Theorem \ref{main-gap-zero-inftyenergy} \\
\hline
$\G, \G_I,\G_{II}, \G_I^o, \G_I^e$ & the set of  gaps of $\Omega_\infty$, see Theorem \ref{symbolic-gap} \\
\hline
$\Sigma_n,\Sigma_\ast, \Sigma_\ast^o,\Sigma_\ast^e, \Sigma_\infty,\Sigma_\infty^{(2)}$ &  binary trees, see \eqref{Sigma-infty}, \eqref{Sigma-ast}, \eqref{Simga-ast-oe},  Sec. \ref{sec-code-B-n}\\
\hline
$\I,A(\I)$ & NS(SNS), limit set, see Sec. \ref{sec-NS}\\
\hline
$\sigma_{\lambda,n}$ & n-th approximation of $\sigma_{\lambda}$, see \eqref{sigman},\eqref{sigma-lamb-n}\\
\hline
$\B_n$ &  family of bands of level $n$, see \eqref{B-n}\\
\hline
$\BB_n$ & optimal covering of level $n$, see \eqref{def-BB-n}\\
\hline
$\widetilde A,\widetilde{\mathbb G},\widetilde{\BB}$ & sub-alphabet, graph, NS, see Sec. \ref{idea-low-dim} and \ref{sec-lower-dim}
\\
\hline
$\mathbb{M}_n,\varpi_n$ & see Sec. \ref{sec-code-B-n}, \eqref{varpi-n} \\
\hline
$\sigma^+,\sigma^-, \sigma|_k,\sigma\wedge\tau,\sigma\lhd\tau, |\sigma|$ & see Sec. \ref{sec-code-B-n} \\
\hline
$B_\sigma, a_\sigma, b_\sigma,z_\sigma$ & band, left and right endpoint, zero, see Sec. \ref{sec-code-B-n}\\
\hline
$I_w$ & rename of band, see Sec. \ref{Sec-code-BB-n}\\
\hline
$\pi: (\Omega_\infty,\preceq)\to (\sigma_\lambda,\le )$  & coding map, see Sec. \ref{sec-main-coding} and \ref{sec-coding-sigma-lambda} \\
\hline
$\Pi:(\Omega_\infty,\preceq)\to (\Sigma_\infty,\le)$ & another coding map, see Sec. \ref{statement-IDS} and \ref{sec-pro-Pi}\\
\hline
$\varepsilon:\Sigma_\infty\to [0,1]$ & evaluation map, see Sec. \ref{statement-IDS}\\
\hline
$N:\sigma_\lambda\to[0,1], N(E), N(G)$ \ \ \ \ & IDS,  see Sec. \ref{statement-IDS}\\
\hline
$\mathcal L, \Pi_\ast$ & see Sec. \ref{sec-pro-Pi}\\
\hline
$\ell, \ell^o,\ell^e$ & see Sec. \ref{sec-gap-symbolic}\\
\hline
\end{longtable}

%

\bigskip

In this paper, two types of orders are defined in various spaces: one is $\le$, another is $\preceq$. When $\le$ is used, it means standard or strong, depending on the context. When $\preceq$ is used, it means non standard or  weak.

\begin{center}
\bf Table II: Index for various orders
\end{center}

\begin{longtable}{ll}
\hline
$I<J, I\prec J$& two orders for bands, see \eqref{order-bands}\\
\hline
$\sigma\leq\tau$ & lexicographical order on  $\Sigma_n$ and  $\Sigma_{\infty}$, see Sec. \ref{sec-code-B-n}, \ref{sec-pro-Pi}\\
\hline
$\sigma\preceq\tau$ & total order on  $\Sigma_\ast$, see \eqref{order-Sigma-ast}\\
\hline
$\alpha<\beta,\alpha\prec\beta$ & two partial orders on $\A$, see \eqref{order}, \eqref{order-strong}\\
\hline
$\omega\le\omega', \omega\preceq \omega'$ & two  partial orders on $\Omega_\infty$, see Sec. \ref{sec-order-Omega}\\
\hline
$(x,y)\leq(x',y')$\ \ \ \ \ \ \  \ \  & partial order on $\mathbb{R}^2$, see Sec. \ref{sec-loc-pro-f}\\
\hline
\end{longtable}


\end{document}